\documentclass[a4paper,11pt]{article}
\usepackage[]{times}
\usepackage{fullpage}
\usepackage{graphicx}
\usepackage{float}
\usepackage{subfigure}
\usepackage{amsmath}
\usepackage{fancyhdr}
\usepackage[colorlinks]{hyperref}
\usepackage{xcolor}
\date{}

\newcommand{\prov}{{\sc Proof}.\hspace*{0mm} }

\newcommand{\QED}{$\rule{2mm}{2mm}$}

\newcommand{\natu}{{\sf I \! N}}
\newtheorem{theorem}{Theorem}[section]
\newtheorem{lemma}[theorem]{Lemma}
\newtheorem{e-proposition}[theorem]{Proposition}

\newtheorem{e-definition}[theorem]{Definition\rm}
\newtheorem{remark}{\it Remark\/}

\title{A variant of the Raviart-Thomas method \\ for smooth domains using straight-edged triangles}
\author{
    Fleurianne Bertrand$^{1}$ \& Vitoriano Ruas$^{2}$\thanks{Sorbonne Universit\'e, Campus Pierre et Marie Curie, 4 place jussieu, Couloir 55-65, 4\`eme \'etage, 75005 Paris, France.}
		\\[1mm]
	{\small $^{1}$ Fakult\"at f\"ur Mathematik, Technische Universit\"at Chemnitz, Germany.}\\[1mm]
  {\small $^{2}$ Institut Jean Le Rond d'Alembert, CNRS UMR 7190, Sorbonne Universit\'e, Paris, France.}\\[1mm]
  {\small e-mail: {\it vitoriano.ruas@upmc.fr}}}

\begin{document}
\maketitle

\begin{abstract}
Several physical problems modeled by second-order elliptic equations can be efficiently solved using mixed finite elements of the Raviart-Thomas family $RT_k$ for $N$-simplexes, introduced in the seventies. In case Neumann conditions are prescribed on a curvilinear boundary, the normal component of the flux variable should preferably not take up values at nodes shifted to the boundary of the approximating polytope in the corresponding normal direction. This is because the method's accuracy downgrades, which was shown in \cite{FBRT}. In that work an order-preserving technique was studied, based on a parametric version of these elements with curved simplexes. In this paper an alternative with straight-edged triangles for two-dimensional problems is proposed. The key point of this method is a Petrov-Galerkin formulation of the mixed problem, in which the test-flux space is a little different from the shape-flux space. 
After describing the underlying variant of $RT_k$, we show that it gives rise to uniformly stable methods taking the Poisson equation as a model problem. Then a priori error estimates of optimal order in the $[{\bf H}(div) \times L^2]$-norm are demonstrated. \\


\end{abstract}

\section{Introduction}

\hspace{4mm} This work deals with a particular type of Petrov-Galerkin formulation, primarily designed to preserve the order of finite element methods in the solution of boundary value problems posed in a two- or a three-dimensional domain having a smooth boundary of arbitrary shape, on which DOFs (degrees of freedom) are prescribed. In \cite{ZAMM} and \cite{IMAJNA} it was applied to conforming Lagrange and Hermite finite element methods with straight-edged triangles and tetrahedra to solve both second and fourth order elliptic equations, in the aim of improving their accuracy, while maintaining orders higher than one in the inherent energy norm. In \cite{JCAM} it was shown that this technique has the very same effect, as applied to nonconforming finite element methods, even those having boundary-prescribed DOFs different from function values at vertexes or any other point of an edge or a face of the elements. Although this approach can be adopted for any type of element geometry, it appears most effective as combined to simplexes. The technique under consideration is simple to implement in both two- and three-dimensional geometries, and optimal orders can be attained very naturally in both cases. Moreover, in the case of simplicial elements it definitively eliminates the need for curved elements, and consequently the use of non affine mappings. \\
\indent Here we further illustrate the application of this formulation, in the framework of the family of Raviart-Thomas mixed finite elements for simplexes to solve second order boundary value problems in a smooth plane domain, with Neumann conditions prescribed, at least on a boundary portion. We recall that in this case the normal components of the primal (flux) variable in the underlying mixed formulation are prescribed on such a boundary portion. As a matter of fact, with the current proposal we attempt to show once more that the approximate flux variable should preferably not take up corresponding prescribed values at nodes shifted to the boundary of the approximating polytope. Actually, even when there is no node shift, the prescribed flux on the boundary should be normal to the true boundary and not to the approximate boundary. This is because in the former case, a significant increase of accuracy occurs. At this point we should report that a technique based on this principle was studied in \cite{FBRT} and \cite{FBSIAM} for Raviart-Thomes mixed elements. It is based on a parametric version of those elements with curved simplexes. In contrast, in this paper the aforementioned Petrov-Galerkin formulation with straight-edged triangles is employed instead to tackle this problem. \\
\indent In this article we focus on the uniform stability and the convergence analysis of our method. With this aim we recall that, according to the well known theory set up by Babu\v{s}ka \cite{Babuska1,Babuska2} and Brezzi \cite{Brezzi}, a linear variational problem is well-posed only if the underlying continuous bilinear form is weakly coercive. In the case of problems defined upon a pair of spaces with the same finite dimension, such as finite-element spaces, this reduces to proving that the underlying bilinear form satisfies an \textit{inf-sup} condition for the pair of spaces under consideration (see e.g. \cite{COAM}). Here, in order to achieve this, we resort to the theoretical framework supplied in \cite{AsyMVF} for the case of asymmetric mixed formulations. By this denomination we mean mixed formulations in which the shape-product space is different from the test-product space.\\ 
\indent In order to bypass non essential difficulties, we study our methodology in connection with the Raviart-Thomas mixed finite elements, taking as a model the Poisson first order system in a bounded plane domain $\Omega$. \\

\indent An outline of the paper is as follows. In Section 2 we define the model problem to be treated by our mixed finite element method. General considerations are made on the kind of Petrov-Galerkin formulation employed in this work to approximate it and corresponding error bounds are given. This material is followed by the presentation in Section 3 of concepts, notations and properties related to the family of meshes considered in the sequel, more particularly in connection with the boundary of the curved domain.
Section 4 is devoted to the description of the specific spaces used here to approximate the model problem with mixed or  Neumann boundary conditions, namely, Raviart-Thomas spaces of arbitrary order based on straight-edged triangles; corresponding feasibility results are proven. In Section 5 the well-posedness and uniform stability analysis of the approximate problem is carried out. Section 6 deals with the interpolation errors in the $H(div)$-norm for the space of shape-fluxes under consideration. This leads to best possible a priori error estimates in the natural norm established in Section 7. In Section 8 we supply a relevant complement, make some important comments and provide a short summary of the article. \\
\indent We left aside lengthy proofs for a result used in Section 5, but for the sake of completeness they are given in Appendix I. Additionally we report in Appendix II some numerical experiments confirming formal results obtained in this work. \\

Throughout our study we shall use the following notations, most related to well known Sobolev function spaces (see e.g. \cite{Adams}), extended to the corresponding vector versions: First of all ${\mathcal A} \cdot {\mathcal B}$ represents the inner product of two tensors ${\mathcal A}$ and ${\mathcal B}$ of order $n$ for $n \in \natu$. $D$ being a bounded subset of $\Re^n$, $\parallel \cdot \parallel_{r,D}$ and $| \cdot |_{r,D}$ denote the standard norm and semi-norm of Sobolev space $H^{r}(D)$, for $r \in \Re^{+}$ with $H^0(D)=L^2(D)$. We further denote by $\parallel \cdot \parallel_{m,p,D}$ and $| \cdot |_{m,p,D}$ the usual norm and semi-norm of $W^{m,p}(D)$ for $m \in \natu^{*}$ and $p \in [1,\infty] \setminus \{2\}$ with $W^{0,p}(D)=L^p(D)$, and eventually for $W^{m,2}(D)=H^m(D)$ as well. In case $D=\Omega$ we omit this subscript. $P_k(D)$ represents the space of polynomials of degree less than or equal to $k$ defined in $D$. 
Finally we introduce the notation $D^j {\mathcal A}$ for the $j$-th (resp. $(j+1)-th$) order tensor, whose components are the $j$-th order partial derivatives with respect to the space variables of a function (resp. a vector field) ${\mathcal A}$ in the strong or in the weak sense.

\section{The model problem and its abstract Petrov-Galerkin approximation}

\hspace{4mm} In our reliability study we shall apply the abstract results given in \cite{AsyMVF} to a variant of the Raviart-Thomas family of mixed elements for triangles, specially designed to handle smooth two-dimensional domain. We take the Poisson equation with mixed or Neumann boundary conditions as a model, in the particular framework that we next describe. \\ 
\indent Let $\Omega$ be a two-dimensional domain assumed to be at least of the piecewise $C^2$-class and $\Gamma$ be its boundary, consisting of two portions $\Gamma_0$ and $\Gamma_1$ with $length(\Gamma_1 ) >0$. We denote by ${\bf n}$ the outer normal vector to $\Gamma$ and by $V$, either the space $L^2(\Omega)$ if $length(\Gamma_0) > 0$, or its subspace $L^2_0(\Omega)$ otherwise, where $L^2_0(D):=\{g|g \in L^2(D), \; \int_D g =0\}$ for a given domain D of $\Re^2$. \\
We consider the following model equation:\\
\begin{equation}
\label{Mixed}
\left\{
\begin{array}{l}
\mbox{Given } f \in V \mbox{ find } ({\bf p}; u) \mbox{ such that,}\\ 
 - \nabla \cdot {\bf p} = f \mbox{ in } \Omega; \\
 {\bf p} - \nabla u = {\bf 0} \mbox{ in } \Omega; \\
 u = 0 \mbox{ on } \Gamma_0 \mbox{ if } length(\Gamma_0) >0 \mbox{ or } \int_{\Omega} u = 0 \mbox{ otherwise}; \\
 {\bf p} \cdot {\bf n}= 0 \mbox{ on } \Gamma_1.
\end{array}
\right.
\end{equation}
Let ${\bf Q}$ be the subspace of ${\bf H}(div,\Omega)$ of those fields ${\bf q}$ such that 
${\bf q} \cdot {\bf n}=0$ on $\Gamma_1$ (cf. \cite{GiraultRaviart}). Denoting by $(\cdot,\cdot)$ the standard inner product of $L^2(\Omega)$ in scalar or vector version, problem (\ref{Mixed}) can be recast in the usual well-posed equivalent variational form \eqref{Dmix}-\eqref{definec} (cf. \cite{RaviartThomas}) that we next recall.\\
\begin{equation}
\label{Dmix}
\mbox{Find } ({\bf p}; u) \in {\bf Q} \times V \mbox{ such that } 
c(({\bf p};u),({\bf q};v)) = (f,v) \; \forall ({\bf q};v) \in {\bf Q} \times V, 
\end{equation}
where
\begin{equation}
\label{definec}
c(({\bf r};w),({\bf q};v)):=  ({\bf r},{\bf q})+( w, \nabla \cdot {\bf q})-(\nabla \cdot {\bf r}, v) \; \forall (({\bf r};w);({\bf q};v)) \in ({\bf Q} \times V) \times ({\bf Q} \times V)
\end{equation} 
In an attempt to handle the boundary conditions for ${\bf p}$ more accurately, in this work we will employ a kind of Petrov-Galerkin finite-element formulation to approximate problem \eqref{Dmix}-\eqref{definec}, which we next describe.\\
\indent To begin with, assume that $\Omega$ is approximated by a polygon $\Omega_h$ defined as the union of triangles in a mesh with maximum edge length equal to $h$. Let us denote by $(\cdot,\cdot)_h$ the standard inner product of $L^2(\Omega_h)$. While $\parallel \cdot \parallel_{0,p,h}$ represents the standard norm of $L^p(\Omega_h)$ with $p \in [1,+\infty] \setminus \{2\}$, let $\| \cdot \|_{0,h}$ stand for the norm of $L^2(\Omega_h)$ associated with $(\cdot,\cdot)_h$. $[| \cdot |]_h$ being the standard norm of ${\bf H}(div,\Omega_h)$ we denote by $||| \cdot |||_h$ the norm of ${\bf H}(div,\Omega_h) \times L^2(\Omega_h)$ given by $||| ({\bf q};v) |||_h := [|{\bf q}|]_h + \| v \|_{0,h}$ $\forall ({\bf q};v) \in {\bf H}(div,\Omega_h) \times L^2(\Omega_h)$.\\
\noindent We consider a pair of spaces $({\bf P}_h;{\bf Q}_h) \subset [{\bf H}(div;\Omega_h)]^2$ together with a subspace $V_h$ of $L^2(\Omega_h)$ and also of $L^2_0(\Omega_h)$ if $length(\Gamma_0)=0$, all of them being finite-dimensional. \\
In our approximation of \eqref{Dmix}-\eqref{definec} the bilinear form $c$ will be replaced  
by $c_h$ given by 
\begin{equation}
\label{bilinearformch}
c_h(({\bf r}; w),({\bf q}; v)):= a_h({\bf r},{\bf q}) + b_h( w, {\bf q}) + d_h({\bf r}, v) 
\forall ({\bf r}; w) \in {\bf P}_{h} \times V_{h} \mbox{ and } \forall ({\bf q}; v) \in {\bf Q}_{h} \times V_{h} 
\end{equation}
where
\begin{equation}
\label{ahbhdh}
a_h({\bf r},{\bf q}):=({\bf r},{\bf q})_h; \; b_h(w,{\bf q}):=(w, \nabla \cdot {\bf q})_h; \; d_h({\bf r},v):= - (\nabla \cdot {\bf r},v)_h.
\end{equation}
Notice that $c_h$ is continuous on $({\bf P}_{h} \times V_{h}) \times ({\bf Q}_{h} \times V_{h})$, for there is a norm $\| c_h \|$  of $c_h$ such that,
\begin{equation}
\label{normc}
c_h(({\bf r}; w),({\bf q}; v)) \leq \| c_h \| |||({\bf r}; w) |||_h |||({\bf q}; v) |||_h \; \forall (({\bf r}; w),({\bf q}; v)) \in  ({\bf P}_{h} \times V_{h}) \times  ({\bf Q}_{h} \times V_{h}),
\end{equation}
Actually we can take $\|c_h\| = \max[\|a_h\|,\|b_h\|,\|d_h\|]$ and observing that $\|a_h\|=\|b_h\|=\|d_h\|=1$, we infer that $c_h$ is uniformly continuous, that is, its norm is independent of the particular subspaces ${\bf P}_h$, ${\bf Q}_h$, $V_h$ under consideration.\\

Now letting $F_h(v)$ be a suitable approximation of $(f,v)$ for any $v \in V_h$, the approximate problem to solve is 
\begin{equation}
\label{Dmixh}
\left\{
\begin{array}{l}
\mbox{Find } ({\bf p}_h; u_h) \in {\bf P}_{h} \times V_{h} \mbox{ such that,}\\
c_h(({\bf p}_h;u_h),({\bf q};v)) = F_h(v) \; \forall ({\bf q};v) \in {\bf Q}_h \times V_h. 
\end{array}
\right.
\end{equation}
\eqref{Dmixh} is a standard linear variational problem defined upon a pair of different Hilbert spaces for the solution and the testing elements. Referring to \cite{Brezzi}, provided the right hand side functional is continuous, problems of this kind have a unique solution if and only if their bilinear form is continuous and weakly coercive over this pair of spaces. According to \cite{COAM}, in the case under study, weak coercivity means that
\begin{itemize} 
\item There exists a constant $\gamma_h >0$ such that $\forall ({\bf r};w) \in {\bf P}_h \times V_h$ there is a non zero $({\bf q};v) \in {\bf Q}_h \times V_h$ such that $c_h(({\bf r};w),({\bf q};v)) \geq \gamma_h ||| ({\bf r};w) |||_h ||| ({\bf q};v) |||_h$;
\item
\vspace{-1mm} 
The dimensions of ${\bf P}_h$ and ${\bf Q}_h$ are the same.
\end{itemize}
Now assume that there exists a function $\tilde{u}$ defined in a domain $\tilde{\Omega}$ large enough to contain both $\Omega$ and $\Omega_h$ for all meshes in use, which coincides with $u$ in $\Omega$ and fulfills $\tilde{u} \in H^1(\tilde{\Omega})$, $\Delta \tilde{u} \in L^2(\tilde{\Omega})$. The field $\tilde{\bf p}:=\nabla \tilde{u}$ belongs to ${\bf H}(div,\tilde{\Omega})$.\\
Next, setting $\tilde{\Omega}_h:= \Omega_h \cup \Omega \subset \tilde{\Omega}$, we observe that $c_h(({\bf r};w),({\bf q};v))$ is well defined $\forall ({\bf r};w)\in {\bf H}(div,\tilde{\Omega}_h) \times L^2(\tilde{\Omega}_h)$ and $\forall ({\bf q};v) \in {\bf Q}_h \times V_h$ and thus $||| ({\bf r};w) |||_h$ as well. Notice that $c_h$ is also continuous in the sense $c_h(({\bf r};w),({\bf q};v)) \leq \|c_h\| ||| ({\bf r};w) |||_h ||| ({\bf q};v) |||_h$ $ \forall 
(({\bf r};w);({\bf q};v)) \in [{\bf H}(div,\tilde{\Omega}_h) \times L^2(\tilde{\Omega}_h)] \times [{\bf Q}_h \times V_h]$ with $\|c_h\|=1$. \\
On the other hand, in the case of curvilinear boundaries, in principle, we are dealing with an external approximation also called a variational crime (cf. \cite{StrangFix}). Hence we resort to the general upper bound for linear variational problems given in \cite{COAM} for weakly coercive bilinear forms defined upon two different spaces. In the case under study this writes
\begin{equation}
\label{unifymajor}
\left\{
\begin{array}{l}
||| (\tilde{\bf p};\tilde{u}) - ({\bf p}_h;u_h) |||_h \leq \displaystyle \frac{1}{\gamma}_h\left[ \inf_{({\bf r};w) \in {\bf P}_{h} \times V_{h}}  ||| (\tilde{\bf p};\tilde{u}) - ({\bf r};w) |||_h  \right. \\
\left. +  \displaystyle \sup_{({\bf q};v) \in {\bf Q}_{h} \times V_{h}/|||({\bf q};v)|||_h=1}| c_h((\tilde{\bf p};\tilde{u}),({\bf q};v))- F_h(v) |. 
\right]
\end{array}
\right.
\end{equation}
From \eqref{unifymajor} the constant $\gamma_h$ is seen to play a key role in error estimates for problem \eqref{Dmixh}. 
Incidentally, the bilinear form $c_h$ is of the type studied in \cite{AsyMVF} with $P={\bf P}_{h}$, $Q={\bf Q}_{h}$ $U \equiv V = V_{h}$. Actually, as recalled in the introductory paragraph of Section 5, the key to establishing its weak coercivity is Theorem 3.1 of that article, whose assumptions are \textit{inf-sup} conditions to be fulfilled by the bilinear forms $a_h$, $b_h$ and $d_h$ addressed in Subsections 5.1, 5.2 and 5.3, respectively. Moreover it is proven in \cite{AsyMVF} that $\gamma_h$ can be expressed in terms of the constants of such conditions that supposedly hold for those bilinear forms. Naturally enough, the \textit{inf-sup} conditions are  established for the three of them in the aforementioned subsections, in connection with the particular form of \eqref{Dmixh} considered here. \\
More specifically, the spaces ${\bf Q}_h$, ${\bf P}_h$ and $V_h$ used in this work are linked to the Raviart-Thomas mixed finite element method of order $k$ for triangles \cite{RaviartThomas}. We recall that (see e.g. \cite{BrezziFortin}), if $\Omega$ is a polygon, provided $u$ belongs to $H^{k+2}(\Omega)$ (cf. \cite{Adams}), this method is of order $k+1$ in the standard norm of the space ${\bf H}(div,\Omega) \times L^2(\Omega)$. Notice however that, in general, if $\Omega$ is a polygon, the above regularity can be expected at most if $k=0$ and $\Omega$ is convex. As a matter of fact, the purpose of problem \eqref{Dmixh} is to take the best advantage of higher order methods in this family, in the framework of problems posed in smooth domains, since in this case the exact solution can be much more regular. 

\section{Meshes and related sets, notations and assumptions} 

\hspace{4mm} The finite element spaces will shall work with are defined in connection with triangulations of $\Omega$, whose main features are next addressed.\\
\indent Let ${\mathcal F} := \{{\mathcal T}_h\}_h$ be a family of meshes consisting of straight-edged triangles, satisfying the well known compatibility conditions for the finite element method (cf. \cite{Ciarlet}), $h$ being the maximum edge length of all the triangles in the mesh ${\mathcal T}_h$. Every element of this mesh is to be viewed as a closed set and ${\mathcal T}_h$ is assumed to fit $\Omega$ in such a way that all the vertexes of the polygon $\Omega_h$ lie on $\Gamma$, where $\Omega_h$ is the interior of $\cup_{T \in {\mathcal T}_h} T$. 
 We recall that $\tilde{\Omega}_h = \Omega \cup \Omega_h$ and further define $\Omega^{'}_h := \Omega \cap \Omega_h$. 
The norms of $L^2(\Omega^{'}_h)$, ${\bf H}(div,\Omega^{'}_h)$ and ${\bf H}(div,\Omega^{'}_h) \times L^2(\Omega^{'}_h)$ are respectively denoted by $\| \cdot \|_{0,h}^{'}$, $[| \cdot |]_{h}^{'}$ and $||| \cdot |||_{h}^{'}$. \\
We assume that the transition points between $\Gamma_1$ and $\Gamma_0$, if any, are always vertexes of the polygon $\Omega_h$.
The boundary of $\Omega_h$ 
is denoted by $\Gamma_h$ and 
${\bf n}_h$ represents the unit outer normal vector to $\Gamma_h$. 
Finally we define $\Gamma_{1,h}$ to be the union of the edges of $\Gamma_h$ having both ends in the union of $\Gamma_1$ and the set of transition points between $\Gamma_0$ and $\Gamma_1$, if any, and let $\Gamma_{0,h}$ be the closure of  $\Gamma_h \setminus \Gamma_{1,h}$.\\
${\mathcal F}$ is assumed to be a regular family of partitions  
in the sense of \cite{Ciarlet} (cf. Section 3.1), though not necessarily quasi-uniform. The maximum edge length of every $T \in {\mathcal T}_h$ is denoted by $h_T$. Throughout this work we make the very reasonable assumption that no triangle $T \in {\mathcal T}_h$ has its three vertexes on $\Gamma$. It follows that every element in ${\mathcal T}_h$ has at most one edge contained in $\Gamma_h$.\\
Let ${\mathcal S}_h$ (resp. ${\mathcal S}_{1,h}$) be the subset of ${\mathcal T}_h$ consisting of triangles $T$ having 
one edge, say $e_T$, on $\Gamma_{h}$ (resp. $\Gamma_{1,h}$). We denote by $O_T$ the vertex of $T \in {\mathcal S}_h$ not belonging to $\Gamma$ and by ${\bf n}_T$ the unit outer normal vector along $e_T$. Notice that the interior of every triangle in ${\mathcal T}_h \setminus {\mathcal S}_h$ has an empty intersection with $\Gamma_{h}$. \\
The following additional considerations are useful: $\forall T \in {\mathcal S}_h$ we denote by $\Delta_T$ the closed subset of $\tilde{\Omega}_h$ not containing $T$ delimited by $\Gamma$ and the edge $e_T$ of $T$, as illustrated in Figure 1. Moreover, $\forall T \in {\mathcal S}_h$, we define $T^{\Gamma}$ to be the set $T \cup \Delta_T$ and $T^{\Gamma^{'}}$ to be the closure of the set $T \setminus \Delta_T$. Notice that, if $\Omega$ is convex $T^{\Gamma^{'}}=T$ and $\Omega_h$ is a proper subset of $\Omega$. Otherwise there are sets $T^{\Gamma}$ which coincide with $T$. Finally we set $\Gamma_T := T^{\Gamma} \cap \Gamma$ for any $T \in {\mathcal S}_h$. \\
\indent Although this is not a real restriction to the practical application of our method, for theoretical purposes and for the sake of simplicity, throughout this article we assume that every triangulation ${\mathcal T}_h \in {\mathcal F}$ is constructed in such a way that all the transition points between the convex and concave portions of $\Gamma$, if any, are vertexes of triangles thereof. In this respect, it is useful to define two subsets ${\mathcal S}_{0,h}^{'}$ and ${\mathcal S}_{1,h}^{'}$ of ${\mathcal S}_{0,h}$ and ${\mathcal S}_{1,h}$, which contain triangles having two vertexes on concave portions of $\bar{\Gamma}_0$ and $\bar{\Gamma}_1$, respectively. Notice that, owing to the above construction, such transition points are the same as those between convex and concave portions of $\Gamma_h$. We also set ${\mathcal S}_h^{'}:={\mathcal S}_{0,h}^{'} \cup {\mathcal S}_{1,h}^{'}$. 
\begin{figure}[h]
\label{fig1}
\centerline{\includegraphics[width=3.8in]{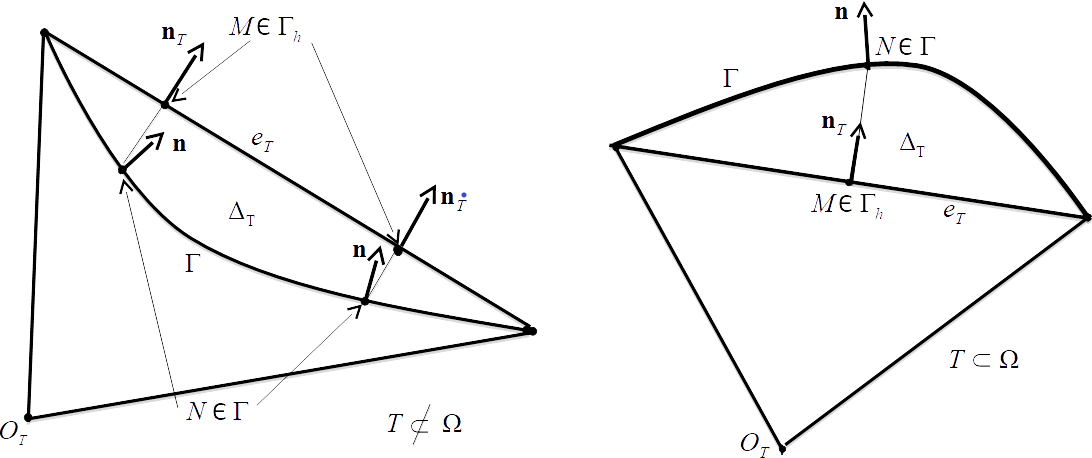}}
\vspace*{8pt}
\caption{Normal-flux DOFs shifted from $\Gamma_{1,h}$ to $\Gamma_1$ for $RT_1(T)$ (left) and $RT_0(T)$ (right) with $T \in {\mathcal S}_{1,h}$}
\end{figure} 
 
\hspace{4mm} For convenience, the reliability analysis of our method will be carried out under the following assumption, which is by no means necessary, neither to define it, nor to implement it:\\
 
\noindent \underline{\textit{Assumption}$^{*}$ :} $h$ is small enough for the intersection $N \in T^{\Gamma}$ with $\Gamma$ of the perpendicular to $e_T$ passing through any point $M \in e_T$ to be unique, $\forall T \in {\mathcal S}_h$ \rule{2mm}{2mm} \\

\noindent The main idea behind \textit{Assumption}$^{*}$ is to ensure that $\Gamma_h$ follows very closely the shape of $\Gamma$, which might not be the case of a coarse mesh.\\

To close this section we prove

\begin{e-proposition}
\label{prop01}
If \textit{Assumption}$^{*}$ holds, there exists a constant $C_{0}$ depending only on $\Gamma$ such that 
$\forall M \in e_T$ the length of the segment perpendicular to $e_T$ joining $M$ and $N \in \Gamma_T$ is bounded above by $C_{0} h_T^2$.
\rule{2mm}{2mm}
\end{e-proposition}
For the proof of Propositions \ref{prop01} we refer to Proposition 2.1 of \cite{ZAMM} or yet to \cite{LewNegri}, Lemma A.2.

\begin{e-proposition}
\label{prop02}
If $\Gamma$ is of the piecewise $C^{2}-class$ and \textit{Assumption}$^{*}$ holds, $\forall T \in {\mathcal S}_{1,h}$ and $\forall N \in \Gamma_T$ the outer normal ${\bf n}(N)$ to $\Gamma_T$ at $N$ satisfies $|{\bf n}(N) - {\bf n}_T| \leq C_1 h_T$, where $C_{1}$ is a constant depending only on $\Gamma_1$.
\end{e-proposition} 
 \prov First we recall Proposition \ref{prop02} of \cite{ZAMM}, according to which there there is a constant $C_{\Gamma}$ depending only on $\Gamma$ such that
\begin{equation}
\label{CGamma}
|\tan{\varphi}| \leq C_{\Gamma} h_T, \mbox{ where } \varphi \mbox{ is the angle between }{\bf n}(N) \mbox{ and }{\bf n}_T.
\end{equation} 
On the other hand we have $|{\bf n}(N) - {\bf n}_T|^2 = 2(1-{\bf n}(N) \cdot {\bf n}_T)= 2(1-cos \varphi)$. Thus 
$|{\bf n}(N) - {\bf n}_T|^2 = 2 \sin^2{\varphi} /(1+\cos{\varphi}) = 2 \tan^2{\varphi}/(\cos^{-2}{\varphi} + \cos^{-1}{\varphi})$, that is,  
$|{\bf n}(N) - {\bf n}_T|^2 = 2 \tan^2{\varphi}/(\tan^2{\varphi} + 1 + \sqrt{1+ \tan^2{\varphi}}) \leq \tan^2{\varphi}$. It trivially follows that 
$|{\bf n}(N) - {\bf n}_T| \leq  C_{\Gamma} h_T$, i.e., $C_1=C_{\Gamma}$. \rule{2mm}{2mm}

\section{The finite-element approximation method}

\hspace{4mm} In this section we give the specific form of the finite element approximation spaces $V_h$, ${\bf Q}_h$ and ${\bf P}_h$ in \eqref{Dmixh}, associated with ${\mathcal T}_h$ and an integer $k \geq 0$, that we advocate for the numerical solution of \eqref{Mixed}.\\
\indent First of all we denote by $RT_k(T)$ the Raviart-Thomas space of vector fields defined in a triangle $T$ (cf. \cite{RaviartThomas}), whose dimension is $(k+1)(k+3)$.\\
\indent $\bullet$ $V_{h}$ is the space of functions $v$ such that $\forall T \in {\mathcal T}_h$ 
$v_{|T} \in P_k(T)$, satisfying $v \in L^2_0(\Omega_h)$ if $length(\Gamma_0)=0$.  \\
\indent $\bullet$ ${\bf Q}_{h}$ is the space consisting of fields ${\bf q} \in {\bf H}(div,\Omega_h)$, whose restriction to each triangle $T \in {\mathcal T}_h$ belongs to $RT_k(T)$, and whose normal component vanishes on $\Gamma_{1,h}$.    
We recall that, classically, the DOFs of $RT_k(T)$ are the $k+1$ moments of order $k$ of the field's normal components along each one of the edges of $T$ and, for $k>0$ only, the $k(k+1)$ moments of order $k-1$ in $T$ of both field's components. 
For convenience we extend every field ${\bf q} \in {\bf Q}_{h}$ by zero in $\bar{\Omega} \setminus \bar{\Omega}_h$. \\
\indent $\bullet$ ${\bf P}_{h}$ in turn is the space of fields ${\bf p}$ defined in $\tilde{\Omega}_h$ having the properties listed below.
\begin{enumerate} 
\item The restriction of ${\bf p}$ to every $T \in {\mathcal T}_h$ belongs to the Raviart-Thomas space $RT_k(T)$;
\item The normal component of every ${\bf p}$ is continuous across all the edges of the elements in ${\mathcal T}_h$ contained in the interior of $\Omega_h$; 
\item In $\bar{\Omega} \setminus \bar{\Omega}_h$ a field ${\bf p}$ is expressed in such a way that $\forall T \in {\mathcal S}_h$ its polynomial expression in $T$ extends to $\Delta_T$ in case $\Delta_T \setminus T$ is not empty;
\item $\forall T \in {\mathcal S}_{1,h}$, the (outer) normal component of ${\bf p}$ vanishes at every point which is the nearest intersection with $\Gamma_1$ of the perpendicular to the edge $e_T$ passing through one of the $k+1$ Gauss-quadrature points of $e_T$. The intersection point on $\Gamma_{T}$ and the corresponding Gauss-quadrature point are generically denoted by $N$ and $M$, respectively, as illustrated in Figure 1.
\end{enumerate}

Before establishing that ${\bf P}_h$ is a non empty finite-dimensional space we recall the following lemma, whose proof is given in 
\cite{ZAMM}, except for the obvious left inequality in \eqref{L2TDelta}:
\begin{lemma}
\label{TDelta}
(\cite{ZAMM}) Provided $h$ satisfies \textit{Assumption}$^{*}$, there exist three constants ${\mathcal C}_{\infty}$, ${\mathcal C}_0^{-}$ and ${\mathcal C}_0^{+}$ depending only on $\Gamma$ and the shape regularity of ${\mathcal T}_h$ (cf. \cite{BrennerScott}, Ch.4, Sect. 4) such that $\forall w \in P_k(T^{\Gamma})$ and $\forall T \in {\mathcal S}_h$ it holds:  
\begin{equation}
\label{LinftyTDelta}
\parallel  w \parallel_{0,\infty,T^{\Gamma}} \leq {\mathcal C}_{\infty} \parallel  w \parallel_{0,\infty,T^{\Gamma^{'}}} 
\end{equation}
and
\begin{equation}
\label{L2TDelta}
{\mathcal C}_0^{-} h_T^{-1} \parallel  w \parallel_{0,T^{\Gamma}} \leq \parallel  w \parallel_{0,\infty,T^{\Gamma}} \leq {\mathcal C}_0^{+} h_T^{-1} \parallel  w \parallel_{0,T^{\Gamma^{'}}} \; \mbox{\QED}.
\end{equation}
\end{lemma}

Now we state
\begin{e-proposition}
\label{prop1}  
Provided $h$ is small enough, given an arbitrary set $B_k$ of $m_k$ real values with $m_k=(k+1)(k+2)$, $\forall T \in {\mathcal S}_{1,h}$ there exists a unique field in $RT_k(T)$, whose normal components across $\Gamma_1$ vanish at the $k+1$ points $N \in \Gamma_T$ defined in accordance with item 4. of the above definition of ${\bf P}_h$, and whose remaining $m_k$ DOFs of fields in the space $RT_k(T)$ take up their values in the set $B_k$ on a one-to-one basis. 
\end{e-proposition}   
 
\prov Let $B_k=\{b_1,b_2,\ldots,b_{m_k}\}$ and $\vec{b}$ be the vector of $\Re^{n_k}$ given by $\vec{b}:=[0,\ldots,0,b_1,\ldots,b_{m_k}]$, i.e. $\vec{b}$ is a vector whose $k+1$ first entries are zero. If the boundary nodes $N$ were replaced by the corresponding $M \in \Gamma_{1,h} \cap T$, it is clear that the result would hold true, according to the well known properties of the space $RT_k(T)$. The vector $\vec{a}$ of coefficients $a_i$ for $i=1,2,\ldots,n_k$ of the canonical basis fields 
${\bf q}_{i,T}$ of $RT_k(T)$ for $1 \leq i \leq n_k$ would be given by $a_i=b_i$ for $1 \leq i \leq n_k$. \\
Let us denote by ${\mathcal H}_i$ the DOFs of $RT_k(T)$ and by ${\bf q}_{i,T}$ the underlying canonical basis fields, for $i=1,2,\ldots,n_k$, in such a way that the first $k+1$ ones correspond to normal components across $e_T$. In doing so, the matrix $E$ whose entries are $e_{ij} := {\mathcal H}_i({\bf q}_{j,T})$ is the identity matrix. Now we define $\tilde{\mathcal H}_i$ to be ${\mathcal H}_i$ if $k+1 < i \leq n_k$ and otherwise to be the normal component across $\Gamma_1$ at one of the $k+1$ points $N$ constructed in accordance with item 4. of the definition of ${\bf P}_h$. For convenience we further number these points $N$ and the corresponding point $M$ on $e_T$, respectively by $N_i$ and $M_i$, for $i=1,\ldots,k+1$. We also denote by ${\bf n}_i$ the outer normal to $\Gamma_1$ at $N_i$, for $i=1,2,\ldots,k+1$. \\
The Proposition will be proved if the $n_k \times n_k$ linear system $\tilde{E} \vec{a} = \vec{b}$ is uniquely solvable, where $\tilde{E}$ is the matrix with entries $\tilde{e}_{ij}:=\tilde{\mathcal H}_i({\bf q}_{j,T})$. Clearly we have $\tilde{E} = E + D_E$, where the entries of $D_E$ are $d_{ij}:= \tilde{\mathcal H}_i({\bf q}_{j,T}) - {\mathcal H}_i({\bf q}_{j,T})$. \\
While $d_{ij} = 0$ $\forall j$ if $k+1 < i \leq n_k$, in case $1 \leq i \leq k+1$ we have 
$$d_{i,j} = {\bf q}_{j,T}(N_i) \cdot {\bf n}_i-{\bf q}_{j,T}(M_i) \cdot {\bf n}_T.$$
This means that,
\begin{equation}
\label{drs}
\left\{
\begin{array}{l}
d_{ij} = r_{ij} + s_{ij}, \mbox{ where} \\
r_{ij} = {\bf q}_{j,T}(N_i) \cdot {\bf n}_i - {\bf q}_{j,T}(N_i) \cdot {\bf n}_T \\
s_{ij} = {\bf q}_{j,T}(N_i) \cdot {\bf n}_T - {\bf q}_{j,T}(M_i) \cdot {\bf n}_T.
\end{array}
\right.   
\end{equation}
From Proposition \ref{prop02} we may write $| {\bf n}_i - {\bf n}_T | \leq C_{\Gamma} h_T$ for $i=1,2,\ldots,k+1$. Taking into account { the fact that there is a mesh-independent constant $\hat{C}_{\bf q}$ such that $\|{\bf q}_{j,T}\|_{0,\infty,T} \leq \hat{C}_{\bf q}$ for $i=1,2,\ldots,k+1$, it follows that,
\begin{equation}
\label{rij}
|r_{ij}| \leq C_R h_T\; \forall \; i,j,\mbox{ with } C_R= \hat{C}_{\bf q} C_{\Gamma}.
\end{equation}
On the other hand, from Proposition \ref{prop01}, the length of the segment $\overline{M_iN_i}$ is bounded 
above by $C_{0} h_T^2$. 
Thus, a simple Taylor expansion of the function ${\bf q}_{j,T} \cdot {\bf n}_T$ implies that,
$$|s_{ij}| \leq C_{0} h_T ^2 \max_{{\bf x} \in T^{\Gamma}} | \nabla\{{\bf q}_{j,T} \cdot {\bf n}_T\}({\bf x})|,$$ 
or yet, recalling \eqref{LinftyTDelta},  
\begin{equation}
\label{boundeij}
|s_{ij}| \leq C_{0} {\mathcal C}_{\infty} h_T^2 \max_{{\bf x} \in T} |\nabla \{{\bf q}_{j,T} \cdot {\bf n}_T\}({\bf x})|.
\end{equation}
Finally, from Lemma \ref{TDelta} and a standard argument we know that the maximum in \eqref{boundeij} is  
bounded above by a mesh-independent constant times $h_T^{-1}$. In short we have 
\begin{equation}
\label{sij}
|s_{ij}| \leq C_S h_T \; \forall \; i,j,
\end{equation} 
where $C_S$ is a constant independent of $T$. \\
It follows from \eqref{drs}, \eqref{rij} and \eqref{sij} that the matrix $\tilde{E}$ equals the identity matrix plus an $O(h_T)$ matrix 
$D_E =[d_{ij}]$. Therefore $\tilde{E}$ is an invertible matrix, as long as $h$ is sufficiently
 small, more precisely, if $h < [\hat{C}_k]^{-1}/2$, with $\hat{C}_k = (C_R + C_S)/ \Theta$ for a suitably small $\Theta$. \QED \\

Thanks to Proposition \ref{prop1} it is legitimate to  
pose problem \eqref{Dmixh} for the spaces ${\bf P}_{h}$, ${\bf Q}_{h}$ and $V_h$ defined in Subsection 3.1, whose solution $({\bf p}_h;u_h)$ approximates the solution of (\ref{Mixed}) with $dim {\bf P}_h = dim {\bf Q}_h$.\\

We close this section by proving a property of space ${\bf P}_{h}$ to be repeatedly used in the sequel.
\begin{lemma}
\label{GammaT}
Assume that $\Gamma$ is of the $C^{2k+2}$-class and of the piecewise $C^{2k+3}$-class. \\
Given ${\bf p}$ in ${\bf P} _{h}$ let $g$ be its normal trace on $\Gamma_1$. Let $v$ be given in $V_{h}$. There exists a constant ${\mathcal C}_{\Gamma_1}$ such that $\forall T \in {\mathcal S}_{1,h}$ we have
\begin{equation}
\label{errgv}  
\displaystyle \left|\int_{\Gamma_T} g \; v \; ds \right| \leq {\mathcal C}_{\Gamma_1} h_T^{k+1} \| {\bf p} \|_{k+1,T} \| v \|_{0,T}.
\end{equation} 
\end{lemma}
\prov
First of all we denote by $\phi$ the function of $x \in [0,l_T]$ representing the $y$-coordinate of a point $N \in \Gamma_T$ corresponding to $M \in e_T$ with abscissa $x$ in the local direct cartesian coordinate system $(A_T,x,y)$, whose origin $A_T$ is chosen in such a manner that $y(x) \geq 0$ $\forall x$. In doing so, for any bounded function $a$ defined in a neighborhood of $\Gamma_T$, we define $\breve{a}$ in $[0,l_T]$ by $\breve{a}(x)=a[x,\phi(x)]$. It is clear that there exists a mesh-independent constant $C_q$ such that for every such functions $a$ we have
\begin{equation}
\label{errgv2}  
\displaystyle \left| \int_{\Gamma_T} a \; ds \right| \leq C_q \displaystyle \left| \int_{0}^{l_T} \breve{a} \; dx \right|. 
\end{equation}
According to our regularity assumption on $\Gamma$ the derivative of order $2k+3$ of $\phi$ and the derivative of order $2k+2$ of the  normal vector ${\bf n}$ are uniformly bounded in $(0,l_T)$. 
Hence $\breve{g}$ belongs to $W^{2k+2,\infty}(0,l_T)$.\\
Now we consider $a$ to be the function $gv$. Noticing that, by construction, $\breve{g}$ vanishes at the $k+1$ Gauss-quadrature points of $e_T$, so does $\breve{a}$. Thus by well known properties of Gauss integration in one-dimensional space (see e.g. \cite{Blum}, \cite{Quarteroni}), it holds for a constant $\breve{C}$ independent of $T$, 
\begin{equation}
\label{errgv3}  
\displaystyle \left| \int_0^{l_T} \breve{a} \; dx \right| \leq \breve{C} l_T^{2k+3} \displaystyle \sup_{x \in (0,l_T)} \left| \frac{d^{2k+2}\breve{a}}{dx^{2k+2}}\left( x \right) \right|
\end{equation}
Combining \eqref{errgv2} and \eqref{errgv3}, setting $C_{a}:= C_q \breve{C}$ we easily come up with
\begin{equation}
\label{errgv4}  
\displaystyle \left| \int_{\Gamma_T} a \; ds \right| \leq C_{a} h_T^{2k+3} | \breve{a} |_{2k+2,\infty,e_T}.
\end{equation}
Since $\breve{a}=\breve{g} \breve{v}$, after straightforward manipulations we obtain for another mesh-independent constant $C_a^{,}$
\begin{equation}
\label{errgv5}  
\displaystyle \left| \int_{\Gamma_T} a \; ds \right| \leq C_{a}^{,} h_T^{2k+3} \| \breve{g} \|_{2k+2,\infty,e_T} \| \breve{v} \|_{2k+2,\infty,e_T}.
\end{equation}
On the other hand, we note that ${\bf n}$ can be expressed only in terms of $x$ in a neighborhood of $\Gamma_T$ and moreover 
${\bf n}_{|(A_T,x,y)} = (-\sin{\varphi};\cos{\varphi})$. Hence we can uniformly bound above the derivatives of ${\bf n}$ up to the $2k+2$-th order as follows.\\
Denoting by ${\mathcal A}^{(j)}(x)$ the $j$-th order derivative for $j>3$ of a function or vector field ${\mathcal A}(x)$, we have 
\[ {\bf n}^{'} = (-\cos{\varphi};-\sin{\varphi}) \varphi^{'}. \]
where $\varphi$ is defined by $\tan{\varphi}:= \phi^{'}$, with $\varphi \in [-\pi/2,\pi/2]$. Therefore  
\[ \displaystyle \varphi^{'} = \displaystyle \frac{\phi^{''}}{[1+(\phi^{'})^2]^2}. \]
Noticing that $|\phi^{'}|$ is uniformly bounded above by a constant times $h_T$ in $e_T$ and the functions $|\phi^{(j)}|$ are all uniformly bounded above by mesh-independent constant $C^{(j)}_{\phi}$ for $2 \leq j \leq 2k+2$, we easily obtain
\begin{equation}
\label{errgv6} 
| {\bf n}^{'} | \leq C^{(1)}_{\bf n}:= C^{(2)}_{\phi}. 
\end{equation}
Performing repeatedly the same type of calculations for higher order derivatives of ${\bf n}$, we conclude with no difficulty that the moduli of the derivatives ${\bf n}^{(j)}$ of any order up to the $2k+2$-th are uniformly bounded above independently of $T$ by constants $C^{(j)}_{\bf n}$, expressed in terms of the constants $C^{(i)}_{\phi}$ for $i=1,2,\ldots,j+1$ that is 
\begin{equation}
\label{errgv7} 
| {\bf n}^{(j)} | \leq C^{(j)}_{\bf n} \mbox{ for } j=1,2,\ldots,2k+2. 
\end{equation} 
Now using the chain rule, we have
\begin{equation}
\label{errgv8}
\breve{g}^{'} = \displaystyle \left[ \frac{\partial {\bf p}}{\partial x} + \frac{\partial {\bf p}}
{\partial y} \phi^{'} \right] \cdot {\bf n} + {\bf p} \cdot {\bf n}^{'}.
\end{equation}
Then the following bound readily derives from \eqref{errgv6} and \eqref{errgv8}:
\begin{equation}
\label{errgv9}
|\breve{g}^{'}(x)| \leq C^{(1)}_g \sqrt{|[\nabla {\bf p}](x)|^2 +  | {\bf p}(x)|^2} \; \forall x \in [0,l_T].
\end{equation}
with $C^{(1)}_g := \sqrt{1+[C^{(1)}_{\phi}]^2 + [C^{(1)}_{\bf n}]^2}$.\\
Further differentiating $\breve{g}$, taking into account \eqref{errgv7} and noticing that $D^i {\bf p}$ 
vanishes identically for $i > k+1$, akin to Theorem 4.4 of \cite{ZAMM} we 
find out that there exists a constant $C^{(2k+2)}_g$ independent of $T$ expressed in terms of $C^{(j)}_{\bf n}$ and $C^{(j)}_{\phi}$ 
such that, for $j=1,2,\ldots,2k+2$,
\begin{equation}
\label{errgv10}
|\breve{g}^{(j)}(x)| \leq C^{(2k+2)}_g \displaystyle \left\{ \sum_{i=0}^{k+1}  |[D^i {\bf p}](x)|^2 \right\}^{1/2} \; \forall x \in [0,l_T]
\end{equation}
As for $\breve{v}$, in the light of the above calculations for $\breve{g}$, and observing that $D^i v$ vanishes identically for $i > k$, it is even easier to infer the existence of a mesh-independent constant $C^{(2k+2)}_v$ expressed in terms of $C^{(j)}_{\phi}$ 
such that, for $j=1,2,\ldots,2k+2$,
\begin{equation}
\label{errgv11}
|\breve{v}^{(j)}(x)| \leq C^{(2k+2)}_v \displaystyle \left\{ \sum_{i=0}^{k}  | [D^i v](x)|^2 \right\}^{1/2} \; \forall x \in [0,l_T].
\end{equation}
Plugging \eqref{errgv10} and \eqref{errgv11} into \eqref{errgv5} we obtain,
\begin{equation}
\label{errgv12}  
\int_{\Gamma_T} {\bf p} \cdot {\bf n}\; v \; ds  \leq C_{a}^{'} C^{(2k+2)}_g C^{(2k+2)}_v (k+1) h_T^{2k+3} \| {\bf p} \|_{k+1,\infty,T^{\Gamma}} \| v \|_{k,\infty,T^{\Gamma}}
\end{equation}
Applying \eqref{L2TDelta} to each norm in \eqref{errgv12} and using $k$ times a classical inverse inequality (see e.g. \cite{ErnGuermond}) to bound above $\| v \|_{k,T}$ by $h_T^{-k} \| v\|_{0,T}$ multiplied by a constant, we immediately come up with \eqref{errgv} for a suitable constant ${\mathcal C}_{\Gamma_1}$. \QED

\section{Well-posedness and uniform stability of the approximate problem}

\hspace{4mm} Besides establishing that problem \eqref{Dmixh} has a unique solution, as a consequence of the fact that the bilinear form $c_h$ given by \eqref{bilinearformch}-\eqref{ahbhdh} is weakly coercive, we will prove in Theorems \ref{theo1} and \ref{nonconvexstab} hereafter that its weak coercivity constant $\gamma_h$ is mesh-independent. This implies the uniform stability of \eqref{Dmixh}. 
According to \cite{AsyMVF}, this result is a consequence of the validity of uniform \textit{inf-sup} conditions for the bilinear forms $a_h$, $b_h$and $d_h$ defined in \eqref{ahbhdh}. Each one of them is worked out in one of the three subsections that follow. \\
\indent At the end of  Subsection 5.3 the \textit{inf-sup} constant $\gamma_h$ of $c_h$ will be denoted by either $\gamma$ or $\gamma^{'}$, depending on whether $\Omega$ is convex or not. Both are exhibited therein as uniform stability constants derived from the resulting uniform \textit{inf-sup} constants of $a_h$, $b_h$ and $d_h$ denoted respectively by $\alpha$, $\beta$ and $\delta$ or $\delta^{'}$ accordingly. More precisely the corresponding \textit{inf-sup} conditions are \eqref{infsupRS}, \eqref{infsupRT} and \eqref{infsupRTc} (or yet \eqref{infsupRTch}) stated in Theorems \ref{lemmars}, \ref{QhVh} and \ref{PhVh} (or yet \ref{PhVhprime}).
  
\subsection{Uniform \textit{inf-sup} condition for bilinear form $a_h$ on $Ker(div)$}

\begin{theorem}
\label{lemmars}
Let 
\begin{equation}
\label{RhkShk}
\begin{array}{l}
{\bf R}_h:=\{ {\bf r}| \; {\bf r} \in {\bf P}_h, \mbox{ and } d_h({\bf r},v) = 0\; \forall v \in V_h\} \mbox{ and }\\
 \\ 
{\bf S}_h:=\{ {\bf s}| \; {\bf s} \in {\bf Q}_h, \mbox{ and } b_h( v, {\bf s}) = 0\; \forall v \in V_h\}.
\end{array}
\end{equation}
Provided $h$ is sufficiently small, there exists a constant $\alpha > 0$ independent of $h$ such that,
\begin{equation}
\label{infsupRS}
\forall {\bf r} \in {\bf R}_h \; \exists {\bf s} \in {\bf S}_h \setminus \{{\bf 0}\} \mbox{ such that } \displaystyle \frac{a_h({\bf r},{\bf s})}{[|{\bf s}|]_{h}} \geq \alpha [|{\bf r}|]_{h}.
\end{equation}
\end{theorem}

\prov 
It is easy to see that the divergence of any field in $RT_k(T)$ is a polynomial in $P_{k}(T)$ whose $k+1$ terms of degree $k$ have the same coefficients as those of the terms of degree $k+1$ of the field itself multiplied by $k+2$. As a matter of fact, the divergence operator is a surjection from $RT_k(T)$ onto $P_k(T)$. Hence the divergence of every field in either ${\bf R}_h$ or ${\bf S}_h$ vanishes identically. We also observe that the restriction to $T$ of a field in either ${\bf R}_h$ or ${\bf S}_h$ is a field in a subspace $RT^{'}_k(T)$ of $[P_k(T)]^2$, whose dimension is $2 dim P_k(T)-k(k+1)/2$, that is, $l_k:=(k+4)(k+1)/2$. \\
For $T \in {\mathcal T}_h$ let $\{ {\bf q}_{i,T}^{'} \}_{i=1}^{l_k}$ be a basis of $RT^{'}_k(T)$ associated with $l_k$ DOFs $\{ {\mathcal G}_{i,T} \}_{i=1}^{l_k}$ such that ${\mathcal G}_{i,T}({\bf q}_{j,T}^{'})=\delta_{ij}$. Since the normal component of a field in ${\bf S}_h$ vanishes identically on ${\Gamma}_{1,h}$, for every $T \in {\mathcal S}_{1,h}$ we choose these DOFs in such a way that the first $k+1$ ones are the field's outer normal components along the edge $e_T$ at the $k+1$ Gauss-quadrature points $M_i$ of $e_T$. \\ 
Now given a field ${\bf r} \in {\bf R}_h$, let ${\bf s} \in {\bf S}_h$ be such that
\begin{itemize}
\item $\forall T \in {\mathcal T}_h \setminus {\mathcal S}_{1,h}$ ${\bf s}_{|T} \equiv {\bf r}_{|T}$ ;
\item $\forall T \in {\mathcal S}_{1,h}$ ${\mathcal G}_{i,T}({\bf s}_{|T})$ = ${\mathcal G}_{i,T}({\bf r}_{|T})$, $i=k+2,\ldots,l_k$. 
\end{itemize}
Since ${\bf s}$ only differs from ${\bf r}$ in a triangle $T \in {\mathcal S}_{1,h}$, we have 
\begin{equation}
\label{differs}
{\bf s}_{|T} - {\bf r}_{|T} = \displaystyle \sum_{i=1}^{k+1} [{\mathcal G}_{i,T}({\bf s}_{|T}) - {\mathcal G}_{i,T}({\bf r}_{|T})] {\bf q}_{i,T}^{'}.
\end{equation}
Abusively denoting ${\bf r}_{|T^{\Gamma}}$ by ${\bf r}_{|T}$, we note that ${\mathcal G}_{i,T}({\bf s}_{|T}) = 0$ and that ${\mathcal G}_{i,T}({\bf r}_{|T}) = 
{\bf r}_{|T}(M_i)\cdot {\bf n}_T$, for $1 \leq i \leq k1$. In view of this and taking into account that ${\bf r}_{|T}(N_i)\cdot {\bf n}_i = 0$, we have,
\begin{equation}
\label{diffG}
{\mathcal G}_{i,T}({\bf s}_{|T}) - {\mathcal G}_{i,T}({\bf r}_{|T}) = [{\bf r}_{|T}(M_i)\cdot {\bf n}_T - {\bf r}_{|T}(N_i)\cdot {\bf n}_T] +
[{\bf r}_{|T}(N_i)\cdot {\bf n}_T - {\bf r}_{|T}(N_i)\cdot {\bf n}_i].
\end{equation}
Now we resort to Propositions \ref{prop01} and \ref{prop02} to write,
\begin{equation}
\label{diffG1}
|{\mathcal G}_{i,T}({\bf s}_{|T}) - {\mathcal G}_{i,T}({\bf r}_{|T})| \leq |{\bf r}_{|T}(M_i) - {\bf r}_{|T}(N_i)| +
|{\bf r}_{|T}(N_i)| |{\bf n}_T - {\bf n}_i| \leq C_0 \| \nabla {\bf r} \|_{0,\infty,T^{\Gamma}} h_T^2 +
C_{\Gamma} h_T \| {\bf r} \|_{0,\infty,T^{\Gamma}}.
\end{equation}
From Lemma \ref{TDelta} and a classical inverse inequality, 
we infer the existence of a mesh-independent constant $C_{\mathcal G}$ such that
\begin{equation}
\label{diffG2}
|{\mathcal G}_{i,T}({\bf s}_{|T}) - {\mathcal G}_{i,T}({\bf r}_{|T})| \leq  C_{\mathcal G} \| {\bf r}_{|T} \|_{0,T}.
\end{equation}
On the other hand, taking into account that there exists a mesh-independent constant $C_{\bf q}^{'}$ such that $\| {\bf q}_{i,T}^{'} \|_{0,T} \leq C_{\bf q}^{'} h_T$ for $i=1,\ldots, k+1$, combining \eqref{differs} and \eqref{diffG2} we easily obtain with $C_k^{'}:= (k+1) C_{\mathcal G} C_{\bf q}^{'}$:
\begin{equation}
\label{differs1}
\| {\bf s}_{|T} - {\bf r}_{|T} \|_{0,T} \leq C_{k}^{'} \| {\bf r}_{|T} \|_{0,T} h_T,
\end{equation}
which immediately yields: 
\begin{equation}
\label{differs2}
[| {\bf s} - {\bf r} |]_{h} \leq C_{k}^{'} h [| {\bf r} |]_{h}.
\end{equation}
From \eqref{differs2} we easily conclude that $[| {\bf s} |]_{h} \leq 
(1+ C_{k}^{'} h) [| {\bf r} |]_{h}$ and $[ {\bf s}, {\bf r}]_{h} \geq (1-C_{k}^{'} h) [| {\bf r} |]^2_h$. Hence, recalling Proposition \ref{prop1}, for $h \leq \max[\hat{C}_k,C_{k}^{'}]^{-1}/2$, the result follows with $\alpha = 1/3$.\rule{2mm}{2mm} \\

In the sequel we will use \eqref{infsupRS} in the following equivalent form: 
\begin{equation}
\label{infsupRSbis}
\forall {\bf r} \in {\bf R}_h \; \exists {\bf s} \in {\bf S}_h \mbox{ with } [| {\bf s} |]_h =1 \mbox{ such that } 
a_h({\bf r},{\bf s}) \geq \alpha [|{\bf r}|]_{h}.
\end{equation} 

\subsection{Uniform \textit{inf-sup} condition for bilinear form $b_h$}

\hspace{4mm} In order to prove this second crucial result} for establishing the uniform stability of problem \eqref{Dmixh}, namely, conditions \eqref{infsupRT} hereafter, we need the following preparatory material: \\ 
Similarly to \cite{BrezziFortin}, we make use of the $L^2(\Omega_h)$-projection operator $\pi_h$ onto the espace $V_h$ and consider as well the interpolation operator $\Pi_h$ onto the space ${\bf Q}_h$. The latter applies to fields in the space ${\bf W}_h$ defined in connection with a real number $s$ strictly greater than two, i.e.,
$${\bf W}_h:= \{ {\bf q} | \; {\bf q} \in {\bf H}(div,\Omega_h) \cap [L^s(\Omega_h)]^2, \; {\bf q} \cdot {\bf n}_h = 0 \mbox{ on } \Gamma_{1,h} \},$$ 
${\bf W}_h$ is equipped with the norm $\| (\cdot) \|_{{\bf W}_h}:= \| (\cdot) \|_{0,s,h} + \| \nabla \cdot (\cdot) \|_{0,h}$. \\
\indent Let $\bar{V}$ be the space $L^2(\Omega_h)$ if $length(\Gamma_0)>0$ or $L^2_0(\Omega_h)$ otherwise. 

\begin{e-proposition}
\label{Omegah}
There exists a constant $\beta_{\bf W} > 0$ independent of $h$ such that,
\begin{equation}
\label{infsupLBB}
\forall v \in \bar{V} \displaystyle \sup_{{\bf q} \in {\bf W}_{h} \mbox{ s.t. } \| {\bf q} \|_{{\bf W}_h}=1} b_h(v,{\bf q}) \geq \beta_{\bf W}  \| v \|_{0,h}.
\end{equation}  
\end{e-proposition}
\prov 
Let $v$ be given in $\bar{V}$ and ${\bf g} := -\nabla z$, where $z$ is the unique solution of the Poisson problem 
\begin{equation}
\label{Poissaux}
\left\{
\begin{array}{l}
-\Delta z = v \mbox{ in } \Omega_h \\
z=0 \mbox{ on } \Gamma_{0,h} \mbox{ if } length(\Gamma_{0,h}) \neq 0 \mbox{ and } \int_{\Omega_h} z =0 \mbox{ otherwise}.\\
\partial z/\partial
 n_h =0 \mbox{ on } \Gamma_{1,h}.
\end{array}
\right.
\end{equation}
Even if $length(\Gamma_0) > 0$, according to \cite{Grisvard}, $z \in H^r(\Omega_h)$ for some $r \in (3/2,2)$ and moreover by the Rellich-Kondrachev Theorem \cite{Adams}, $\nabla z \in L^{s}(\Omega_h)$ for some $s > 2$. Actually, due to both arguments we know that there exists a constant $C_{s,h}$ depending on $\Omega_h$ such that 
\begin{equation}
\label{Csh}
\| \nabla z \|_{0,s,h} \leq C_{s,h} \| v \|_{0,h}.
\end{equation}
 It trivially follows that  
\begin{equation}
\displaystyle \frac{b_h(v,{\bf g})}{\| {\bf g} \|_{0,s,h}} \geq C_{s,h}^{-1} \| v \|_{0,h}
\end{equation}
and, as one can easily check, 
\begin{equation}
\displaystyle \frac{b_h(v,{\bf g})}{\| {\bf g} \|_{{\bf W}_h}} \geq \beta_h \| v \|_{0,h},
\end{equation}
with $\beta_h = (C_{s,h}+1)^{-1}$.\\ 
 In APPENDIX I we claim and prove that the constant $C_{s,h}$ has an upper bound $C_s$ independent of $h$. Hence the constant $\beta_h$ has a lower bound $\beta_{W}>0$ independent of $h$. 
Setting $\tilde{\bf q} = {\bf g}/\| {\bf g} \|_{{\bf W}_h}$ \eqref{infsupLBB} follows. \QED \\

Now we are able to establish the uniform \textit{inf-sup} condition for $b_h$ in 
\begin{theorem}
\label{QhVh}
There exists a constant $\beta>0$ independent of $h$ such that 
\begin{equation}
\label{infsupRT}
\forall v \in V_h \displaystyle \sup_{{\bf q} \in {\bf Q}_h \mbox{ s.t. } [| {\bf q} |]_{h}=1} b_h(v,{\bf q}) \geq \beta \| v \|_{0,h}.
\end{equation}  
\end{theorem}

\prov
According to \eqref{infsupLBB} we know that the divergence operator has a continuous lifting from $\bar{V}$ to ${\bf W}_h$. Moreover, following the same steps as in Chapter III of \cite{BrezziFortin}, we can assert that 
there exists a constant $C_{\Pi}$ independent of $h$ such that 
$$[| \Pi_h {\bf q} |]_h \leq C_{\Pi} \| {\bf q} \|_{{\bf W}_h} \; \forall {\bf q} \in {\bf W}_h.$$  
On the other hand, quoting again \cite{BrezziFortin}, for every $v \in V_h$ there exists ${\bf q} \in {\bf Q}_h$ such that 
$\nabla \cdot {\bf q} =v$. After a quick check, all the other conditions specified on page 138, Chapter IV of \cite{BrezziFortin} for the commuting diagram (I.23) given therein, apply to the pairs of spaces ${\bf W}_h-\bar{V}$ on the upper part and ${\bf Q}_h-V_h$ on the lower part, while $\Pi_h$ and $\pi_h$ operate from one space to the other on its left and right parts.\\
Summarizing, $\forall {\bf q } \in {\bf W}_h$ it holds:
\begin{equation}
\label{preinfsup}
\left\{
\begin{array}{ll}
b_h(v,\Pi_h {\bf q}) = b_h(v,{\bf q}) & \forall v \in V_h \\
\mbox{ } [| \Pi_h {\bf q} |]_h \leq C_{\Pi} \| {\bf q} \|_{{\bf W}_h}. &
\end{array}
\right.
\end{equation}
Now, taking an arbitrary $v \in V_h$ and recalling \eqref{infsupLBB}, 
let ${\bf q} \in {\bf W}_h$ be such that $\| {\bf q} \|_{{\bf W}_h} =1$ and $b_h(v,{\bf q}) \geq \beta_{\bf W} \| v \|_{0,h}$ 
Then, from \eqref{preinfsup} we immediately conclude that \eqref{infsupRT} holds with $\beta = \beta_{\bf W}/C_{\Pi}$. \QED
 
\subsection{Uniform \textit{inf-sup} conditions for bilinear forms $d_h$ and $c_h$}

\hspace{4mm} Taking into account the results of Subsections 4.1 and 4.2, according to \cite{AsyMVF} a uniform \textit{inf-sup} holds for bilinear form $c_h$, as long as a uniform \textit{inf-sup} condition also holds for bilinear form $d_h$.    
In order to establish the latter, the following additional preparatory material is of paramount importance.\\
\indent First of all we introduce for convenience the subset ${\mathcal F}^{'}$ of \textit{triangulations of interest} in the family ${\mathcal F}$. Henceforth the latter expression will qualify meshes ${\mathcal T}_h$ which are not too coarse. We refrain from further elaborating about this concept, since it implicitly goes hand in hand with an argument repeatedly invoked throughout this paper, namely, the assumption that $h$ is sufficiently small. \\
We shall need an interpolation operator onto the space ${\bf P}_h$ for fields defined beyond $\Omega_h$. The other way around, fields defined in $\Omega$ are not necessarily defined in $\Omega_h$. Thus, in the remainder of this section and on we will work with extensions thereof outside $\Omega$, say, in a smooth domain $\tilde{\Omega}$ strictly containing $\bar{\Omega}$. More specifically we shall use the following auxiliary material: \\
We require $\tilde{\Omega}$ to be sufficiently large to strictly contain the sets $\tilde{\Omega}_h$ for all ${\mathcal T}_h \in {\mathcal F}^{'}$. We denote by $\tilde{\Gamma}$ the boundary of $\tilde{\Omega}$ and define $\Delta_{\Omega} = \tilde{\Omega} \setminus \bar{\Omega}$.  Notice that the boundary $\partial \Delta_{\Omega}$ of $\Delta_{\Omega}$ equals $\Gamma \cup \tilde{\Gamma}$. \\
Now for any ${\bf m} \in [W^{t,r}(\Omega)]^2 \cap {\bf H}(div,\Omega)$, with $t \in \Re^{+}$ and $r \geq 1$, let us assume that there exists an extension $\tilde{\bf m}$ of ${\bf m}$ to the whole $\tilde{\Omega}$ such that $\tilde{\bf m} \in [W^{t,r}(\tilde{\Omega})]^2 \cap {\bf H}(div,\tilde{\Omega})$. Referring to an argument highlighted in \cite{BrezziFortin}, as long as $L^{s}(\tilde{\Omega})$ is continuously embedded into $W^{t,r}(\tilde{\Omega})$ for some $s >2$, we may define an interpolation operator $\tilde{\Pi}_h :[W^{t,r}(\tilde{\Omega})]^2 \cap {\bf H}(div,\tilde{\Omega}) \rightarrow {\bf P}_h$ by the following relations, where 
$\partial T$ represents the boundary of a triangle $T \in {\mathcal T}_h$:
\begin{equation}
\label{Phkinterpolate}
\left\{
\begin{array}{l}
\forall T \in {\mathcal T}_h \setminus {\mathcal S}_{1,h}, \;\int_e [\tilde{\Pi}_h \tilde{\bf m} - \tilde{\bf m}] \cdot {\bf n}_T \; v = 0 \; \forall v \in P_k(e) \; \forall \mbox{ edge } e \subset \partial T; \\
\forall T \in {\mathcal S}_{1,h} \; \int_e [\tilde{\Pi}_h \tilde{\bf m} - \tilde{\bf m}] \cdot {\bf n}_T \; v = 0 \; \forall v \in P_k(e) \; \forall \mbox{ edge } e \subset \partial T \setminus e_T ; \\
\mbox{For } k>0, \; \int_T [\tilde{\Pi}_h \tilde{\bf m} - \tilde{\bf m}] \cdot {\bf o} = 0 \; \forall {\bf o} \in [P_{k-1}(T)]^2.  
\end{array}
\right.
\end{equation}
We can assert that $\tilde{\Pi}_h \tilde{\bf m} \in {\bf H}(div,\tilde{\Omega})$ $\forall \tilde{\bf m} \in [W^{t,r}(\tilde{\Omega})]^2 \cap {\bf H}(div,\tilde{\Omega})$. Moreover, as long as for every bounded Lipschitz domain $D$ of $\Re^2$, $W^{t,r}(D)$ is continuously embedded into $W^{q,s}(D)$ for some $q \geq 0$ and $s \geq 1$, there exists a constant $C^{t,r}_{q,s}$ not depending on $T$ such that the following estimate holds
\begin{equation}
\label{extension}
\| \tilde{\Pi}_h \tilde{\bf m} - \tilde{\bf m} \|_{q,s,T^{\Gamma}} \leq C_{q,s}^{t,r} h_T^{2/s-2/r} h_T^{-q+t} | \tilde{\bf m} |_{t,r,T^{\Gamma}} \; \forall T \in {\mathcal S}_{1,h} \; \forall {\bf m} \in [W^{t,r}(\Omega)]^2,
\end{equation}
The proof of \eqref{extension} for $q=0$, $s=2$, $t=k+1$ and $r=2$ will be given in Section 6 in the framework of an estimate for the interpolation error in the norm $[| \cdot |]_h$. However, since only Sobolev norms are involved in \eqref{extension}, it is easy to establish that it holds true, thanks to standard estimates (cf. \cite{ErnGuermond,Arcangeli}).\\

\noindent \underline{\textbf{A) The convex case}}\\
\indent We pursue the study of the uniform stability of $d_h$ by first assuming that the problem-definition domain is convex. This is because the analysis, though still complex, is simpler in this case.\\ 
First of all we introduce a suitable counterpart of space ${\bf W}_h$, namely,

$${\bf W}:= \{ {\bf q} | \; {\bf q} \in {\bf Q} \cap [H^{1/2+\varepsilon}(\Omega)]^2\},$$
 
\noindent where $\varepsilon$ is a small strictly positive real number. ${\bf W}$ is equipped with the norm $\| (\cdot) \|_{\bf W}:= \| (\cdot) \|_{1/2+\varepsilon} + \| \nabla \cdot (\cdot) \|_{0}$ (for details on Sobolev spaces of fractional order we refer to \cite{DiNezza}).\\
We next prove

\begin{e-proposition}
\label{MixedPoisson}
For every $v \in V_h$ define a function $\tilde{v}$ in $\Omega$ associated with $v$, as the extension of $v$ by zero in $\bar{\Omega} \setminus \bar{\Omega}_h$. In doing so, it holds
\begin{equation}
\label{infsupOmega}
\exists \delta_{\bf W}>0 \mbox{ such that } \; \forall v \in V_h \displaystyle \sup_{{\bf p} \in {\bf W} \setminus {\bf 0}_{\bf W}} 
\frac{(\nabla \cdot {\bf p},\tilde{v})}{\| {\bf p} \|_{\bf W}} \geq \delta_{\bf W} \| \tilde{v} \|_{0}.
\end{equation}
\end{e-proposition}
   
\prov
Noting that $\tilde{v} \in L^2_0(\Omega)$ if $length(\Gamma_0)=0$, we consider a function $\tilde{z}$ defined in $\Omega$ to be the unique solution of the following Poisson problem:
\begin{equation}
\label{Poissontilde}
\left\{
\begin{array}{l}
- \Delta \tilde{z} = \tilde{v} \mbox{ in } \Omega \\
\tilde{z} = 0 \mbox{ on } \Gamma_0 \mbox{ if } length(\Gamma_0) >0 \mbox{ and } \int_{\Omega} \tilde{z} = 0 \mbox{ if } length(\Gamma_0) =0 \\
\partial \tilde{z}/\partial n = 0 \mbox{ on } \Gamma_1.
\end{array}
\right.
\end{equation}
Setting $\tilde{\bf p}:= -\nabla \tilde{z}$ we have $\tilde{\bf p} \in {\bf Q}$ and, according to \cite{Grisvard},   
$\tilde{\bf p} \in [H^{1/2+\varepsilon}]^2$. Furthermore, there exists a constant $\tilde{C}$ depending only on $\Omega$ and $\Gamma_1$ such that
\begin{equation}
\label{Continuity}
\| \tilde{\bf p} \|_{1/2+\varepsilon} \leq \tilde{C} \| \tilde v \|_0.
\end{equation}
Since $(\nabla \cdot \tilde{\bf p},\tilde{v}) = (-\Delta \tilde{z},\tilde{v}) = \| \tilde{v} \|_{0}^2$, using 
\eqref{Continuity}, \eqref{infsupOmega} is trivially seen to hold with $\delta_{\bf W} = (1+\tilde{C})^{-1}$.
\QED \\

As much as $\Pi_h$ was used in the study of $b_h$, the uniform \textit{inf-sup} condition for $d_h$ relies on the properties of the interpolation operator $\tilde{\Pi}_h$ defined in \eqref{Phkinterpolate}. Here we take $t=1/2+\varepsilon$ and $r=2$, to have $\tilde{\Pi}_h : {\bf W} \rightarrow {\bf P}_h$ defined as in \eqref{Phkinterpolate} by replacing $\tilde{\bf m}$ with ${\bf m} \in {\bf W}$. This is because we can do without extensions in $\tilde{\Omega}$, since $\Omega$ is convex and thus the estimate \eqref{extension} also holds for any ${\bf m} \in {\bf W}$ instead of $\tilde{\bf m}$.\\
\indent Let us prove two preliminary results for the interpolation operator $\tilde{\Pi}_h$, namely, 

\begin{e-proposition}
\label{tildePih}
There exists a constant $C_{\tilde{\Pi}}$ independent of $h$ such that
\begin{equation}
\label{boundtildePih}
[| \tilde{\Pi}_h {\bf m} |]_{h} \leq C_{\tilde{\Pi}} \| {\bf m} \|_{\bf W} \; \forall {\bf m} \in {\bf W}.
\end{equation}
\end{e-proposition}

\prov
First of all, from the same arguments as those applying to the operator $\Pi_h$ (see e.g. \cite{BrezziFortin}), it is easy to figure out that, 
for a suitable constant $C^{'}_{\tilde{\Pi}}$ independent of $h$, it holds
\begin{equation}
\label{boundtilPih1}
\| \tilde{\Pi}_h {\bf m} \|_{0,h} \leq C^{'}_{\tilde{\Pi}} \| {\bf m} \|_{0,s,h} \; \forall {\bf m} \in {\bf W}.
\end{equation}
Indeed, although zero values are prescribed to $\tilde{\Pi}_h {\bf m} \cdot {\bf n}$ at some points of $\Gamma_1$, we are not assuming that $\tilde{\Pi}_h {\bf m}$ takes up values of ${\bf m}$ anywhere in a pointwise sense.\\
On the other hand, the Rellich-Kondratchev Theorem ensures that ${\bf W}$ is continuously embedded in $[L^s(\Omega)]^2$ for $2 < s \leq 4$. Hence, from \eqref{boundtilPih1} and the convexity of $\Omega$, we immediately infer the existence of another mesh-independent constant $C^{0}_{\tilde{\Pi}}$ such that
\begin{equation}
\label{boundtilPih2}
\| \tilde{\Pi}_h {\bf m} \|_{0,h} \leq C^{0}_{\tilde{\Pi}} \| {\bf m} \|_{\bf W} \; \forall {\bf m} \in {\bf W}.
\end{equation} 
Now we turn our attention to $\nabla \cdot \tilde{\Pi}_h {\bf m}$, which we denote by $g_h$ for simplicity. \\
First of all, owing to the definition of $\tilde{\Pi}_h$, $\forall T \in {\mathcal T}_h \setminus {\mathcal S}_{1,h}$ we successively have
\begin{equation}
\label{divPihcalTh0}
\left\{
\begin{array}{l}
\| g_h \|_{0,T}^2 = \int_{T} (\nabla \cdot \tilde{\Pi}_h {\bf m} - \nabla \cdot {\bf m})\nabla \cdot \tilde{\Pi}_h {\bf m} + \int_{T} \nabla \cdot \tilde{\Pi}_h {\bf m} \nabla \cdot {\bf m} \\
= \int_{\partial T} (\tilde{\Pi}_h {\bf m} - {\bf m}) \cdot {\bf n}_T \nabla \cdot \tilde{\Pi}_h {\bf m} - \int_T (\tilde{\Pi}_h {\bf m} - {\bf m}) \cdot \nabla \nabla \cdot \tilde{\Pi}_h {\bf m} + \int_{T} \nabla \cdot \tilde{\Pi}_h {\bf m} \nabla \cdot {\bf m} \\
= \int_{T} \nabla \cdot \tilde{\Pi}_h {\bf m} \nabla \cdot {\bf m},
\end{array}
\right.
\end{equation}
which yields
\begin{equation}
\label{divPihcalTh}
\| g_h \|_{0,T} \leq \| \nabla \cdot {\bf m} \|_{0,T} \; \forall T \in {\mathcal T}_h \setminus {\mathcal S}_{1,h}. 
\end{equation}
Now taking $T \in {\mathcal S}_{1,h}$, again by the definition of $\tilde{\Pi}_h$ we successively have 
\begin{equation}
\label{boundivPih1}
\| g_h \|_{0,T}^2 = \oint_{\partial T} \tilde{\Pi}_h {\bf m} \cdot {\bf n}_T g_h - \int_T \tilde{\Pi}_h {\bf m} \cdot \nabla g_h,
\end{equation}
\begin{equation}
\label{boundivPih2} 
\| g_h \|_{0,T}^2 = \oint_{\partial T} {\bf m} \cdot {\bf n}_T g_h - \int_T {\bf m} \cdot \nabla g_h + \int_{e_T} [\tilde{\Pi}_h {\bf m} - {\bf m}]\cdot {\bf n}_T g_h,
\end{equation}

\begin{equation}
\label{boundivPih3} 
\| g_h \|_{0,T}^2 = \int_T \nabla \cdot {\bf m} g_h - \int_{e_T} {\bf m} \cdot {\bf n}_T g_h + \int_{e_T} \tilde{\Pi}_h {\bf m} \cdot {\bf n}_T g_h.
\end{equation}
Since ${\bf m} \cdot {\bf n}$ vanishes on $\Gamma_1$, referring to Figure 1 and using the set $T^{\Gamma} = T \cup \Delta_T$ we can write

\begin{equation}
\label{boundivPih4} 
\left\{
\begin{array}{l}
\| g_h \|_{0,h}^2 = \displaystyle \sum_{T \in {\mathcal S}_{1,h}} [J_{1,T}({\bf m},g_h)+J_{2,T}({\bf m},g_h)+J_{3,T}({\bf m},g_h)] \mbox{ where }\\
J_{1,T}({\bf m},g_h):= \int_{T^{\Gamma}} \nabla \cdot {\bf m} g_h, \\
J_{2,T}({\bf m},g_h):= \int_{\Delta_T} {\bf m} \cdot \nabla g_h, \\
J_{3,T}({\bf m},g_h):= \int_{e_T} \tilde{\Pi}_h {\bf m} \cdot {\bf n}_T g_h.
\end{array}
\right.
\end{equation}
By the Cauchy-Schwarz inequality, followed by \eqref{L2TDelta} and Young's inequality, we obtain for 
$C_{J1}:= {\mathcal C}_0^{+}/{\mathcal C}_0^{-}$
\begin{equation}
\label{boundJ1}
|J_{1,T}({\bf m},g_h)| \leq \displaystyle \frac{\| g_h \|_{0,T}^2}{4} + C_{J1}^2 \| \nabla \cdot {\bf m} \|_{0,T^{\Gamma}}^2.
\end{equation}
Now, since ${\bf m} \in [L^4(\Omega)]^2$, using H\"older's inequality with $s=4$ and $r=4/3$ we come up with
\begin{equation}
\label{bound1J2}
|J_{2,T}({\bf m},g_h)| \leq \displaystyle \left[ \int_{\Delta_T} |{\bf m}|^s \right]^{1/s} \displaystyle 
\left[ \int_{\Delta_T} |\nabla g_h|^r \right]^{1/r} \leq \displaystyle \left[ \int_{\Delta_T} |{\bf m}|^s \right]^{1/s}  
C_0^{1/r} h_T^{3/r} \| \nabla g_h \|_{0,\infty,T^{\Gamma}}.
\end{equation}
Then resorting again to \eqref{L2TDelta} we further have
\begin{equation}
\label{bound2J2}
|J_{2,T}({\bf m},g_h)|
\leq C_0^{1/r} {\mathcal C}_0^{+} \displaystyle \left[ \int_{T^{\Gamma}} |{\bf m}|^s \right]^{1/s} h_T^{3/r-1} \| \nabla g_h \|_{0,T}.
\end{equation}
From \eqref{bound2J2} together with well known inverse inequalities (see e.g. \cite{ErnGuermond}), we infer the existence of a mesh-independent constant $C^{'}_{J2}$ such that 
\begin{equation}
\label{boundJ2}
|J_{2,T}({\bf m},g_h)|
\leq C^{'}_{J2} \displaystyle \left[ \int_{T^{\Gamma}} |{\bf m}|^s \right]^{1/s} h_T^{3/r-2} \| g_h \|_{0,T}.
\end{equation}
As for $J_{3,T}({\bf m},g_h)$, denoting by $\omega_i$ the weights of the Gauss quadrature formula on $e_T$ we have
\begin{equation}
\label{bound1J3}
\begin{array}{l}
|J_{3,T}({\bf m},g_h)| \leq l_T \displaystyle \sum_{i=1}^{k+1} \omega_i |[\tilde{\Pi}_h {\bf m}](M_i) \cdot {\bf n}_T| \; | g_h(M_i) | \\
\leq h_T \| g_h \|_{0,\infty,T} \displaystyle \sum_{i=1}^{k+1} |[\tilde{\Pi}_h {\bf m}](M_i) \cdot {\bf n}_T |.
\end{array}
\end{equation}
Then, recalling the definition of the points $N_i$ on $\Gamma_1$ and the notation ${\bf n}_i = {\bf n}(N_i)$, \eqref{bound1J3} yields
\begin{equation}
\label{bound2J3}
|J_{3,T}({\bf m},g_h)| \leq h_T \| g_h \|_{0,\infty,T} \displaystyle \sum_{i=1}^{k+1} |[\tilde{\Pi}_h {\bf m}](M_i) \cdot {\bf n}_T - [\tilde{\Pi}_h {\bf m}](N_i) \cdot {\bf n}_i|.
\end{equation}
Adding and subtracting $[\tilde{\Pi}_h {\bf m}](N_i) \cdot {\bf n}_T$ inside the absolute value on the right hand side of \eqref{bound2J3}, owing to Propositions \ref{prop01} and \ref{prop02} we easily obtain, for a suitable mesh-independent constant $C_{J3}^{'}$
\begin{equation}
\label{bound3J3}
|J_{3,T}({\bf m},g_h)| \leq C_{J3}^{'} h_T \| g_h \|_{0,\infty,T} [ h_T \| \tilde{\Pi}_h {\bf m} \|_{0,\infty,T^{\Gamma}} + h_T^2 
\| \nabla \tilde{\Pi}_h {\bf m} \|_{0,\infty,T^{\Gamma}}.
\end{equation}
Finally, applying standard inverse inequalities (see e.g. \cite{ErnGuermond}) and \eqref{L2TDelta} to \eqref{bound3J3}, followed by Young's inequality, we easily come up with a constant $C_{J3}$ independent of $h$ such that 
\begin{equation}
\label{boundJ3}
|J_{3,T}({\bf m},g_h)| \leq \displaystyle \frac{\| g_h \|_{0,T}^2}{4} + C_{J3}^2 \| \tilde{\Pi}_h {\bf m} \|_{0,T}^2.
\end{equation}
Now all that is left to do is to sum up over ${\mathcal S}_{1,h}$ the three upper bounds \eqref{boundJ1}, \eqref{boundJ2} and \eqref{boundJ3}.
While on the one hand the first and the third ones are trivially handled, on the other hand we still have some work to do, as far as \eqref{boundJ2} is concerned. \\
Using H\"older's inequality with $s=4$ and $r=4/3$, together with a trivial upper bound, we obtain
\begin{equation}
\label{sum1J2}
\displaystyle \sum_{T \in {\mathcal S}_{1,h}} |J_{2,T}({\bf m},g_h)|
\leq C^{'}_{J2} \| {\bf m} \|_{0,s} \displaystyle \left[ \sum_{T \in {\mathcal S}_{1,h}} h_T^{3-2r} \| g_h \|_{0,T}^r \right]^{1/r}
\end{equation}
Next we use again H\"older's inequality in the summation on the right hand side with conjugate exponents $s^{'}$ and $r^{'}$ for  
$r^{'}=2/r=3/2$, to obtain 
\begin{equation}
\label{sum2J2}
\begin{array}{l}
\displaystyle \sum_{T \in {\mathcal S}_{1,h}} |J_{2,T}({\bf m},g_h)|
\leq C^{'}_{J2} \| {\bf m} \|_{0,s} \displaystyle \left[ \sum_{T \in {\mathcal S}_{1,h}} h_T^{(3-2r)s^{'}} \right]^{1/(r s^{'})} \| g_h \|_{0,h} \\
= C^{'}_{J2} \| {\bf m} \|_{0,s} \displaystyle \left[ \sum_{T \in {\mathcal S}_{1,h}} h_T \right]^{1/4} \| g_h \|_{0,h}.
\end{array}
\end{equation}
Obviously enough, there exists a constant $C_{\Gamma_1}$ independent of $h$ such that $\displaystyle \sum_{T \in {\mathcal S}_{1,h}} h_T 
\leq C_{\Gamma_1} length(\Gamma_1)$. Thus using Young's inequality, for a suitable mesh-independent constant $C_{J2}$, \eqref{sum2J2} easily yields 
\begin{equation}
\label{sumJ2}
\displaystyle \sum_{T \in {\mathcal S}_{1,h}} |J_{2,T}({\bf m},g_h)|
\leq \displaystyle \frac{\| g_h \|_{0,h}^2}{4} + C_{J2}^2 \| {\bf m} \|_{0,4}^2.
\end{equation}
Summing also up all the terms in \eqref{boundJ1} and \eqref{boundJ3} and plugging into \eqref{boundivPih4} the resulting inequalities together with \eqref{sumJ2}, we readily obtain
\begin{equation}
\label{boundivPih5} 
\displaystyle \sum_{T \in {\mathcal S}_{1,h}} \| \nabla \cdot \tilde{\Pi}_h {\bf m} \|_{0,T}^2 \leq 2 \{\max[C_{J1},C_{J2},C_{J3}]\}^2 \{ \| {\bf m} \|_{0}^2 + \| {\bf m} \|_{0,4}^2 + \| \tilde{\Pi}_h {\bf m} \|_{0,h}^2 \}.
\end{equation}
Now, taking into account\eqref{divPihcalTh}, the following upper bound for $\| \nabla \cdot \tilde{\Pi}_h {\bf m} \|_{0,h}$ derives from \eqref{boundivPih5}:
\begin{equation}
\label{boundivPih6}
 \| \nabla \cdot \tilde{\Pi}_h {\bf m} \|_{0,h}^2 \leq 2 \{\max[C_{J1},C_{J2},C_{J3}]\}^2 \{ \| {\bf m} \|_{0}^2 + \| {\bf m} \|_{0,4}^2 + \| \tilde{\Pi}_h {\bf m} \|_{0,h}^2 \} + \| \nabla \cdot {\bf m} \|_{0,h}^2.
\end{equation}
Finally, recalling \eqref{boundtilPih2} and the fact that $H^{1/2+\varepsilon}(\Omega)$ is continuously embedded into $L^4(\Omega)$, we come up with a constant $C^D_{\tilde{\Pi}}$ such that
\begin{equation}
\label{boundivPih} 
\| \nabla \cdot \tilde{\Pi}_h {\bf m} \|_{0,h} \leq  C^D_{\tilde{\Pi}} \| {\bf m} \|_{\bf W} \; \forall {\bf m} \in {\bf W},
\end{equation}
which completes the proof with $C_{\tilde{\Pi}}\sqrt{(C_{\tilde{\Pi}}^0)^2+ (C^D_{\tilde{\Pi}})^2}$.\QED \\

\begin{e-proposition}
\label{divPih}
There exists a constant $C_{\delta}$ independent of $h$ such that
\begin{equation}
\label{errdivPih}
d_h(\tilde{\Pi}_h {\bf m} - {\bf m},v)_{h} \leq C_{\delta} h^{1/4} \| {\bf m} \|_{\bf W} \| v \|_{0,h} \; \forall {\bf m} \in {\bf W} \mbox{ and } \forall v \in V_h.
\end{equation}
\end{e-proposition}

\prov
From the definition of $\tilde{\Pi}_h$ we may write $\forall v \in V_h$
\begin{equation}
\label{errdiv}
\begin{array}{l}
d_h(\tilde{\Pi}_h {\bf m} - {\bf m}, v) = \displaystyle \sum_{T \in {\mathcal T}_h} \left\{ \oint_{\partial T} 
[\tilde{\Pi}_h {\bf m} - {\bf m}] \cdot {\bf n}_T v - ( \tilde{\Pi}_h {\bf m} - {\bf m},\nabla v)_{0,T} \right\} \\
= \displaystyle \sum_{T \in {\mathcal S}_{1,h}} \int_{e_T} [\tilde{\Pi}_h {\bf m} - {\bf m}] \cdot {\bf n}_T v.
\end{array}
\end{equation}
Since ${\bf m} \cdot {\bf n}$ vanishes identically on $\Gamma_1$, extending $v$ to $\Delta_T$ in the natural way $\forall T \in {\mathcal S}_{1,h}$ the following expression derives from \eqref{errdiv}: 
\begin{equation}
\label{errdivbis}
d_h(\tilde{\Pi}_h {\bf m} - {\bf m}, v) = \displaystyle \sum_{T \in {\mathcal S}_{1,h}} \left\{ 
\int_{e_T} [\tilde{\Pi}_h {\bf m} - {\bf m}] \cdot {\bf n}_T v - \int_{\Gamma_T}  [\tilde{\Pi}_h {\bf m} - {\bf m}] \cdot {\bf n} v 
+ \int_{\Gamma_T} \tilde{\Pi}_h {\bf m} \cdot {\bf n} v \right\}.
\end{equation} 
Now applying the Divergence Theorem in $\Delta_T$, from \eqref{errdivbis} we obtain rather easily 
\begin{equation}
\label{errdivter}
\left\{
\begin{array}{l}
d_h(\tilde{\Pi}_h {\bf m} - {\bf m}, v) 
=  \displaystyle \sum_{T \in {\mathcal S}_{1,h}} [K_{1,T}({\bf m},v)+K_{2,T}({\bf m},v)+K_{3,T}({\bf m},v)] \\
\mbox{with} \\
K_{1,T}({\bf m},v) :=  \int_{\Gamma_T} \tilde{\Pi}_h {\bf m} \cdot {\bf n} v \\
K_{2,T}({\bf m},v) := -\int_{\Delta_T} \nabla \cdot (\tilde{\Pi}_h {\bf m} - {\bf m})v \\
K_{3,T}({\bf m},v) := -\int_{\Delta_T} (\tilde{\Pi}_h {\bf m}
 - {\bf m}) \cdot \nabla v.
\end{array}
\right.
\end{equation}
First of all the functional $K_{1,T}$ can be dealt with by means of Lemma \ref{GammaT}, that is,
\begin{equation}
\label{errK11}
|K_{1,T}({\bf m},v)| \leq  C_{\Gamma_1} h_T^{k+1} \| \tilde{\Pi}_h {\bf m} \|_{k+1,T} \| v \|_{0,T}. 
\end{equation}
Now we repeatedly resort to fractional inverse inequalities (cf. \cite{Georgoulis}) to infer the existence of a constant $C_{1/2}$ independent of $T$ such that
\begin{equation}
\label{errK12}
|K_{1,T}({\bf m},v)| \leq  C_{1/2} h_T^{1/2} \| \tilde{\Pi}_h {\bf m} \|_{1/2+\varepsilon,T} \| v \|_{0,T} 
\end{equation}
Owing to the fact that, by construction, to each triangle $T \in {\mathcal S}_{1,h}$ exactly $(k+1)(k+2)$ linearly independent DOFs not associated with $\Gamma_T$ are assigned to $\tilde{\Pi}_h$, this operator is invariant as applied to fields in $[P_k(T)]^2$ independently of boundary DOFs. Therefore, on the basis of the well established interpolation theory 
in $\Re^n$, including the case of vector fields (see e.g. \cite{BrennerScott,BrezziFortin}), and recalling also approximation results in fractional norms \cite{Arcangeli}, we can assert that there exists another constant $C_{\varepsilon}$ independent of $T$ such that  
\begin{equation}
\label{errK13}
\| \tilde{\Pi}_h {\bf m} \|_{1/2+\varepsilon,T} \leq C_{\varepsilon} \| {\bf m} \|_{1/2+\varepsilon,T}.
\end{equation}
Setting $C_{K,1}:= C_{1/2} C_{\varepsilon}$, this immediately yields the following estimate for $K_{1,T}({\bf m},v)$: 
\begin{equation}
\label{errK1}
|K_{1,T}({\bf m},v)| \leq C_{K,1} h_T^{1/2} \| {\bf m} \|_{1/2+\varepsilon,T} \| v \|_{0,T}.
\end{equation} 
Next we observe that both $\nabla \cdot \tilde{\Pi}_h {\bf m}$ and $v$ are polynomials in $P_k(T^{\Gamma})$. In doing so we resort to a well known argument (see e.g. \cite{StrangFix}), according to which 
\begin{equation}
\label{DeltaT}
\exists C_{\Delta} \mbox{ independent of } T \mbox{ such that } \| p \|_{0,\Delta_T} \leq C_{\Delta} h_T^{1/2} \|  p \|_{0,T} \; 
\forall p \in P_k(T^{\Gamma}) \mbox{ and } \forall T \in {\mathcal S}_{1,h}.
\end{equation}
Applying \eqref{DeltaT} to $K_{2,T}({\bf m},v)$ we readily come up with
\begin{equation}
\label{errK2}
|K_{2,T}({\bf m},v)| \leq C_{\Delta} h_T^{1/2} ( C_{\Delta} h_T^{1/2} \| \nabla \cdot \tilde{\Pi}_h {\bf m} \|_{0,T} + \| \nabla \cdot {\bf m} \|_{0,T^{\Gamma}}) \| v \|_{0,T}.
\end{equation}
As for $K_{3,T}({\bf m},v)$, we first employ H\"older's inequality to obtain
\begin{equation}
\label{errK31}
|K_{3,T}({\bf m},v)| \leq \| \tilde{\Pi}_h {\bf m} - {\bf m} \|_{0,4,T^{\Gamma}} \| \nabla v \|_{0,4/3,\Delta_T}
\end{equation} 
Next we note that, according to the Rellich-Kondratchev Theorem (cf. \cite{Adams}) $H^{1/2+\varepsilon}(D)$ is continuously embedded in $L^4(D)$ for all domains $D$ in the aforementioned class. Thus, \eqref{extension} yields
\begin{equation}
\label{L4}
\| \tilde{\Pi}_h {\bf m} - {\bf m} \|_{0,4,T^{\Gamma}} \leq C_{0,4} h_T^{\varepsilon} | {\bf m}|_{1/2+\varepsilon,T^{\Gamma}}.
\end{equation}
Taking into account \eqref{L2TDelta}, plugging \eqref{L4} into \eqref{errK31} we come up with 
\begin{equation}
\label{errK32}
|K_{3,T}({\bf m},v)| \leq C_{0,4} C_0^{3/4} h_T^{9/4+\varepsilon}  | {\bf m} |_{1/2+\varepsilon,T^{\Gamma}} \| \nabla v \|_{0,\infty,T^{\Gamma}}
\leq {\mathcal C}_0^{+} C_{0,4} C_0^{3/4} h_T^{5/4+\varepsilon}  | {\bf m} |_{1/2+\varepsilon,T^{\Gamma}} \| \nabla v \|_{0,T}.
\end{equation}
Owing to \eqref{errK32} and a classical inverse inequality, we guarantee the existence of a constant $C_{1/4}$ independent of $T$ such that
\begin{equation}
\label{errK3}
|K_{3,T}({\bf m},v)| \leq C_{1/4} h_T^{1/4} | {\bf m} |_{1/2+\varepsilon,T^{\Gamma}} \| v \|_{0,T}.
\end{equation} 
Now plugging \eqref{errK1}, \eqref{errK2} and \eqref{errK3} into \eqref{errdivter}, we easily obtain for a suitable constant $C_{\Delta}^{'}$ 
\begin{equation}
\label{errdivqua}
|(\nabla \cdot [\tilde{\Pi}_h {\bf m} - {\bf m}] v)_h | \leq C_{\Delta}^{'} h_T^{1/4} (\| \nabla \cdot \tilde{\Pi}_h {\bf m} \|_{0,h} + \|  {\bf m} \|_{\bf W} + | {\bf m} |_{1/2+\varepsilon}) \| v \|_{0,h}, 
\end{equation}
Finally, taking into account \eqref{boundivPih}, the result follows. \QED \\

As a direct consequence of Propositions \ref{tildePih} and \ref{divPih} we have,
 
\begin{theorem}
\label{PhVh}
Provided $h$ is sufficiently small and $\Gamma$ is of the piecewise $C^{2k+3}$-class, there exists a constant $\delta > 0$ independent of $h$, such that,
\begin{equation}
\label{infsupRTc}
\forall v \in V_h \displaystyle \sup_{{\bf p} \in {\bf P}_{h} \mbox{ s.t. } [| {\bf p} |]_h=1} d_h({\bf p},v) \geq \delta \| v \|_{0,h}.
\end{equation}  
\end{theorem}
\prov
According to Proposition \ref{MixedPoisson}, for a given $v \in V_h$ $\exists \tilde{\bf p} \in {\bf W}$ such that
\begin{equation}
\label{equation1}
d_h(\tilde{\bf p},v) \geq \delta_{\bf W} \| \tilde{\bf p} \|_{\bf W} \| v \|_{0,h}.   
\end{equation}
Combining \eqref{errdivPih} and \eqref{equation1} we can write
\begin{equation}
\label{equation2}
d_h(\tilde{\Pi}_h \tilde{\bf p},v) \geq (\delta_{\bf W} - C_{\delta} h^{1/4}) \| \tilde{\bf p} \|_{\bf W} \| v \|_{0,h}.
\end{equation}
Assume that $h^{1/4} < \displaystyle \frac{\delta_{\bf W}}{2C_{\delta}}$. In this case \eqref{boundtildePih} allows us to write 
\begin{equation}
\label{equation3}
\displaystyle \frac{d_h(\tilde{\Pi}_h \tilde{\bf p},v)}{[| \tilde{\Pi}_h \tilde{\bf p} |]_h} \geq \displaystyle \frac{\delta_{\bf W} - C_{\delta} h^{1/4}}{C_{\tilde{\Pi}}} \| v \|_{0,h}.
\end{equation}
Thus defining ${\bf p}_1:= \displaystyle 
\frac{\tilde{\Pi}_h \tilde{\bf p}}{[| \tilde{\Pi}_h \tilde{\bf p} |]_h} \in {\bf P}_h$, it holds 
\begin{equation}
\label{equation4}
d_h({\bf p}_1,v) \geq \displaystyle \frac{\delta_{\bf W}}{2C_{\tilde{\Pi}}} \| v \|_{0,h},  
\end{equation}
with $[|{\bf p}_1|]_h=1$. Then the result follows with $\delta = \delta_{\bf W}/(2C_{\tilde{\Pi}})$. \QED \\

Finally, taking into account \eqref{normc}, we can prove the main result in this section, namely,

\begin{theorem}
\label{theo1}
Provided $h$ is sufficiently small and $\Omega$ is a convex domain smooth enough, (\ref{Dmixh}) has a unique solution. Moreover there exists a constant $\gamma > 0$ independent of $h$ such that,
\begin{equation}
\label{infsup}
\forall ({\bf p}; u) \in {\bf P}_h \times V_h \; \displaystyle \sup_{({\bf q};v) \in {\bf Q}_h \times V_h \mbox{ s.t. } ||| ({\bf q};v) |||_h = 1}
c_h(({\bf p}; u),({\bf q}; v)) \geq \gamma ||| ({\bf p}; u) |||_h. 
\end{equation}  
\end{theorem}

\prov 
According to \cite{AsyMVF} and owing to Theorems \ref{lemmars}, \ref{QhVh} and \ref{PhVh}, the bilinear form $c_h$ defined by \eqref{bilinearformch}-\eqref{ahbhdh} is weakly coercive over the pair of spaces ${\bf P}_h \times
 U_h$ and ${\bf Q}_h \times V_h$, with $U_h = V_h$. This allows us to assert that problem \eqref{Dmixh} has a unique solution (cf. \cite{COAM}). \\
Furthermore, again by \cite{AsyMVF} (cf. Theorem 3.1), the weak coercivity constant $\gamma$ for bilinear form $c_h$ given by \eqref{bilinearformch}-\eqref{ahbhdh} is also independent of $h$, since the three constants $\alpha$, $\beta$ and $\delta$ are independent of $h$. 
Actually, according to \cite{AsyMVF}, $\gamma$ may be expressed by
\begin{equation}
\label{gamma}
\left\{
\begin{array}{l}
\gamma:= \beta \nu/\sqrt{2[(\nu +1)^2+ \beta^2]} \\
\mbox{with} \\
\nu = \alpha \delta/\sqrt{(1+ \delta)^2  + \alpha^2}. \mbox{ \QED}
\end{array}
\right.
\end{equation}
 
\noindent \underline{\textbf{B) The non convex case}}\\
\indent At the price of more or less lengthier technicalities, we next extend the study of the uniform stability of problem \eqref{Dmixh} for a convex $\Omega$ to non convex domains.\\ 
\indent The non convex counterpart of condition \eqref{infsupRTc} will be established by exploiting the so-called extension technique.  
With this aim, we first observe that it does make sense to define an interpolate in ${\bf P}_h$ of a field $\tilde{\bf m}$ in a certain analog $\tilde{\bf W}$ of ${\bf W}$ with respect to $\tilde{\Omega}$, namely,\\
$\tilde{\bf W} := \{\tilde{\bf m}| \; \tilde{\bf m} \in {\bf H}(div,\tilde{\Omega}), \tilde{\bf m}_{|D} \in [H^{1/2+\varepsilon}(D)]^2, \mbox{ for } D=\Omega \mbox{ and } D = \Delta_{\Omega} \mbox{ and } \tilde{\bf m} \cdot {\bf n}_{|\Gamma_1}=0 \mbox{ on } \Gamma_1\}$. \\
We equip $\tilde{\bf W}$ with the norm $\| \cdot \|_{\tilde{\bf W}}$ defined as follows $\forall \tilde{\bf m} \in \tilde{\bf W}$:
\begin{equation}
\label{normtildeW}
\| \tilde{\bf m}  \|_{\tilde{\bf W}} := \| \tilde{\bf m} \|_{0,\tilde{\Omega}} + \| \nabla \cdot \tilde{\bf m} \|_{0,\tilde{\Omega}} + 
\| \tilde{\bf m}_{|\Omega} \|_{1/2+\varepsilon} + \| \tilde{\bf m}_{|\Delta_{\Omega}} \|_{1/2+\varepsilon,\Delta_{\Omega}}.
\end{equation}
We assume that $\forall {\bf m} \in {\bf W}$ there exists a suitable extension $\tilde{\bf m} \in \tilde{\bf W}$ of ${\bf m}$ to the whole $\tilde{\Omega}$ satisfying for a suitable constant $\tilde{C}_W$ 
\begin{equation}
\label{extendm}
\| \tilde{\bf m} \|_{\tilde{\bf W}} \leq \tilde{C}_{\bf W} \| {\bf m} \|_{\bf W} \; \forall {\bf m} \in {\bf W}.
\end{equation}
Once the extension $\tilde{\bf m} \in \tilde{\bf W}$ of ${\bf m}$ is fixed, the interpolate of $\tilde{\bf m}$ in ${\bf P}_h$ denoted by $\Pi_h^{'} \tilde{\bf m}$ is the field satisfying  
\begin{equation}
\label{Phkprime}
\left\{
\begin{array}{l}
\forall T \in {\mathcal T}_h \setminus {\mathcal S}_{1,h}, \;\int_e [\Pi_h^{'} \tilde{\bf m} - \tilde{\bf m}] \cdot {\bf n}_T \; v = 0 \; \forall v \in P_k(e) \; \forall \mbox{ edge } e \subset \partial T ; \\
\forall T \in {\mathcal S}_{1,h} \; \int_e [\Pi_h^{'} \tilde{\bf m} - {\bf m}] \cdot {\bf n}_T \; v = 0 \; \forall v \in P_k(e) \; \forall \mbox{ edge } e \subset \partial T \setminus e_T. \\
\mbox{For } k>0, \\
\forall T \in {\mathcal T}_h \setminus {\mathcal S}^{'}_{1,h} \; \int_T [\Pi_h^{'} \tilde{\bf m} - \tilde{\bf m}] \cdot {\bf o} = 0 \; \forall {\bf o} \in [P_{k-1}(T)]^2 ; \\
\forall T \in {\mathcal S}^{'}_{1,h} \; \int_{T^{{\Gamma}^{'}}} [\Pi_h^{'} \tilde{\bf m} - {\bf m}] \cdot {\bf o} = 0 \; \forall {\bf o} \in [P_{k-1}(T^{{\Gamma}^{'}})]^2.
\end{array}
\right.
\end{equation}
The existence and uniqueness of $\Pi_h^{'} \tilde{\bf m}$ relies basically upon the same arguments that apply to the interpolation operator $\Pi_h$ from ${\bf W}_h$ onto ${\bf Q}_h$, together with Proposition \ref{prop1}. In this respect, it is noticeable that for $T \in {\mathcal S}_{0,h}^{'}$ the intersection of $e_T$ with $\Gamma$ is necessarily the set of its end-points, by assumption. Therefore the condition on the normal trace of $\Pi_h^{'} \tilde{\bf m}$ over $e_T$ makes sense, since not only $\tilde{\bf m}$ belongs to $[L^s(\Delta_T)]^2$ for some $s>2$ (cf. \cite{BrezziFortin}), but also the embedding of $H^{1/2+\varepsilon}(\Delta_T)$ into $L^s(\Delta_T)$ is ensured by the Rellich-Kondratchev Theorem \cite{Adams}.  \\
The counterparts of \eqref{boundtildePih} and \eqref{errdivPih} for operator $\Pi_h^{'}$ will be established hereafter. 
Before pursuing however, we note that a condition analogous to \eqref{infsupOmega} in Proposition \ref{MixedPoisson} holds here, namely,

\begin{e-proposition}
\label{MixPoissonprim}
\begin{equation}
\label{LBBprime}
\left\{
\begin{array}{l}
\exists \delta_{\bf W}^{'} >0 \mbox{ such that }\forall {\mathcal T}_h \in {\mathcal F}^{'} \mbox{ and } \\
\forall v \in V_h \; \exists {\bf p} \in {\bf W} \setminus {\bf 0}_{{\bf W}} \mbox{ fulfilling } 
(\nabla \cdot {\bf p},\tilde{v}) \geq \delta_{\bf W}^{'} \| {\bf p} \|_{\bf W} \| \tilde{v} \|_{0},
\end{array}
\right.
\end{equation}
where $\tilde{v}$ is the function identical to $v$ in $\Omega_h$ and extended by zero in $\bar{\Omega} \setminus \bar{\Omega}_h$.
\end{e-proposition}

\prov 
First we address the case where $length(\Gamma_0)=0$. \\
We observe that, in principle $\tilde{v} \notin L^2_0(\Omega)$. Thus we define a function $v^{'} \in L^2_0(\Omega)$ by $v^{'} = 0$ in $\bar{\Omega} \setminus \bar{\Omega}_h$ and $v^{'}:=\tilde{v} + c^{'}$ in $\tilde{\Omega}_h$ where $c^{'}$ is a constant. Since $\int_{\Omega_h} \tilde{v}=0$, we trivially have
\begin{equation}
\label{cprime}
0 = c^{'} area(\Omega) + \int_{\Omega} \tilde{v} \implies c^{'} = \displaystyle \frac{\int_{\Omega_h \setminus \Omega} \tilde{v}}{area(\Omega)}.
\end{equation}
Notice that 
\begin{equation}
\label{tilde2}
\| v^{'} \|_{0}^2 \geq \| \tilde{v} \|_{0}^2 - 2 |c^{'}| \int_{\Omega} |\tilde{v}| + [c^{'}]^2 area(\Omega).
\end{equation}
On the other hand we have
\begin{equation}
\label{tilde3}
|c^{'}| \leq \displaystyle \sum_{T \in {\mathcal S}_h^{'}} \frac{\int_{\Delta_T} |\tilde{v}|}{area(\Omega)}.
\end{equation}
It is clear that
\begin{equation}
\label{tilde4}
\int_{\Delta_T} |\tilde{v}| \leq C_0 h_T^2 l_T \| \tilde{v} \|_{0,\infty,T} \; \forall T \in {\mathcal S}_h^{'}.
\end{equation}
Hence, using \eqref{L2TDelta} together with a straightforward argument, we come up with the constant $C_{\Delta}:=C_0{\mathcal C}_0^{+}$ such that 
\begin{equation}
\label{tilde5}
\int_{\Delta_T} |\tilde{v}| \leq C_{\Delta} h_T^{3/2} l_T^{1/2} \| \tilde{v} \|_{0,T^{{\Gamma}^{'}}} \; \forall T \in {\mathcal S}_h^{'}.
\end{equation}
Plugging \eqref{tilde5} into \eqref{tilde3} we readily obtain
\begin{equation}
\label{tilde6}
|c^{'}| \leq [area(\Omega)]^{-1} C_{\Delta} h^{3/2} \displaystyle \sqrt{ \sum_{T \in {\mathcal S}_h^{'}} \; l_T} \; 
\| \tilde{v} \|_{0} \leq [area(\Omega)]^{-1}  C_{\Delta} h^{3/2} [length(\Gamma)]^{1/2} \| \tilde{v} \|_{0}, 
\end{equation}
from which we infer the existence of a constant $C^{'}(\Omega)$ depending only on $\Omega$ such that
\begin{equation}
\label{tilde7}
|c^{'}| \leq C^{'}(\Omega) h^{3/2} \| \tilde{v} \|_{0}. 
\end{equation}
Plugging the above upper bound of $|c^{'}|$ into equation \eqref{tilde2} we successively obtain  
\begin{equation}
\label{tilde8}
\| v^{'} \|_{0}^2 \geq \| \tilde{v} \|_{0}^2 - 2 C^{'}(\Omega) h^{3/2} \| \tilde{v} \|_{0} \int_{\Omega} |\tilde{v}| + [C^{'}(\Omega)]^2 h^3 \| \tilde{v} \|_{0}^2  area(\Omega), 
\end{equation}
\begin{equation}
\label{tilde9}
\| v^{'} \|_{0}^2 \geq \| \tilde{v} \|_{0}^2 - 2 C^{'}(\Omega)  h^{3/2} \| \tilde{v} \|_{0}^2 \sqrt{area(\Omega)} + [C^{'}(\Omega)]^2 h^3 \| \tilde{v} \|_{0}^2  area(\Omega). 
\end{equation}
Hence setting $C(\Omega)=2C^{'}(\Omega) \sqrt{area(\Omega)}$ we conclude that 
\begin{equation}
\label{tilde}
\| v^{'} \|_{0}^2 \geq [1- C(\Omega) h^{3/2}] \| \tilde{v} \|_{0}^2.
\end{equation}
It follows from \eqref{tilde} that there exists a mesh-independent constant $C_{{\mathcal F}^{'}}$ such that $\forall {\mathcal T}_h \in {\mathcal F}^{'}$ it holds 
\begin{equation}
\label{calF}
\| v^{'} \|_{0} \geq C_{{\mathcal F}^{'}} \| \tilde{v} \|_{0} \; \forall v \in V_h.
\end{equation}
To complete this introductory part of the proof, we set $v{'} = \tilde{v}$ in $\tilde{\Omega}_h$ in case $length(\Gamma_0)>0$.\\

Now we consider a function $z^{'}$ defined in $\Omega$ to be the unique solution of the following Poisson problem:
\begin{equation}
\label{Poissonprime}
\left\{
\begin{array}{l}
- \Delta z^{'} = v^{'} \mbox{ in } \Omega \\
z^{'} = 0 \mbox{ on } \Gamma_0 \mbox{ if } length(\Gamma_0) >0 \mbox{ and } \int_{\Omega} z^{'} = 0 \mbox{ if } length(\Gamma_0) =0 \\
\partial z^{'}/\partial n = 0 \mbox{ on } \Gamma_1.
\end{array}
\right.
\end{equation}
Setting ${\bf p}^{'}:= -\nabla z^{'}$ we have ${\bf p}^{'} \in {\bf Q}$ and, according to \cite{Grisvard},   
${\bf p}^{'} \in [H^{1/2+\varepsilon}(\Omega)]^2$. Furthermore, there exists a constant $C^{'}_{\bf W}$ depending only on $\Omega$ and $\Gamma_1$ such that
\begin{equation}
\label{Continuityprime}
\| {\bf p}^{'} \|_{1/2+\varepsilon} \leq C^{'}_{\bf W} \| v^{'} \|_0.
\end{equation}
Since $(\nabla \cdot {\bf p}^{'},c^{'})=c^{'} \oint_{\Gamma} {\bf p}^{'} \cdot {\bf n} = 0$ if $length(\Gamma_0)=0$, using \eqref{calF} we obtain
\begin{equation}
\label{numerator}
(\nabla \cdot {\bf p}^{'},\tilde{v})=(\nabla \cdot {\bf p}^{'},v^{'}) = (-\Delta z^{'},v^{'}) = \| v^{'} \|_{0}^2 \geq C_{{\mathcal F}^{'}}^2 \| \tilde{v} \|_0^2. 
\end{equation}
Thus, choosing ${\bf p}={\bf p}^{'}$ and taking into account \eqref{Continuityprime}, \eqref{LBBprime} is trivially seen to hold with $\delta_{\bf W}^{'} = C_{{\mathcal F}^{'}}^2/(1+C^{'}_{\bf W})$.
\QED \\

Next we prove two propositions for operator $\Pi_h^{'}$ analogous to Propositions \ref{tildePih} and \ref{divPih} for $\Pi_h$. 
\begin{e-proposition}
\label{Pihprime}
There exists a constant $C_{\Pi}^{'}$ independent of $h$ such that
\begin{equation}
\label{boundPihprime}
[| \Pi_h^{'} \tilde{\bf m} |]_{h} \leq C_{\Pi}^{'} \| \tilde{\bf m} \|_{\tilde{\bf W}} \; \forall \tilde{\bf m} \in \tilde{\bf W}.
\end{equation}
\end{e-proposition}

\prov
First of all, in regard of \eqref{Phkprime}, it is not difficult to figure out that $\Pi_h^{'} \tilde{\bf m} \equiv {\bf m}$ in $T^{\Gamma{'}}$ $\forall {\bf m} \in [P_k(T^{{\Gamma}^{'}})]^2$. In addition to this the interpolate $\Pi_h^{'} \tilde{\bf m}$ restricted to $T^{{\Gamma}^{'}}$ fulfills all the other conditions specified in classical interpolation theory (cf. \cite{BrennerScott,Arcangeli}) to satisfy for a constant $C^{'}_{\varepsilon}$ independent of $h$
\begin{equation}
\label{L2Pihprime}
\forall q \in [0,1/2+\varepsilon] \; \| \Pi_h^{'} \tilde{\bf m} - {\bf m} \|_{q,T^{\Gamma^{'}}} \leq C^{'}_{\varepsilon} h^{1/2+\varepsilon-q} \| {\bf m} \|_{1/2+\varepsilon,T^{\Gamma^{'}}} \; \forall T \in {\mathcal T}_h \mbox{ and } \forall {\bf m} \in {\bf W}.
\end{equation}
It trivially follows from \eqref{L2Pihprime} that 
\begin{equation}
\label{boundPihprime1}
\| \Pi_h^{'} \tilde{\bf m} \|^{'}_{0,h} \leq (1+C^{'}_{\varepsilon} h^{1/2+\varepsilon}) \| {\bf m} \|_{1/2+\varepsilon} \; \forall {\bf m} \in {\bf W}.
\end{equation}
On the other hand, recalling \eqref{L2TDelta}, we have 
\begin{equation}
\label{boundPihprime2}
\| \Pi_h^{'} \tilde{\bf m} \|_{0,T}^2 \leq ({\mathcal C}_{0}^{-})^{-2}({\mathcal C}_{0}^{+})^2 \| \Pi_h^{'} \tilde{\bf m} \|_{0,T^{\Gamma^{'}}}^2 \; \forall T \in {\mathcal T}_h.
\end{equation}
Therefore, summing up both sides of \eqref{boundPihprime2} over the mesh, taking into account \eqref{boundPihprime1} and assuming that $h<1$, the following bound is seen to hold with $C^{0}_{\Pi^{'}}=(1+C^{'}_{\varepsilon}){\mathcal C}_{0}^{+}/{\mathcal C}_{0}^{-}$:
\begin{equation}
\label{boundPihprime3}
\| \Pi_h^{'} \tilde{\bf m} \|_{0,h} \leq C^{0}_{\Pi^{'}} \| {\bf m} \|_{\bf W} \; \forall {\bf m} \in {\bf W}.  
\end{equation}

As for the divergence term, first we note that the case of $\| \nabla \cdot \Pi_h^{'} \tilde{\bf m} \|_{0,T}$ for $T \in {\mathcal T}_h \setminus {\mathcal S}_{h}^{'}$ was already addressed in Section 5. Indeed, from \eqref{divPihcalTh} together with \eqref{boundivPih5} and 
recalling \eqref{boundivPih}, we easily infer that 
\begin{equation}
\label{partdivPih} 
\displaystyle \sum_{T \in {\mathcal T}_h \setminus {\mathcal S}^{'}_{h}} \| \nabla \cdot \Pi^{'}_h \tilde{\bf m} \|_{0,T}^2 \leq  (C^D_{\tilde{\Pi}})^2 \| {\bf m} \|_{\bf W}^2 \; \forall {\bf m} \in {\bf W}.
\end{equation}
Now for triangles in ${\mathcal S}_{h}^{'}$ we first have
\begin{equation}
\label{S1hprime}
\left\{
\begin{array}{l}
\forall T \in {\mathcal S}_{1,h}^{'} \\
\| \nabla \cdot \Pi_h^{'} \tilde{\bf m} \|_{0,T}^2 = J^1_{1,T}({\bf m}) + J^1_{2,T}({\bf m}) + J^1_{3,T}({\bf m}) \\
\mbox{with} \\
J^1_{1,T}({\bf m}) := \int_{T^{\Gamma^{'}}} \nabla \cdot \Pi_h^{'} \tilde{\bf m} (\nabla \cdot \Pi_h^{'} \tilde{\bf m} -
\nabla \cdot {\bf m} );\\
J^1_{2,T}({\bf m}) := \| \nabla \cdot \Pi_h^{'} \tilde{\bf m} \|_{0,\Delta_T}^2; \\
J^1_{3,T}({\bf m}) := \int_{T^{\Gamma^{'}}} \nabla \cdot \Pi_h^{'} \tilde{\bf m} \nabla \cdot {\bf m}.
\end{array}
\right.
\end{equation}
From the definition of the operator $\Pi_h^{'}$, integrating by parts we obtain after straightforward calculations
\begin{equation}
\label{J11h1}
J^1_{1,T}({\bf m}) = (\Pi_h^{'} \tilde{\bf m} \cdot {\bf n},\nabla \cdot \Pi_h^{'} \tilde{\bf m})_{0,\Gamma_T}.
\end{equation}
Then Lemma \ref{GammaT} yields 
\begin{equation}
\label{errK11prime}
|J^1_{1,T}({\bf m})| \leq  C^{'}_{\Gamma_1} h_T^{k+1} \| \Pi^{'}_h \tilde{\bf m} \|_{k+1,T^{\Gamma^{'}}} \| \nabla \cdot \Pi^{'}_h \tilde{\bf m} \|_{0,T}, 
\end{equation}
from which we obtain in turn for a typical inverse inequality constant $C^{'}_{1/2}$ (cf. \cite{Georgoulis})
\begin{equation}
\label{J11h2}
|J^1_{1,T}({\bf m})| \leq  C_{1/2}^{'} h_T^{1/2} \| \Pi^{'}_h \tilde{\bf m} \|_{1/2+\varepsilon,T^{\Gamma^{'}}} \| \nabla \cdot \Pi^{'}_h \tilde{\bf m} \|_{0,T}.
\end{equation} 
Further using \eqref{L2Pihprime} with $q = 1/2+\varepsilon$ and setting $ C_{J,1}^{1}=C_{1/2}^{'}(1+C^{'}_{\varepsilon})$, from \eqref{J11h2} we infer that
\begin{equation}
\label{J11h3}
|J^1_{1,T}({\bf m})| \leq C_{J,1}^{1} h_T^{1/2} \| {\bf m} \|_{1/2+\varepsilon,T^{\Gamma^{'}}} \| \nabla \cdot \Pi^{'}_h \tilde{\bf m} \|_{0,T}.
\end{equation}
Then using Young's inequality we come up with
\begin{equation}
\label{J11h}
|J^1_{1,T}({\bf m})| \leq (C_{J,1}^{1})^2 h_T \| \nabla \cdot \Pi^{'}_h \tilde{\bf m} \|_{0,T}^2/2 + \| {\bf m} \|_{1/2+\varepsilon,T^{\Gamma^{'}}}^2/2.
\end{equation}
As for $J^1_{2,T}({\bf m})$, by a classical inverse inequality, we obtain for a suitable constant $C^{1}_{J,2}$ 
\begin{equation}
\label{J12h}
|J^1_{2,T}({\bf m})| \leq C_0 h_T^3 \| \nabla \cdot \Pi_h^{'} \tilde{\bf m} \|_{0,\infty,T}^2 \leq C^{1}_{J,2} h_T 
\| \nabla \cdot \Pi_h^{'} \tilde{\bf m} \|_{0,T}^2.
\end{equation}
Finally, using Young's inequality, we easily obtain for $J^1_{3,T}({\bf m})$
\begin{equation}
\label{J13h}
|J^1_{3,T}({\bf m})| \leq \| \nabla \cdot \Pi_h^{'} \tilde{\bf m} \|_{0,T}^2/2 + \| \nabla \cdot  {\bf m} \|_{0,T^{\Gamma^{'}}}^2/2 \leq 
\| \nabla \cdot \Pi_h^{'} \tilde{\bf m} \|_{0,T}^2/2 + \| \nabla \cdot  \tilde{\bf m} \|_{0,T}^2/2. \\
\end{equation}
Collecting \eqref{J11h}, \eqref{J12h} and \eqref{J13h}, \eqref{S1hprime} is readily seen to yield
\begin{equation}
\label{J1h}
\forall T \in {\mathcal S}_{1,h}^{'} \;
\| \nabla \cdot \Pi_h^{'} \tilde{\bf m} \|_{0,T}^2 \leq \{1 - 2h_T[C^{1}_{J,2}+(C_{J,1}^{1})^2]\}^{-1} ( \| {\bf m} \|_{1/2+\varepsilon,T^{\Gamma^{'}}}^2 + \| \nabla \cdot  \tilde{\bf m} \|_{0,T}^2).
\end{equation}
Thus, as long as $h$ is sufficiently small, there exists a mesh-independent constant $C_J^1$ such that
\begin{equation}
\label{divPihprime1}
\forall T \in {\mathcal S}_{1,h}^{'} \;
\| \nabla \cdot \Pi_h^{'} \tilde{\bf m} \|_{0,T}^2 \leq C_J^1 ( \| {\bf m} \|_{1/2+\varepsilon,T^{\Gamma^{'}}}^2 +  
\| \nabla \cdot  \tilde{\bf m} \|_{0,T}^2).
\end{equation}
Next we consider triangles in ${\mathcal S}_{0,h}^{'}$. We have 
\begin{equation}
\label{S0hprime}
\left\{
\begin{array}{l}
\forall T \in {\mathcal S}_{0,h}^{'} \\
\| \nabla \cdot \Pi_h^{'} \tilde{\bf m} \|_{0,T}^2 = J^0_{1,T}({\bf m}) + J^0_{2,T}({\bf m}) \\
\mbox{with} \\
J^0_{1,T}({\bf m}) := \int_T (\nabla \cdot \Pi_h^{'} \tilde{\bf m} -
\nabla \cdot \tilde{\bf m} )\nabla \cdot \Pi_h^{'} \tilde{\bf m};\\
J^0_{2,T}({\bf m}) := \int_T \nabla \cdot \tilde{\bf m} \nabla \cdot \Pi_h^{'} \tilde{\bf m}.
\end{array}
\right.
\end{equation}
Taking into account the definition $\Pi_h^{'}$ we can assert that
\begin{equation}
\label{J01h}
J^0_{1,T}({\bf m}) = \int_{\partial T} (\Pi_h^{'} \tilde{\bf m} -
\tilde{\bf m} ) \cdot {\bf n}_T \nabla \cdot \Pi_h^{'} \tilde{\bf m} - \int_T (\Pi_h^{'} \tilde{\bf m} -
\tilde{\bf m}) \cdot \nabla \nabla \cdot \Pi_h^{'} \tilde{\bf m} = 0. 
\end{equation}
On the other hand, we have
\begin{equation}
\label{J02h}
|J^0_{2,T}({\bf m})| \leq  \| \nabla \cdot \Pi_h^{'} \tilde{\bf m} \|_{0,T}^2/2 + \| \nabla \cdot \tilde{\bf m} \|_{0,T}/2 
\end{equation}
Combining \eqref{S0hprime}, \eqref{J01h} and \eqref{J02h} we obtain
\begin{equation}
\label{divPihprime0}
\forall T \in {\mathcal S}_{0,h}^{'} \; \| \nabla \cdot \Pi_h^{'} \tilde{\bf m} \|_{0,T} \leq \| \nabla \cdot \tilde{\bf m} \|_{0,T}
\end{equation}
Finally, as an immediate consequence of \eqref{partdivPih}, \eqref{divPihprime1} and \eqref{divPihprime0}, it holds for a suitable mesh-  independent constant $C^D_{\Pi^{'}}$ 
\begin{equation}
\label{boundivPihprime}
\| \nabla \cdot \Pi_h^{'} \tilde{\bf m} \|_{0,h} \leq C^D_{\Pi^{'}} \| \tilde{\bf m} \|_{\tilde{\bf W}} \; \forall {\bf m} \in {\bf W}.
\end{equation}
The addition of \eqref{boundivPihprime} and \eqref{boundPihprime3} yields \eqref{boundPihprime}
with $C_{\Pi}^{'} = \sqrt{(C^{0}_{\Pi^{'}})^2 + (C^D_{\Pi^{'}})^2}$. \QED \\
 
\begin{e-proposition}
\label{divPihprime}
There exists a constant $C_{\delta}^{'}$ independent of $h$ such that $\forall {\bf m} \in {\bf W}$ and $\forall v \in V_h$
\begin{equation}
\label{errdivPihprime}
d_h(\Pi_h^{'} \tilde{\bf m},v) \geq d_h(\tilde{\bf m},v) 
- \displaystyle \sum_{T \in {\mathcal S}_{1,h}^{'}} |(\nabla \cdot \tilde{\bf m},v)_{\Delta_T}|  
- C^{'}_{\delta} h^{1/4} \| {\bf m} \|_{{\bf W}} \| v \|_{0,h}.
\end{equation}
\end{e-proposition}

\prov
Using the extension $\tilde{\bf m} \in \tilde{\bf W}$ of ${\bf m}$ satisfying \eqref{extendm}, first we note that by construction we have $\forall {\bf m} \in {\bf W}$ and $\forall v \in V_h$
\begin{equation}
\label{TnonS1}
\displaystyle \sum_{T \in {\mathcal T}_h \setminus {\mathcal S}_{1,h} } (\nabla \cdot {\Pi}^{'}_h \tilde{\bf m},v)_{T} = 
\displaystyle \sum_{T \in {\mathcal T}_h \setminus {\mathcal S}_{1,h} } (\nabla \cdot \tilde{\bf m},v)_{T}
\end{equation}
Furthermore, on the basis of Proposition \ref{divPih}, we may write $\forall {\bf m} \in {\bf W}$ and $\forall v \in V_h$
\begin{equation}
\label{S1nonprime}
\displaystyle \sum_{T \in {\mathcal S}_{1,h} \setminus {\mathcal S}_{1,h}^{'}} (\nabla \cdot {\Pi}^{'}_h \tilde{\bf m},v)_{T} \geq 
\displaystyle \sum_{T \in {\mathcal S}_{1,h} \setminus {\mathcal S}_{1,h}^{'}} (\nabla \cdot \tilde{\bf m},v)_{T} - C_{\delta} h^{1/4} \| {\bf m} \|_{\bf W} \| v \|_{0,h} \; 
\end{equation}
As for $T \in {\mathcal S}_{1,h}^{'}$, we may write $\forall {\bf m} \in {\bf W}$ and $\forall v \in V_h$
\begin{equation}
\label{prime1}
\left\{
\begin{array}{l}
\forall T \in {\mathcal S}_{1,h}^{'} \; 
(\nabla \cdot \Pi_h^{'} \tilde{\bf m},v)_{T} = K^{'}_{1,T}({\bf m},v)+K^{'}_{2,T}({\bf m},v)+K^{'}_{3,T}({\bf m},v) \\
\mbox{with} \\
K^{'}_{1,T}({\bf m},v) = (\nabla \cdot \Pi_h^{'} \tilde{\bf m}, v)_{\Delta_T} \\
K^{'}_{2,T}({\bf m},v) = (\nabla \cdot [\Pi_h^{'} \tilde{\bf m}-{\bf m}], v)_{T^{\Gamma^{'}}} \\
K^{'}_{3,T}({\bf m},v) = (\nabla \cdot {\bf m}, v)_{T^{\Gamma^{'}}}.
\end{array}
\right.
\end{equation}
Using arguments repeatedly exploited in this work we may write for a constant $C_{K,1}^{'}$ independent of $T$
\begin{equation}
\label{K11}
|K^{'}_{1,T}({\bf m},v)| \leq C_{K,1}^{'} h_T \|\nabla \cdot \Pi_h^{'} \tilde{\bf m}\|_{0,T} \| v \|_{0,T} \; \forall T \in {\mathcal S}^{'}_{1,h}
\end{equation}
Moreover by construction we have $K^{'}_{2,T}({\bf m},v)| = (\Pi^{'}_h \tilde{\bf m},v)_{\Gamma_T}$. Thus Lemma \ref{GammaT} and \eqref{L2TDelta} yield 
\begin{equation}
\label{errgvprime}
|K^{'}_{2,T}({\bf m},v)| = \displaystyle \left| \int_{\Gamma_T} \Pi^{'}_h \tilde{\bf m} v \; ds \right| \leq {\mathcal C}_{\Gamma_1}^{'} h_T^{k+1} \| \Pi_h^{'} \tilde{\bf m} \|_{k+1,T^{\Gamma^{'}}} \| v \|_{0,T^{\Gamma^{'}}} \; \forall T \in {\mathcal S}_{1,h}^{'},
\end{equation}
where ${\mathcal C}_{\Gamma_1}^{'}$ is a constant independent of $T$. Then  mimicking \eqref{J11h2} and taking $q=1/2$ in \eqref{L2Pihprime}, we can assert the existence of another mesh-independent constant $C_{K,2}^{'}$ such that
\begin{equation}
\label{K12}
|K^{'}_{2,T}({\bf m},v)| \leq C_{K,2}^{'} h_T^{1/2} \| {\bf m}\|_{1/2+\varepsilon,T^{\Gamma^{'}}} \| v \|_{0,T} \; \forall 
T \in {\mathcal S}_{1,h}^{'}.
\end{equation}
Putting together \eqref{K11}, \eqref{K12} and \eqref{prime1}, $\forall {\bf m} \in {\bf W}$ and $\forall v \in V_h$ we have
\begin{equation}
\label{S1prime}
\displaystyle \sum_{T \in {\mathcal S}_{1,h}^{'}}
(\nabla \cdot \Pi_h^{'} \tilde{\bf m},v)_{T} \geq \displaystyle \sum_{T \in {\mathcal S}_{1,h}^{'}} \left[ (\nabla \cdot {\bf m},v)_{T^{\Gamma^{'}}} - (C_{K,1}^{'} h_T \| \nabla \cdot \Pi_h^{'} \tilde{\bf m} \|_{0,T} + C_{K,2}^{'} h_T^{1/2} \|{\bf m} \|_{1/2+\varepsilon,T^{\Gamma^{'}}}) \| v \|_{0,T} \right]
\end{equation}  
Thus, taking into account \eqref{TnonS1}, \eqref{S1nonprime} and \eqref{S1prime} and setting $C^{'}_K:= C_{K,2}^{'} + C_{\delta}$, since $h<1$ we obtain $\forall {\bf m} \in {\bf W}$ and $\forall v \in V_h$
\begin{equation}
\label{Pihprimev}
\left\{
\begin{array}{l}
d_h(\Pi_h^{'} \tilde{\bf m},v) \geq d_h(\tilde{\bf m},v) 
- \displaystyle \sum_{T \in {\mathcal S}_{1,h}^{'}} |(\nabla \cdot \tilde{\bf m},v)_{\Delta_T}| \\ 
- \displaystyle \left[ C_{K}^{'} h^{1/4} \| {\bf m} \|_{\bf W}  + C^{'}_{K,1} h \| \nabla \cdot \Pi_h^{'} \tilde{\bf m} \|_{0,h} \right] \| v \|_{0,h}.
\end{array}
\right.
\end{equation}
Finally, recalling \eqref{boundivPihprime} and setting $C^{'}_{\delta} = C_{K}^{'}+\tilde{C}_{\bf W} C_{\Pi^{'}}^D C^{'}_{K,1}$, \eqref{errdivPihprime} readily derives from \eqref{Pihprimev}. \QED \\ 

Propositions \ref{MixPoissonprim}, \ref{Pihprime} and \ref{divPihprime} are all that is needed to establish the validity of an \textit{inf-sup} condition analogous to \eqref{infsupRTc} in the non convex case. More precisely we have 
 
\begin{theorem}
\label{PhVhprime}
If $\Gamma$ is of the piecewise $C^{2k+3}$-class, there exists a constant $\delta^{'}>0$ independent of $h$, such that $\forall {\mathcal T}_h \in {\mathcal F}^{'}$ it holds
\begin{equation}
\label{infsupRTch}
\forall v \in V_h \displaystyle \sup_{{\bf p} \in {\bf P}_{k,h} \mbox{ s.t. } [| {\bf p} |]_h=1} d_h({\bf p},v)_h \geq  \delta^{'} \| v \|_{0,h}.
\end{equation}
\end{theorem}

\prov
Recalling Proposition \ref{MixPoissonprim}, given $v \in V_h$, let ${\bf p}$ be the field in ${\bf W}$ fulfilling \eqref{LBBprime}. The first step to prove \eqref{infsupRTch} is to construct a suitable extension $\tilde{\bf p} \in \tilde{\bf W}$ of ${\bf p}$ to $\tilde{\Omega}$ satisfying \eqref{extendm}. With this aim we pose the following Poisson problem in the outer domain $\Delta_{\Omega}$, in terms of the solution $z^{'}$ of \eqref{Poissonprime}: 
\begin{equation}
\label{Poissonout}
\left\{
\begin{array}{l}
\mbox{Find } z_{\Delta} \mbox{ such that}\\
- \Delta z_{\Delta} = 0 \mbox{ in } \Delta_{\Omega} \\
z_{\Delta} =0 \mbox{ on } \tilde{\Gamma} \\
-\partial z_{\Delta}/\partial n = \partial z^{'}/\partial n \mbox{ on } \Gamma,
\end{array}
\right.
\end{equation}
. \\
Now we set $\tilde{\bf p} = {\bf p} (= -\nabla z^{'}$) in $\Omega$ and $\tilde{\bf p}= -\nabla z_{\Delta}$ in $\Delta_{\Omega}$. It is clear that $\tilde{\bf p} \in {\bf H}(div,\tilde{\Omega})$, owing to the continuity of the normal trace of $\tilde{\bf p}$ across $\Gamma$ by construction.
Moreover, since $\partial z^{'}/\partial n \in H^{\varepsilon}(\Gamma)$, we have $z_{\Delta} \in H^{3/2+\varepsilon}(\Delta_{\Omega})$, and 
the harmonic function $z_{\Delta}$ fulfills for a constant $C_{out}$
\begin{equation}
\label{Cout}
\| z_{\Delta} \|_{3/2+\varepsilon,\Delta_{\Omega}} \leq C_{out} \|\partial z^{'}/\partial n \|_{\varepsilon,\Gamma}.
\end{equation}
Thus, by the Trace Theorem in $\Omega$, there is another constant $C_{in}$ such that 
\begin{equation}
\label{Cin}
\| z_{\Delta} \|_{3/2+\varepsilon,\Delta_{\Omega}} \leq C_{in} \| z^{'} \|_{3/2+\varepsilon}.
\end{equation}
Using classical Friedrichs-Poincar\'e inequalities (cf. \cite{DuvautLions}), we can assert the existence of a constant $C_{\Delta}$ depending only on $\Omega$ and $\Delta_{\Omega}$ such that
\begin{equation}
\label{poutin}
\| \tilde{\bf p} \|_{1/2+\varepsilon,\Delta_{\Omega}} \leq C_{\Delta} \| {\bf p}^{'} \|_{1/2+\varepsilon}.
\end{equation}
Since $\nabla \cdot \tilde{\bf p}=0$ in $\Delta_{\Omega}$, the norm of $\tilde{\bf p}$ in ${\bf H}(div,\Delta_{\Omega})$ is bounded above by $\| \tilde{\bf p} \|_{1/2+\varepsilon,\Delta_{\Omega}}$. Hence for $\tilde{C}_{\bf W}:=(1+C_{\Delta}^2)^{1/2}$ it holds
\begin{equation}
\label{boundpDelta}
\| \tilde{\bf p} \|_{\tilde{\bf W}} \leq  \tilde{C}_{\bf W} \| {\bf p} \|_{\bf W}.\\
\end{equation}
Now for an arbitrary function in $v \in V_h$, ${\bf p}$ being the field in ${\bf W}$ fulfilling \eqref{LBBprime} and $\tilde{\bf p} \in \tilde{\bf W}$ being the extension of ${\bf p}$ fulfilling \eqref{poutin}, we have
\begin{equation}
\label{LBBprimeh1}
d_h( \Pi_h^{'} \tilde{\bf p},v) = d_h(\Pi_h^{'} \tilde{\bf p}-\tilde{\bf p},v) +\displaystyle \sum_{T \in {\mathcal S}_{h}^{'}} (\nabla \cdot \tilde{\bf p},v)_{\Delta_T} + (\nabla \cdot {\bf p},\tilde{v}).
\end{equation} 
Since $\nabla \cdot \tilde{\bf p} \equiv 0$ in $\Delta_{\Omega}$, and thus in $\Delta_T$ too $\forall T \in {\mathcal S}^{'}_h$, taking into account \eqref{errdivPihprime}, \eqref{LBBprimeh1} yields
\begin{equation}
\label{LBBprimeh2}
d_h(\Pi_h^{'} \tilde{\bf p},v) \geq (\nabla \cdot {\bf p},\tilde{v}) 
- C^{'}_{\delta} h^{1/4} \| \tilde{\bf p} \|_{\tilde{\bf W}} \| v \|_{0,h}.
\end{equation} 
Furthermore, using \eqref{LBBprime} and \eqref{boundpDelta}, from \eqref{LBBprimeh2} we obtain
\begin{equation}
\label{LBBprimeh3}
d_h(\Pi_h^{'} \tilde{\bf p},v) \geq \delta_{\bf W}^{'} \| {\bf p} \|_{\bf W} \| \tilde{v} \|_{0}
- C^{'}_{\delta} \tilde{C}_W h^{1/4} \| {\bf p} \|_{\bf W} \| v \|_{0,h}.
\end{equation}
On the other hand, by a straightforward argument similar to those already exploited in this work, we can assert that, provided $h$ is small  enough, there exists a mesh-independent constant $C_{V}^{'}$ such that
\begin{equation}
\label{CVprime}
\| \tilde{v} \|^2_{0} = \| v \|_{0,h}^2 - \displaystyle \sum_{T \in {\mathcal S}^{'}_h} \int_{\Delta_T} v^2 \geq \| v \|_{0,h}^2 -
C_0 h_T^3 \displaystyle \sum_{T \in {\mathcal S}^{'}_h} \| v \|_{0,\infty,T}^2 \geq (C_{V}^{'})^2 \| v \|_{0,h}^2.
\end{equation}
Plugging \eqref{CVprime} into \eqref{LBBprimeh3}, it follows that 
\begin{equation}
\label{LBBprimeh4}
d_h(\Pi_h^{'} \tilde{\bf p},v) \geq (\delta_{\bf W}^{'} C_V^{'}  
- C^{'}_{\delta} \tilde{C}_W h^{1/4}) \| {\bf p} \|_{\bf W} \| v \|_{0,h}.
\end{equation}
On the other hand, recalling \eqref{boundPihprime} and \eqref{boundpDelta} we have
\begin{equation}
\label{Pihprimetildep}
[| \Pi_h^{'} \tilde{\bf p} |]_{h} \leq C_{\Pi}^{'} \| \tilde{\bf p} \|_{\tilde{\bf W}} \leq C_{\Pi}^{'} \tilde{C}_W \| {\bf p} \|_{\bf W},
\end{equation}
and hence from \eqref{LBBprimeh4} we obtain
\begin{equation}
\label{LBBprimeh5}
d_h(\Pi_h^{'} \tilde{\bf p},v) \geq \displaystyle \frac{\delta_{\bf W}^{'} C_V^{'}  
- C^{'}_{\delta} \tilde{C}_W h^{1/4}}{C_{\Pi}^{'} \tilde{C}_W} [| \Pi_h^{'} \tilde{\bf p} |]_{h} \| v \|_{0,h}.
\end{equation}
Thus, provided $h$ is small enough and in particular $h^{1/4} \leq \displaystyle \frac{\delta_{\bf W}^{'} C_V^{'}}{  
2 C^{'}_{\delta} \tilde{C}_W}$, \eqref{infsupRTch} immediate follows from \eqref{LBBprimeh5} with $\delta^{'} = \displaystyle  
\delta_{\bf W}^{'} C_V^{'}/\left(2 C_{\Pi}^{'} \tilde{C}_W \right)$. 
\QED\\
To conclude we have
 
\begin{theorem}
\label{nonconvexstab}
Assume that $\Omega$ is a non convex domain of the piecewise $C^{2k+3}$-class. Provided $h$ is sufficiently small problem \eqref{Dmixh} has a unique solution. Moreover we can assert that a uniform stability condition of the same type as \eqref{infsup} holds true for this problem.
\end{theorem}

\prov
Recalling that the validity of both \eqref{infsupRT} and \eqref{infsupRSbis} is guaranteed, irrespective of the fact that $\Omega$ be convex or not, from \eqref{infsupRTch} we can assert that the approximate problem \eqref{Dmixh} has a unique solution. Furthermore recalling \eqref{gamma}, a uniform stability condition of the type \eqref{infsup} does hold with a constant $\gamma^{'}$ instead of $\gamma$, given by 
\begin{equation}
\label{gammaprime}
\left\{
\begin{array}{l}
\gamma^{'}:= \beta \nu^{'}/\sqrt{2[(\nu^{'} +1)^2+ \beta^2]} \\
\mbox{with} \\
\nu^{'} = \alpha \delta^{'}/\sqrt{(1+ \delta^{'})^2  + \alpha^2}. \mbox{ \QED}
\end{array}
\right.
\end{equation}

\section{Interpolation errors} 

\hspace{4mm} Generalizing classical estimates for interpolation errors applying to the space ${\bf Q}_h$ (cf. \cite{RaviartThomas}), we 
establish here new ones for the space ${\bf P}_h$. \\

In this section the notations $D^1w$ and $D^2w$ will eventually be replaced by $\nabla w$ and $H(w)$ for the gradient and the Hessian of $w$, respectively. Furthermore in the sequel several arguments will rely upon classical inverse inequalities for polynomials defined in $T \in {\mathcal T}_h$ (cf. \cite{Verfuerth,ErnGuermond}) and their extensions to neighboring sets. Besides \eqref{LinftyTDelta} and \eqref{L2TDelta} we shall use the following one: \\
There exists a constant ${\mathcal C}_I$ depending only on 
$k$ and the shape regularity of ${\mathcal T}_h$ (cf. \cite{BrennerScott}, Ch.4, Sect 4) such that for any 
set $T^{'}$ satisfying $T^{\Gamma^{'}} \subseteq T^{'} \subseteq T^{\Gamma}$, it holds for $1 \leq j \leq k$:  
\begin{equation}
\label{inverse}
\parallel D^j w \parallel_{0,T^{'}} \leq {\mathcal C}_I h_T^{-1} \parallel D^{j-1} w \parallel_{0,T^{'}} \; 
\forall w \in {\mathcal P}_k(T^{'}).
\end{equation}
\indent Now let $\varepsilon > 0$ be a small real number. We define
$${\mathcal M}^{k}(\Omega):= \{{\bf m}| \;{\bf m} \in [H^{k+1+\varepsilon(1-\min[k,1])}(\Omega)]^2, \; \nabla \cdot {\bf m} \in H^{k+1}(\Omega), \mbox{ and } [{\bf m} \cdot {\bf n}]_{|\Gamma_1} = 0 \}.$$
We shall also use extensions of fields in ${\mathcal M}^k(\Omega)$ to a domain $\tilde{\Omega}$ containing both $\Omega$ and all the sets $\Omega_h$ for ${\mathcal T}_h \in {\mathcal F}^{'}$. More precisely, here we consider $\tilde{\Omega}$ to be a domain whose boundary $\tilde{\Gamma}$ is such that the minimum of the distance from all of its points to $\Gamma$ is not less than $h_{max}^2$, where $h_{max} = 2 \displaystyle \max_{{\mathcal T}_h \in {\mathcal F}^{'}}[h]$. Let $${\mathcal M}^k(\tilde{\Omega}):=\{\tilde{\bf m}| \;\tilde{\bf m} \in [H^{k+1+\varepsilon(1-\min[k,1])}(\tilde{\Omega})]^2, \; \nabla \cdot \tilde{\bf m} \in H^{k+1}(\tilde{\Omega}), \mbox{ and } [\tilde{\bf m} \cdot {\bf n}]_{|\Gamma_1} = 0 \}.$$
With every ${\bf m} \in {\mathcal M}^k(\Omega)$ we associate its extension $\tilde{\bf m}$ to $\tilde{\Omega} \setminus \Omega$ belonging to ${\mathcal M}^k(\tilde{\Omega})$ constructed, for instance, as advocated by Stein \& al. \cite{Stein} (i.e. using polyharmonic functions).\\  
Now, since the restriction to $\Omega$ of any field in ${\mathcal M}^k(\tilde{\Omega})$ necessarily belongs to ${\bf W}$ for all $k \geq 0$, given $\tilde{\bf m}$ in ${\mathcal M}^k(\tilde{\Omega})$ we may define its interpolate in ${\bf P}_h$ inside $\tilde{\Omega}_h$ in the same manner as \eqref{Phkinterpolate} for fields in ${\bf W}$ and a convex $\Omega$. In short, we still denote such an interpolation operator from ${\mathcal M}^{k}(\tilde{\Omega})$ onto ${\bf P}_h$ by $\tilde{\Pi}_h$, which is constructed as prescribed in \eqref{Phkinterpolate}.\\
Notice that here we have $\tilde{\Pi}_h \tilde{\bf m} \in {\bf H}(div,\tilde{\Omega}_h)$. 
It is also noteworthy that every $\tilde{\bf m} \in {\mathcal M}^k(\tilde{\Omega})$ 
belongs to $[H^{2}(\tilde{\Omega})]^2$ for every $k \geq 1$ and hence, from the Sobolev Embedding Theorem, it makes sense to prescribe $[\tilde{\Pi}_h \tilde{\bf m} \cdot {\bf n}](P) = [\tilde{\bf m} \cdot {\bf n}](P)$ at any point $P \in \Gamma$ (cf. \cite{Adams}). Notice that, in case $k=0$, such a prescription is also admissible, thanks to the assumption that ${\bf m} \in [H^{1+\varepsilon}(\Omega)]^2$. 
Moreover, since $\forall \tilde{\bf m} \in {\mathcal M}^k(\tilde{\Omega})$ and $\forall k \geq 1$ we have $\nabla \cdot \tilde{\bf m} \in H^{2}(\Omega)$, for the same reason it is legitimate to associate with $\nabla \cdot \tilde{\bf m}$ a function $d_h(\tilde{\bf m}) \in L^2(\Omega_h)$, whose restriction to every $T \in {\mathcal T}_h$ 
is the standard Lagrange interpolate of $[\nabla \cdot \tilde{\bf m}]_{|T}$ in $P_k(T)$ at the $(k+2)(k+1)/2$ lattice points of $T$ (cf. \cite{Zienkiewicz}). On the other hand, for $k=0$, $d_h(\tilde{\bf m})$ is defined to be the mean value in $T$ of $[\nabla \cdot \tilde{\bf m}]_{|T}$ $\forall T \in {\mathcal T}_h$. For every $T \in {\mathcal S}_h$ such that $\Delta_T$ is not contained in $T$, we extend $d_h(\tilde{\bf m})$ to $\Delta_T$ in the natural manner.\\
\indent Now, using the affine Piola transformation for $T \in {\mathcal T}_h$ (see e.g. \cite{Ciarlet3D}) from ${\bf H}(div,T)$ onto ${\bf H}(div,\hat{T})$ where $\hat{T}$ is the standard unit master triangle, we may extend the estimate given in Proposition 3.6 of Chapter III of \cite{BrezziFortin} to the case of Sobolev spaces of order greater than one. In a more standard way, as in \cite{Ciarlet} for the second estimate, this leads to the existence of two constants $\bar{C}_{k,\tilde{\Omega}}$ and $\bar{C}_{k,d}$ such that, $\forall T \in {\mathcal T}_h \setminus {\mathcal S}_{1,h}$, it holds $\forall \tilde{\bf m} \in {\mathcal M}^{k}(\tilde{\Omega})$
\begin{equation}
\label{interpolerro}
\left\{
\begin{array}{l}
\parallel D^j[\tilde{\bf m} - \tilde{\Pi}_h \tilde{\bf m}] \parallel_{0,T} \leq \bar{C}_{k,\tilde{\Omega}} h^{k+1-j} | \tilde{\bf m} |_{k+1,T} \mbox{ for } j=0,1,\ldots,k+1\\
\parallel \nabla \cdot \tilde{\bf m} - d_h(\tilde{\bf m}) \parallel_{0,T}  \leq \bar{C}_{k,d} h^{k+1} | \nabla \cdot \tilde{\bf m} |_{k+1,T}.
\end{array}
\right.
\end{equation}
The estimates \eqref{interpolerro} are completed by means of the following technical lemma:
\begin{lemma}
\label{interperror}
Provided $h$ is sufficiently small, there exist two mesh-independent constants $\tilde{C}_{k,\tilde{\Omega}}$ and $\tilde{C}_{k,d}$ such that $\forall \tilde{\bf m} \in {\mathcal M}^{k}(\tilde{\Omega})$ fulfilling $\tilde{\bf m} \cdot {\bf n}_{|\Gamma_1} = 0$, it holds $\forall T \in {\mathcal S}_{1,h}$ 
\begin{equation}
\label{interpolerror}
\left\{
\begin{array}{l}
\parallel D^j[\tilde{\bf m} - \tilde{\Pi}_h \tilde{\bf m}] \parallel_{0,T} \leq \tilde{C}_{k,\tilde{\Omega}} h^{k+1-j} | \tilde{\bf m} |_{k+1+\varepsilon(1-\min[k,1]),T^{\Gamma}} \mbox{ for } j=0,1,\ldots,k+1\\
\parallel \nabla \cdot \tilde{\bf m} - d_h(\tilde{\bf m}) \parallel_{0,T} \leq \tilde{C}_{k,d} h^{k+1} |\nabla \cdot \tilde{\bf m}|_{k+1,T^{\Gamma}}.
\end{array}
\right.
\end{equation}
\end{lemma}
\prov
First of all we point out that the second estimate in \eqref{interpolerror} is standard, for it results from well known arguments applying to affine transformations between domains with aspect ratio uniformly bounded away from zero, such as $T^{\Gamma}$ (cf. \cite{BrennerScott} and the proof Lemma 4.1 of \cite{ZAMM}). \\
As for the first estimate, we use the inverse of the affine Piola transformation ${\mathcal P}_{T^{\Gamma}}$ from $[L^2(T^{\Gamma})]^2$ onto $[L^2(\hat{T}^{\Gamma})]^2$, where $\hat{T}^{\Gamma}$ is a triangle having two perpendicular straight edges with unit length and at most one curved edge. Then the remainder of the proof follows the same principles as in Lemma 4.1 of \cite{ZAMM}, by extending the approximation results given in Proposition 3.6 of Chapter III of \cite{BrezziFortin} to hilbertian Sobolev spaces of order greater than one and using standard interpolation theory (cf. \cite{Ciarlet,Arcangeli}). \QED \\  

An important consequence of Lemma \ref{interperror} is
\begin{theorem}
\label{Hdiverror}
Provided $h$ is sufficiently small there exists a mesh-independent constant $C_H$ such that $\forall \tilde{\bf m} \in {\mathcal M}^k(\tilde{\Omega})$ it holds
\begin{equation}
\label{HdivInterperror}
[| \tilde{\bf m} - \tilde{\Pi}_h \tilde{\bf m} |]_h \leq C_H h^{k+1} [\| \tilde{\bf m} \|_{k+1+\varepsilon(1-\min[k,1]),\tilde{\Omega}}^2 + | \nabla \cdot \tilde{\bf m} |_{k+1,\tilde{\Omega}}^2 ]^{1/2}.
\end{equation}
\end{theorem}
\prov
First of all, taking $j=0$ in the first equation of both \eqref{interpolerro} and \eqref{interpolerror} and setting $C_{k,\tilde{\Omega}}=max[\bar{C}_{k,\tilde{\Omega}},\tilde{C}_{k,\tilde{\Omega}}]$, we immediately infer that
\begin{equation}
\label{Hherror}
\parallel \tilde{\bf m} - \tilde{\Pi}_h \tilde{\bf m} \parallel_{0,h} \leq C_{k,\tilde{\Omega}} h^{k+1} | \tilde{\bf m} |_{k+1+\varepsilon(1-\min[k,1],\tilde{\Omega}}.
\end{equation}
Next we estimate the term $\parallel \nabla \cdot [\tilde{\bf m} - \tilde{\Pi}_h \tilde{\bf m}] \parallel_{0,h}$. \\
From standard results for the space $RT_k(T)$ (cf. \cite{BrezziFortin}, Chapter III, Proposition 3.6), we know that there is a suitable constant $\bar{C}_{k,H}$ independent of $h$ such that
\begin{equation}
\label{errdiv0}
\parallel \nabla \cdot [\tilde{\bf m} - \tilde{\Pi}_h \tilde{\bf m}] \parallel_{0,T} \leq \bar{C}_{k,H} h^{k+1} \| \nabla \cdot \tilde{\bf m} \|_{k+1,T} \; 
\forall T \in {\mathcal T}_h \setminus {\mathcal S}_{1,h}.
\end{equation}
On the other hand, $\forall T \in {\mathcal S}_{1,h}$ we first write 
\begin{equation}
\label{errdiv1} 
\left\{
\begin{array}{l}
\| \nabla \cdot [\tilde{\bf m} - \tilde{\Pi}_h \tilde{\bf m} \|^2_{0,T} = ( \nabla \cdot [\tilde{\bf m} - \tilde{\Pi}_h \tilde{\bf m}],\nabla \cdot \tilde{\bf m} - d_h(\tilde{\bf m}))_T \\
+  ( \nabla \cdot [\tilde{\bf m} - \tilde{\Pi}_h \tilde{\bf m}],d_h(\tilde{\bf m})-\tilde{\Pi}_h \tilde{\bf m})_T.
\end{array}\textit{}
\right.
\end{equation}
Application of Young's inequality to \eqref{errdiv1}, followed by the use of the second equation of \eqref{interpolerror}, easily yields
\begin{equation}
\label{errdiv1bis}
\left\{
\begin{array}{l} 
\| \nabla \cdot [\tilde{\bf m} - \tilde{\Pi}_h \tilde{\bf m} \|^2_{0,T} \leq \tilde{C}_{k,d}^2 h^{2(k+1)} |\nabla \cdot \tilde{\bf m}|^2_{k+1,T^{\Gamma}}+ 
2 \int_T \nabla \cdot [\tilde{\bf m} - \tilde{\Pi}_h \tilde{\bf m}]v \\
\mbox{with } v:=d_h(\tilde{\bf m}) - \nabla \cdot \tilde{\Pi}_h \tilde{\bf m}.
\end{array}
\right. 
\end{equation}
Now, from the construction of $\tilde{\Pi}_h \tilde{\bf m}$, we readily establish that $\forall T \in {\mathcal S}_h$ and $\forall v \in P_k(T)$ it holds
\begin{equation}
\label{errdiv2}
\int_T \nabla \cdot [\tilde{\bf m} - \tilde{\Pi}_h \tilde{\bf m}] v = \int_{e_T} [\tilde{\bf m} - \tilde{\Pi}_h \tilde{\bf m}] \cdot {\bf n}_T v \; ds \; \forall v \in P_k(T).
\end{equation}
Let us set $\sigma(T) = 1$ if $T \in {\mathcal S}_{1,h}^{'}$ and $\sigma(T) = -1$ if $T \in {\mathcal S}_{1,h} \setminus {\mathcal S}_{1,h}^{'}$. Using the Divergence Theorem and taking into account that $\tilde {\bf m} \cdot {\bf n} \equiv 0$ on $\Gamma_T$, after straightforward calculations starting from \eqref{errdiv2}, we obtain $\forall v \in P_k(T^{\Gamma})$,
\begin{equation}
\label{errdiv3}
\left\{
\begin{array}{l}
\int_T \nabla \cdot [\tilde{\bf m} - \tilde{\Pi}_h \tilde{\bf m}] v = I_{1,T}(\tilde{\bf m},v) + I_{2,T}(\tilde{\bf m},v) + I_{3,T}(\tilde{\bf m},v) \\
\mbox{where} \\
I_{1,T}(\tilde{\bf m},v): = - \sigma(T) \int_{\Gamma_T} [\tilde{\Pi}_h \tilde{\bf m}] \cdot {\bf n} \; v \; ds\\
I_{2,T}(\tilde{\bf m},v): = \sigma(T) \int_{\Delta_T} [\tilde{\bf m} - \tilde{\Pi}_h \tilde{\bf m}] \cdot \nabla v \\
I_{3,T}(\tilde{\bf m},v): = \sigma(T) \int_{\Delta_T} \nabla \cdot [\tilde{\bf m} - \tilde{\Pi}_h \tilde{\bf m}] v.
\end{array}
\right.
\end{equation}
\indent Let us estimate $|I_{i,T}(\tilde{\bf m},v)|$, $i=1,2,3$.\\

First of all, recalling Lemma \ref{GammaT} and using the first equation of \eqref{interpolerror} with $j=k+1$, we may write for a mesh-independent constant $C_{I1}$
\begin{equation}
\label{errI1}
|I_{1,T}(\tilde{\bf m},v)| \leq C_{I1} h_T^{k+1} \| \tilde{\bf m} \|_{k+1+\varepsilon(1-\min[k,1]),T^{\Gamma}} \| v \|_{0,T}.
\end{equation}

Next we turn our attention to $I_{2,T}(\tilde{\bf m},v)$. Since this quantity vanishes for $k=0$ $\forall v \in P_k(T)$, we handle it only for $k \geq 1$ as follows: \\ 
First we denote by $\hat{\Pi}_T$ the counterpart in $\hat{T}^{\Gamma}$ of the interpolation operator $\tilde{\Pi}_h$ restricted to $T^{\Gamma}$ (cf. \cite{BrezziFortin,Ciarlet}). In doing so, recalling the Piola transformation ${\mathcal P}_{T^{\Gamma}}$ from $[L^2(T^{\Gamma})]^2$ onto $[L^2(\hat{T}^{\Gamma})]^2$, we set $\hat{\bf m}:= 
{\mathcal P}_{T^{\Gamma}}\tilde{\bf m}$. Then denoting by $\rho_T$ the radius of the circle inscribed in $T$, referring to \cite{Ciarlet,Ciarlet3D}, we observe that    \\ 
\begin{equation}
\label{errdiv4}
\| \tilde{\bf m} - \tilde{\Pi}_h \tilde{\bf m} \|_{0,\infty,\Delta_T} \leq [(2-\sqrt{2}) \rho_T]^{-1} \| \hat{\bf m} - \hat{\Pi}_T \hat{\bf m}\|_{0,\infty,\hat{T}^{\Gamma}} 
\end{equation}
Now let $\widehat{\tilde{T}}$ be the smallest isosceles right triangle that contains all $\hat{T}^{\Gamma}$ corresponding to $T \in {\mathcal S}_{1,h}$ for all ${\mathcal T}_h \in {\mathcal F}$ having the same right-angle vertex as all of them . It is not difficult to figure out that the diameter $\hat{h}$ of $\widehat{\tilde{T}}$ is bounded above by $1+\hat{c} h$, where $\hat{c}$ is a constant independent of $h$, and $area(\widehat{\tilde{T}})$ is bounded below by $1/2$. Next, applying the inverse of ${\mathcal P}_{T^{\Gamma}}$ to $[L^2(\widehat{\tilde{T}})]^2$ as well, we consider the extension $\widehat{\tilde{\bf m}}$ of $\hat{\bf m}$ to $\widehat{\tilde{T}} \setminus \hat{T}^{\Gamma}$ belonging to $[H^{k+1}(\widehat{\tilde{T}})]^2$ given by $\widehat{\tilde{\bf m}} = {\mathcal P}_{T^{\Gamma}} \tilde{\bf m}$. \\
Let us extend $\hat{\Pi}_T \hat{\bf m}$ to the whole $\widehat{\tilde{T}}$. This extension in $\widehat{\tilde{T}} \setminus \hat{T}^{\Gamma}$ equals $\hat{\Pi}_T \widehat{\tilde{\bf m}}$ by construction. Indeed, the expression of $\hat{\Pi}_T \widehat{\tilde{\bf m}}$ depends only on DOFs attached to $T^{\Gamma}$. Thus it is the same polynomial field as $\hat{\Pi}_T \hat{\bf m}$ in $\hat{T}^{\Gamma}$. It trivially follows that for a constant $C^{\Gamma}$ independent of $T^{\Gamma}$ it holds
\begin{equation}
\label{errdiv5}
\| \hat{\Pi}_T \hat{\bf m} \|_{0,\infty,\hat{T}^{\Gamma}} \leq C^{\Gamma} \| \widehat{\tilde{\bf m}} \|_{0,\infty,\hat{T}^{\Gamma}} \leq 
C^{\Gamma} \| \widehat{\tilde{\bf m}} \|_{0,\infty,\widehat{\tilde{T}}}.
\end{equation}
Incidentally we observe that, by the Rellich-Kondratchev Theorem $H^{k+1}(\widehat{\tilde{T}})$ is continuously embedded in $L^{\infty}(\widehat{\tilde{T}})$ for $k \geq 1$. Moreover, since $\widehat{\tilde{\bf m}} \equiv \hat{\Pi}_T \widehat{\tilde{\bf m}}$ whenever $\widehat{\tilde{\bf m}} \in P_{k}(\widehat{\tilde{T}})$, by the Bramble-Hilbert Lemma we can ensure the existence of a constant $\widehat{\tilde{C}}$ independent of $h$, such that
\begin{equation}
\label{errdiv6}
\| \hat{\bf m} - \hat{\Pi}_T \hat{\bf m} \|_{0,\infty,\hat{T}^{\Gamma}}  \leq 
\| \widehat{\tilde{\bf m}} - \hat{\Pi}_T \widehat{\tilde{\bf m}} \|_{0,\infty,\widehat{\tilde{T}}} \leq \widehat{\tilde{C}} | \widehat{\tilde{\bf m}} |_{k+1,\widehat{\tilde{T}}}.
\end{equation}
Now proceeding to the estimate of $I_{2,T}(\tilde{\bf m},v)$, we first have 
\begin{equation}
\label{errI21}
|I_{2,T}(\tilde{\bf m},v)| \leq C_0 h_T^3 \| \tilde{\bf m} - \tilde{\Pi}_h \tilde{\bf m} \|_{0,\infty,\Delta_T} \| \nabla v \|_{0,\infty,T^{\Gamma}}.
\end{equation}
Then, recalling \eqref{L2TDelta}, we obtain 
\begin{equation}
\label{errI22}
|I_{2,T}(\tilde{\bf m},v)| \leq {\mathcal C}_0^{+} C_0 h_T^{2}  \| \tilde{\bf m} - \tilde{\Pi}_h \tilde{\bf m} \|_{0,\infty,\Delta_T} \| \nabla v \|_{0,T}.
\end{equation}
Further resorting to a standard inverse inequality, for a suitable constant $\bar{C}_0$ independent of $T$ it holds
\begin{equation}
\label{errI23} 
|I_{2,T}(\tilde{\bf m},v)| \leq \bar{C}_{0} h_T \| \tilde{\bf m} - \tilde{\Pi}_h \tilde{\bf m} \|_{0,\infty,\Delta_T} \| v \|_{0,T}.
\end{equation}
Combining \eqref{errdiv4} and \eqref{errdiv6} and plugging the resulting inequality into \eqref{errI23} we come up with
\begin{equation}
\label{errI24} 
|I_{2,T}(\tilde{\bf m},v)| \leq \bar{C}_{0} \widehat{\tilde{C}} h_T [(2-\sqrt{2})\rho_T]^{-1} | \widehat{\tilde{\bf m}} |_{k+1,\widehat{\tilde{T}}} \| v \|_{0,T}
\end{equation}
Finally, let $\tilde{T} \subset \tilde{\Omega}$ be the reciprocal image of $\widehat{\tilde{T}}$ under the same affine transformation from $T$ onto the standard master triangle $\hat{T}$. Using standard results for affine Piola transformations such as ${\mathcal P}_{T^{\Gamma}}$ (see e.g. \cite{Ciarlet,Ciarlet3D}), for a suitable mesh-independent constant $C_{I2}$, \eqref{errI24} easily leads to 
\begin{equation}
\label{errI2} 
|I_{2,T}(\tilde{\bf m},v)| \leq  C_{I2} h_T^{k+1} | \tilde{\bf m} |_{k+1,\tilde{T}} \| v \|_{0,T}. 
\end{equation} 
Notice that several pairs of triangles $\tilde{T}$ overlap, but this does not affect the final result given below.\\

As for $I_{3,T}(\tilde{\bf m},v)$ we trivially have
\begin{equation}
\label{errI31}
|I_{3,T}(\tilde{\bf m},v)| \leq [C_{0} h_T^3]^{1/2} 
\| \nabla \cdot [\tilde{\bf m} - \tilde{\Pi}_h \tilde{\bf m} \|_{0,T^{\Gamma}} \| v \|_{0,\infty,T^{\Gamma}}.
\end{equation}
Recalling \eqref{L2TDelta}, for a suitable constant $\tilde{C}_{I3}$ the following estimate easily derives from \eqref{errI31}:
\begin{equation}
\label{errI33}
|I_{3,T}(\tilde{\bf m},v)| \leq \tilde{C}_{I3}  h_T^{1/2} \| \nabla \cdot [\tilde{\bf m} - \tilde{\Pi}_h \tilde{\bf m}] \|_{0,T^{\Gamma}} \| v \|_{0,T}.
\end{equation}
Now applying the triangle inequality and the second estimate of \eqref{interpolerror} to \eqref{errI33} we easily come up with a constant $C_{I3}$ such that 
\begin{equation}
\label{errI3}
|I_{3,T}(\tilde{\bf m},v)| \leq C_{I3} h_T^{1/2} ( h_T^{k+1} | \nabla \cdot \tilde{\bf m} |_{k+1,T^{\Gamma}} \| v\|_{0,T} + \| v \|^2_{0,T} ).
\end{equation}
Plugging \eqref{errI1}, \eqref{errI2}, \eqref{errI3} into \eqref{errdiv3} and setting $\bar{C}_I:=\max[C_{I1},C_{I2},C_{I3}]$, we obtain $\forall T \in {\mathcal S}_{1,h}$
\begin{equation}
\label{errdiv7}
\begin{array}{l}
\| \nabla \cdot (\tilde{\bf m} - \tilde{\Pi}_h \tilde{\bf m}) \|_{0,T}^2 \leq \bar{C}_I h_T^{k+1} \{ \| \tilde{\bf m} \|_{k+1+\varepsilon(1-\min[k,1]),T^{\Gamma}} + | \tilde{\bf m} |_{k+1,\tilde{T}} + h_T^{1/2} | \nabla \cdot \tilde{\bf m} |_{k+1,T^{\Gamma}} \} \\
\| v \|_{0,T} + C_{I3} h_T^{1/2} \| v \|_{0,T}^2.
\end{array}
\end{equation}
Then, since $h_T<1$, after straightforward manipulations we further establish the existence of a constant $\tilde{C}_I$ independent of $T$ such that
\begin{equation}
\label{errdiv8}
(1-4C_{I3}h_T^{1/2}) \| \nabla \cdot (\tilde{\bf m} - \tilde{\Pi}_h \tilde{\bf m}) \|_{0,T}^2 \leq \tilde{C}_I^2 h_T^{2k+2} 
\{ \| \tilde{\bf m} \|_{k+1+\varepsilon(1-\min[k,1]),T^{\Gamma}}^2 + | \tilde{\bf m} |_{k+1,\tilde{T}}^2 + | \nabla \cdot \tilde{\bf m} |_{k+1,T^{\Gamma}}^2 \}.
\end{equation}
Assuming that $h_T \leq 1/(8C_{I3})^2$, we combine \eqref{errdiv8} and \eqref{errdiv0} to conclude that the interpolation error of $\nabla \cdot \tilde{\bf m}$ can be estimated by 
\begin{equation}
\label{Hderror}
\parallel \nabla \cdot [\tilde{\bf m} - \tilde{\Pi}_h \tilde{\bf m}] \parallel_{0,h} \leq C_{k,div} h^{k+1} \{ \| \tilde{\bf m} \|_{k+1+\varepsilon(1-\min[k,1]),\tilde{\Omega}}
+ \| \tilde{\bf m} \|_{k+1,\tilde{\Omega}} + | \nabla \cdot \tilde{\bf m} |_{k+1,\tilde{\Omega}} \},
\end{equation}
for a suitable constant $C_{k,div}$ independent of $h$. \\
Owing to \eqref{Hderror} and \eqref{Hherror}, the result immediately follows with $C_H= [ C_{k,div}^2+C_{k,\tilde{\Omega}}^2 ]^{1/2}$\QED \\

\section{Error estimates}

\hspace{4mm}In this section we endeavor to obtain error estimates for problem \eqref{Dmixh} based on the error bound \eqref{unifymajor}. For the sake of clarity we distinguish the cases of convex and non convex domains.

\subsection{The convex case} 

\hspace{4mm} Here \eqref{unifymajor} can be applied with $\gamma_h=\gamma$ and $F_h(v)=(f,v)_h$. 
We consider separately the two possibilities, namely,\\
 
\noindent \underline{\textbf{I) $\Omega$ is convex and $length(\Gamma_0) > 0$}} \\
\indent In this case we may take $(\tilde{\bf p};\tilde{u})=({\bf p};u)$.\\
The $inf$ term in \eqref{unifymajor} can be bounded above by $||| ({\bf p};u) - (\tilde{\Pi}_h{\bf p};\tilde{\pi}_h u) |||_h$, where $\tilde{\pi}_h u$ is the standard Lagrange interpolate of $u$ for $k \geq 1$ and the mean value of $u$ in each triangle of the mesh for $k=0$. We know that, as long as $u \in H^{k+1}(\Omega)$, there is a mesh-independent constant $C_{\pi}$ such that, for $k=0,1, \ldots$ (cf. \cite{Ciarlet}) 
\begin{equation}
\label{pih}
     \| u - \tilde{\pi}_h u \|_{0,h} \leq C_{\pi} h^{k+1} | u |_{k+1}.
\end{equation}

On the other hand, the variational residual in the $sup$ term of  \eqref{unifymajor} is given by
\begin{equation}
\label{varesidual2}
{\mathcal R}_h(u,{\bf q}):=({\bf p},{\bf q})_h - (u,\nabla \cdot {\bf q})_h = -(u,{\bf q} \cdot {\bf n}_h)_{0,\Gamma_{0,h}} \; \forall {\bf q} \in {\bf Q}_h, 
\end{equation}
In principle, ${\mathcal R}_h(u,{\bf q})$ does not vanish, since $u \equiv 0$ on $\Gamma_{0}$ but not necessarily on $\Gamma_{0,h}$. Nevertheless, we can give a rough estimate of this residual, in order to determine our method's real order in terms of $h$ if $length(\Gamma_0) > 0$.
\begin{e-proposition}
\label{neq0convex}
Let $length(\Gamma_0)> 0$ and ${\mathcal R}_h(u,{\bf q})$ be given by \eqref{varesidual2}. If the first (resp. second) component ${\bf p}$ (resp. u) of the solution of \eqref{Dmix}-\eqref{definec} belongs to $[H^{1}(\Omega)]^2$ (resp. $H^{2}(\Omega)$) there exists a suitable mesh-independent constant $C_{\Gamma_0}$ such that for any $k \geq 0$ it holds 
\begin{equation}
\label{estimgamma0}
 |{\mathcal R}_h(u,{\bf q})| \leq C_{\Gamma_0} h^{3/2} \| u \|_ {2} [| {\bf q} |]_h \; \forall {\bf q} \in {\bf Q}_h. 
\end{equation} 
\end{e-proposition}
\prov
First of all we have
\begin{equation}
\label{Gamma01}
(u,{\bf q} \cdot {\bf n}_h)_{e_T} \leq \| u \|_{0,e_T} \| {\bf q} \cdot {\bf n}_h \|_{0,e_T} \; \forall T \in {\mathcal S}_{0,h}.
\end{equation}
By a straightforward manipulation already performed several times in this work, we come up with a constant $C_{0I}$ independent of $T$ such that
\begin{equation}
\label{Gamma02} 
\| {\bf q} \cdot {\bf n}_h \|_{0,e_T} \leq l_T^{1/2} \| {\bf q} \|_{0,\infty,T} \leq C_{0I} h_T^{-1/2} \| {\bf q} \|_{0,T}.
\end{equation}
Moreover, since $u$ vanishes at both ends of $e_T$, we know that, for a certain constant $C_{e0}$, it holds 
\begin{equation}
\label{Gamma03}
\|u\|_{0,e_T} \leq C_{e0} h_T^{j} \| \partial^j u/\partial x^j \|_{0,e_T} \mbox{ for } 1 \leq j \leq 2.
\end{equation}
Hereafter we take $j=1$. Recalling the local coordinate system $(A_T;x,y)$ attached to $e_T$ and the function $\phi$ such that $M=(x,0) \in e_T$ corresponds to $N=(x,\phi(x)) \in \Gamma_T$. We have
\begin{equation}
\label{Gamma04}
(\partial u/\partial x)(x,0) = (\partial u/\partial x)(x,\phi(x)) - 
\int_M^N \partial^{2} u/\partial x \partial y \; dy.
\end{equation}
${\bf t}$ being the tangent unit vector along $\Gamma$ we introduce the direct reference frame ${\mathcal R}:=(N;{\bf t},{\bf n})$. Then recalling the angle $\varphi$ between ${\bf n}_h$ and ${\bf n}$, let the unit vector along $e_T$ be ${\bf l}_T$. We have 
${\bf l}_T = \cos{\varphi} {\bf t} - \sin{\varphi} {\bf n}$ and hence we have $\partial u/\partial x = \nabla u_{|{\mathcal R}} \cdot (\cos{\varphi};- \sin{ \varphi})$. We know that $|\sin{\varphi} | \leq C_{\Gamma} h_T$ owing to \eqref{CGamma} and also that all the tangential derivatives of $u$ of any order vanish on $\Gamma_0$. This easily yields 
\begin{equation}
\label{Gamma05}  
|(\partial u/\partial x)(x,\phi(x))| \leq C_{\Gamma} h_T |(\partial u/\partial n)(x,\phi(x))| \; \forall x \in [0,l_T].
\end{equation}
Plugging \eqref{Gamma05} into \eqref{Gamma04} we obtain 
\begin{equation}
\label{Gamma06}
|(\partial u/\partial x)(x,0)| \leq C_{\Gamma} h_T | (\partial u/\partial n)(x,\phi(x))| + C_0^{1/2} h_T \displaystyle 
\left[ \int_M^N |\partial^{2} u/\partial x \partial y|^2 dy \right]^{1/2} \; \forall x \in [0,l_T].
\end{equation}
It is easy to see that $\cos{\varphi} \; ds= dx$ where $s$ is the curvilinear abscissa along $\Gamma_T$, that is, $ds$ is just a small perturbation of $dx$. Thus plugging in turn \eqref{Gamma06} into \eqref{Gamma03}, it trivially holds for a certain constant $C_{e1}$ independent of $T$:
\begin{equation}
\label{Gamma07}
\|u\|^2_{0,e_T}  \leq C^2_{e1} h_T^{4} [\| \partial u/\partial n \|^2_{0,\Gamma_T} + | u |^2_{2,\Delta_T}]
\end{equation}
Finally recalling \eqref{Gamma01},\eqref{Gamma02} and sweeping the triangles $T \in {\mathcal S}_{0,h}$, straightforward manipulations yield 
\begin{equation}
\label{Gamma8}
(u,{\bf q}\cdot {\bf n}_h)_{0,\Gamma_{0,h}} \leq C_{0I}C_{e1}  h^{3/2} \{ \| \partial u/ \partial n \|_{0,\Gamma_0} +  
| u |_{2} \}[| {\bf q} |]_h.
\end{equation}
Finally, thanks to the Trace Theorem for $H^2(\Omega)$ (cf. \cite{Adams}), \eqref{estimgamma0} derives from \eqref{Gamma8}.  \QED \\

As a consequence of Proposition \ref{neq0convex} we trivially have
\begin{theorem}
\label{Estimix}
Assume that $length(\Gamma_0) > 0$, $\Omega$ is smooth enough, $h$ is sufficiently small and $k \geq 0$. If for a small $\varepsilon >0$ $\tilde{\bf p} \in [H^{k+1+\varepsilon(1-\min[1,k])}(\tilde{\Omega})]^2$, $\nabla \cdot {\bf p} \in H^{k+1}(\Omega)$ and $u \in H^{\max[2,k+1]}(\Omega)$, the following error estimate holds with a mesh-independent constant $C_{01}$, for the approximation of the solution of \eqref{Dmix}-\eqref{definec} by the solution of \eqref{Dmixh}.
\begin{equation}
\label{Esticonvexmix}
||| ({\bf p};u) - ({\bf p}_h;u_h) |||_h \leq  \displaystyle C_{01} \left\{ h^{3/2} \| u \|_{2} + h^{k+1}
\left[ \| \tilde{\bf p} \|_{k+1+\varepsilon(1-\min[1,k]),\tilde{\Omega}} + | \nabla \cdot {\bf p} |_{k+1} + | u |_{k+1} \right] \right\}. 
\mbox{\QED}
\end{equation}
\end{theorem}

\begin{remark}
According to the estimate \eqref{Esticonvexmix}, in case $length(\Gamma_0) > 0$, the order of our method is optimal (up to a very small number $\varepsilon$) only for $k=0$. But there should be no grievance about this fact, since in general the solution $u$ of the mixed Poisson problem does not even belong to $H^2(\Omega)$, but rather to $H^{3/2+\varepsilon}(\Omega)$ \cite{Grisvard}. Therefore ${\bf p} \in [H^{1/2+\varepsilon}(\Omega)]^2$ in general, and hence the best order we can hope for here is a little more than $1/2$. Nevertheless there is a rather simple way out, if the exact solution happens to be more regular, as we will see in Subsection 7.3. Incidentally, we will also see later on that, by requiring more regularity from $u$, our method's optimal (second) order can be recovered for $k=1$, in the case of mixed boundary conditions.   \QED
\end{remark}

\noindent \underline{\textbf{II) $\Omega$ is convex and $length(\Gamma_0) = 0$}} \\
\indent Here it is advisable to take $(\tilde{\bf p};\tilde{u})=({\bf p};\bar{u})$ where $\bar{u}$ differs from $u$ by a suitable additive constant. Indeed, since $u \in L^2_0(\Omega)$ in this case, a natural way to define the interpolate of $u$ in $V_h$ is as follows: We start with the operator $\tilde{\pi}_h$ from $V^{k}(\Omega):=H^{k+1}(\Omega) \cap V$ into $L^2(\Omega_h)$ already defined in I). 
Noticing that $\tilde{\pi}_h u$ does not necessarily belong to $L^2_0(\Omega_h)$, $\tilde{\pi}_h u$ may not belong to $V_h$. Thus we set $\bar{\pi}_h u := \tilde{\pi}_h u + C_h \in V_h$ with $C_h:=- \int_{\Omega_h} \tilde{\pi}_h u/area(\Omega_h)$. While on the one hand the estimate of the interpolation error in $V_h$ of $u$ is not optimal, we note that the estimate $\| u - \tilde{\pi}_h u \|_{0,h} \leq C_{\pi}  h^{k+1} | u |_{k+1}$ holds if $u \in V^{k}(\Omega)$. Now we define $\bar{u} = u + C_h$, which leads to the desired optimal estimate of the best approximation in ${\bf P}_h \times V_h$ of the pair $({\bf p};\bar{u})$ instead of $({\bf p};u)$, as seen below. \\
It is noticeable that, since $\int_{\Omega} \nabla \cdot {\bf q} =0$ $\forall {\bf q} \in {\bf Q}$, the replacement of $u$ by $\bar{u}$ does not change the essence of problem \eqref{Dmix}-\eqref{definec}; it is just a way like any other to fix the additive constant up to which $u$ is defined. \\
Summarizing, the error bound \eqref{unifymajor} applies to the pair $(\tilde{\bf p};\tilde{u})=({\bf p},\bar{u})$. 
It turns out that  
\begin{equation} 
\label{barpih}
\displaystyle \inf_{w \in V_h} \| \bar{u} - w \|_{0,h} \leq \| \bar{u} - \bar{\pi}_h u \|_{0,h} =  
\| u - \tilde{\pi}_h u \|_{0,h} \leq C_{\pi} h^{k+1} | u |_{k+1}.
\end{equation}
Thus recalling \eqref{HdivInterperror}, provided $\tilde{\bf p} \in [H^{k+1+\varepsilon(1-\min[k,1])}(\tilde{\Omega})]^2$ and $\nabla \cdot {\bf p} \in  H^{k+1}(\Omega)$, setting $C_{I}:=C_H + C_{\pi}$, we come up with, 
\begin{equation}
\label{pairerror}
\inf_{({\bf r};w) \in {\bf P}_h \times V_h}  ||| ({\bf p};\bar{u}) - ({\bf r};w) |||_h \leq   C_I h^{k+1} \displaystyle 
\left\{ \| \tilde{\bf p} \|_{k+1+\varepsilon(1-\min[k,1]),\tilde{\Omega}} + | \nabla \cdot {\bf p} |_{k+1} + | u |_{k+1} \right \}.
\end{equation}
As for the $sup$ term in \eqref{unifymajor}, we can assert that it vanishes identically if $length(\Gamma_0)=0$. Indeed, $\nabla \cdot  {\bf p} = - f$ in $\Omega_h$ and on the other hand, as long as $u \in H^{k+1}(\Omega)$, 
\begin{equation}
\label{varesidual1}
({\bf p},{\bf q})_h - (\bar{u},\nabla \cdot {\bf q})_h = ({\bf p}-\nabla \bar{u},{\bf q})_h = 0  \; \forall {\bf q} \in {\bf Q}_h, 
\end{equation}
since ${\bf q} \cdot {\bf n}_h$ vanishes everywhere in $\Gamma_h$ in this case.\\
This leads to the following result.
\begin{theorem}
\label{EstiNeumann}
Assume that $length(\Gamma_0) = 0$, $\Omega$ is smooth enough and $h$ is sufficiently small. If $\tilde{\bf p} \in [H^{k+1+\varepsilon(1-\min[k,1])}(\tilde{\Omega})]^2$, $\nabla \cdot {\bf p} \in  H^{k+1}(\Omega)$ and $u \in H^{k+1}(\Omega)$, the following error estimate holds for the approximation of the solution of \eqref{Dmix}-\eqref{definec} by the solution of \eqref{Dmixh}.
\begin{equation}
\label{EsticonvexNeumann}
\left\{
\begin{array}{l}
||| ({\bf p};\bar{u}) - ({\bf p}_h;u_h) |||_h \leq  \displaystyle \frac{C_I}{\gamma} h^{k+1}
\left\{ \| \tilde{\bf p} \|_{k+1+\varepsilon(1-\min[k,1]),\tilde{\Omega}} + | \nabla \cdot {\bf p} |_{k+1} + | u |_{k+1} \right\} \\ 
\mbox{where } u - \bar{u} = \mbox{ constant, for } k=0,1,2,\ldots . \mbox{\QED}
\end{array}
\right.
\end{equation}
\end{theorem}
\subsection{The non convex case}

\hspace{4mm} Here we can apply \eqref{unifymajor} with $\gamma_h=\gamma^{'}$, (cf. \eqref{gammaprime}) but $F_h(v)$ cannot be $(f,v)_h$. Besides this one, new issues must be overcome in the non convex case as far as error estimates are concerned. \\
First of all, whenever $\Omega$ is not convex, the error cannot be measured in norms associated with the set $\Omega_h$, since the exact solution is not defined in $\Omega_h \setminus \bar{\Omega} \neq \emptyset$. Akin to \cite{ZAMM} the error estimate is now expressed for $||| ({\bf p};u) - ({\bf p}_h;u_h) |||_h^{'}$ which is not greater than  
$||| (\tilde{\bf p};\tilde{u}) - ({\bf p}_h;u_h) |||_h$ for a suitable pair $(\tilde{\bf p};\tilde{u})$ defined hereafter. This is in turn bounded above as in \eqref{unifymajor}. \\
As for functional 
$F_h$, it must be defined in such a way that 
no unknown values of an extension $\tilde{f}$ of $f$ to $\Omega_h \setminus \Omega$ come into play. Notice that the replacement of $(\tilde{f},v)_h$ by $F_h(v)$ on the right hand side of \eqref{Dmixh} is an old trick considered in \cite{CiarletRaviart} to handle curved non convex domains.\\  
Incidentally we refer to Subsection 5.3 for some definitions and concepts evoked hereafter concerning extensions to $\tilde{\Omega}$, i.e.,  
a smooth open set sufficiently large to strictly contain not only $\Omega$, but also the sets $\tilde{\Omega}_h$ for all ${\mathcal T}_h \in {\mathcal F}^{'}$.\\
\indent Actually the pair $(\tilde{\bf p};\tilde{u})$ extending $({\bf p};u)$ to $\tilde{\Omega}$ must have the required regularity for method's order to remain the same as in the convex case. Akin to \cite{ZAMM}, the construction of the extension $\tilde{u}$ of $u$ - which incidentally induces regular extensions $\tilde{\bf p}:= \nabla \tilde{u}$ and $\tilde{f}:=-\Delta \tilde{u}$ of ${\bf p}$ and $f$ -, is the one advocated in \cite{Stein}. Of course, the explicit knowledge of neither $(\tilde{\bf p};\tilde{u})$ nor $\tilde{f}$ is necessary. The point here is that all these quantities do exist, and can thus be employed in error estimates for the non convex case, which we next complete. \\
Notice that the extension technique allows for establishing 
an estimate for the interpolation error in the non convex case analogous to \eqref{HdivInterperror}.
More precisely if $\tilde{\bf m} \in {\mathcal M}^{k}(\tilde{\Omega})$ we have
\begin{equation}
\label{HdivInterperro}
[| \tilde{\bf m} - \tilde{\Pi}_h \tilde{\bf m} |]_h \leq C_H h^{k+1} [\| \tilde{\bf m} \|_{k+1+\varepsilon(1-\min[k,1]),\tilde{\Omega}}^2 + | \nabla \cdot \tilde{\bf m} |_{k+1,\tilde{\Omega}}^2 ]^{1/2}.
\end{equation}
\indent Furthermore, we may exploit the extension $\tilde{f}$ in the framework of the approximate right hand side $F_h(v)$ of problem \eqref{Dmixh}. 
Among many other possibilities, for the sake of simplicity we consider here a bilinear form $F_h$ defined as follows. First of all, assuming that $\tilde{f} \in H^{k+1+\varepsilon(1-\min[k,1])}(\tilde{\Omega}_h)$, let $f_h$ be the Lagrange $P_k(T)$-interpolate of $\tilde{f}$ at the set of $(k+2)(k+1)/2$ points in $T$, $\forall T \in {\mathcal T}_h$, which form the standard set of lattice points of $T_{\chi}$, where $T_{\chi}$ is the reduction of $T$ by homothety with ratio $\chi <1$ centered at the centroid of $T$. $\chi$ is chosen in such a way that $T_{\chi} \subset \Omega$ $\forall T \in {\mathcal S}_h$ for all ${\mathcal T}_h$ in ${\mathcal F}^{'}$. \\
Then we set 
\begin{equation}
\label{Fh}
F_h(v):= (f_h,v)_h \; \forall v \in V_h.
\end{equation}
Since by standard results \cite{Ciarlet} there is a constant $C_{\chi,k}$ independent of $h$ such that
$$\| \tilde{f} - f_h \|_{0,h} \leq C_{\chi,k} h^{k+1} \| \tilde{f} \|_{k+1+\varepsilon(1-\min[k,1]),\Omega_h},$$
we obtain from the Cauchy-Schwarz inequality, together with \eqref{Fh}
\begin{equation}
\label{errf}
|(f,v)_h - F_h(v)| \leq C_{\chi,k} h^{k+1} \| \tilde{f} \|_{k+1+\varepsilon(1-\min[k,1]),\Omega_h} \| v \|_{0,h} \forall v \in V_h
\end{equation}
Henceforth we definitively take $\tilde{f}=-\Delta \tilde{u}$ in $\tilde{\Omega} \setminus \Omega$ as the extension of $f$ with the required smoothness.\\
In view of all the above considerations, we are now ready to establish error estimates for problem \eqref{Dmixh} in the case  of non convex domains. Akin to the previous subsection we study separately the two possible cases, namely,\\

\noindent \underline{\textbf{III) $\Omega$ is not convex and $length(\Gamma_0) > 0$}} \\
\indent We next focus on the estimation of the contribution to the approximation error of the variational residual incorporating the functional $F_h(v)$, namely, $|((\tilde{\bf p};\tilde{u}),({\bf q};v))-F_h(v)|$. With this aim we first write
\begin{equation}
\label{residual1}
| c_h((\tilde{\bf p};\tilde{u}),({\bf q};v))-F_h(v)| \leq | c_h((\tilde{\bf p};\tilde{u}),({\bf q};v))-(\tilde{f},v)_h| + |(\tilde{f},v)_h-F_h(v)|.
\end{equation} 
Now we observe that, either by construction, or due to \eqref{Mixed}, it holds $-\nabla \cdot \tilde{\bf p} = \tilde{f}$ and $\tilde{\bf p}=\nabla \tilde{u}$ almost everywhere in $\tilde{\Omega}$. Thus assuming that $\tilde{u} \in H^2(\tilde{\Omega})$, the first term on the right hand side of \eqref{residual1} is ${\mathcal R}_h(\tilde{u},{\bf q}):=(\tilde{u},{\bf q} \cdot {\bf n}_h)_{0,\Gamma_{0,h}}$. It is easy to see that ${\mathcal R}_h(\tilde{u},{\bf q})$ can be estimated in the same manner as in 
\eqref{estimgamma0}, that is,
\begin{equation}
\label{tildestimgamma0} 
 |{\mathcal R}_h(\tilde{u},{\bf q})| \leq C_{\Gamma_0} h^{3/2} \| \tilde{u} \|_{2,\tilde{\Omega}_h} [| {\bf q} |]_h \; \forall {\bf q} \in {\bf Q}_h.
\end{equation}
Combining \eqref{errf} and \eqref{tildestimgamma0} and setting $\tilde{C}_{r}:=\max[C_{\Gamma_0},C_{\chi,k}]$, we immediately come up with
\begin{equation}
\label{residual} 
 \displaystyle \sup_{({\bf q};v) \in {\bf Q}_h \times V_h/|||({\bf q};v)|||_h=1}| c_h((\tilde{\bf p};\tilde{u}),({\bf q};v))- F_h(v) | 
\leq \tilde{C}_{r} \displaystyle \left[ h^{\frac{3}{2}} \! \| \tilde{u} \|_{2,\tilde{\Omega}}+ h^{k+1} \| \tilde{f} \|_{k+1+\varepsilon(1-\min[k,1]),\tilde{\Omega}} \right]. 
\end{equation}
Now taking $\tilde {\bf m} = \tilde{\bf p}$ in \eqref{HdivInterperro} and letting $\tilde{\pi}_h \tilde{u}$ be the $V_h$-interpolate of $\tilde{u}$ for $k>0$ and the mean value of $\tilde{u}$ in every triangle of the mesh for $k=0$, we can assert that there exists a constant $\tilde{C}_I$ independent of $h$ such that the following estimate holds: 
\begin{equation}
\label{bestapprox}
\begin{array}{l}
\inf_{({\bf r};w) \in {\bf P}_h \times V_h}  ||| (\tilde{\bf p};\tilde{u}) - ({\bf r};w) |||_h \leq 
||| (\tilde{\bf p};\tilde{u}) - (\tilde{\Pi}_h \tilde{\bf p}; \tilde{\pi}_h \tilde{u}) |||_h \\
\leq \tilde{C}_I h^{k+1} \displaystyle 
\left\{ \| \tilde{\bf p} \|_{k+1+\varepsilon(1-\min[k,1]),\tilde{\Omega}} + | \nabla \cdot \tilde{\bf p} |_{k+1,\tilde{\Omega}} + | \tilde{u} |_{k+1,\tilde{\Omega}} \right\}.
\end{array}
\end{equation}
Finally, plugging \eqref{residual} and \eqref{bestapprox} into \eqref{unifymajor}, we obtain an error estimate for the non convex case with mixed boundary conditions. More precisely, we have
\begin{theorem}
\label{Estimixnonc}
Assume that $\Omega$ is not convex and $length(\Gamma_0)>0$. 
If ${\bf p} \in [H^{k+1+\varepsilon(1-\min[k,1])}(\Omega)]^2$, $u \in H^{\max[2,k+1]}(\Omega)$ and $\nabla \cdot {\bf p} \in H^{k+1}(\Omega)$, then
\begin{equation}
\label{Estimixnonconv}
\begin{array}{l}
||| ({\bf p};u) - ({\bf p}_h;u_h) |||_h^{'} \leq \max[\tilde{C}_I,\tilde{C}_{r}] \displaystyle \left\{ h^{3/2} \| \tilde{u} \|_{2,\tilde{\Omega}} \right. \\
\left.  + h^{k+1} \displaystyle \left[ \| \tilde{f} \|_{k+1+\varepsilon(1-\min[k,1]),\tilde{\Omega}} + 
 \| \tilde{\bf p} \|_{k+1+\varepsilon(1-\min[k,1]),\tilde{\Omega}} + | \nabla \cdot \tilde{\bf p} |_{k+1,\tilde{\Omega}} + | \tilde{u} |_{k+1,\tilde{\Omega}} \right] \right\}. \mbox{ \QED} 
\end{array}
\end{equation}
\end{theorem}
From \eqref{Estimixnonconv} we infer that, in case mixed conditions are prescribed on the boundary of a non convex domain, our method in original version enjoys an optimal order of convergence only for $k=0$, even if at the price of the stronger regularity assumption $u \in H^2(\Omega)$. \\

We next explain why an optimal order can also be expected of our method for $k=1$, in case $u \in H^3(\Omega)$. With this aim we first note that 
\begin{equation}
\label{residual2}
|{\mathcal R}_h(\tilde{u},{\bf q})| \leq \| \tilde{u} \|_{1/2,\Gamma_{0,h}} \| {\bf q} \cdot {\bf n}_h \|_{-1/2,\Gamma_{0,h}},
\end{equation}
where $\| \cdot \|_{-1/2,D}$ represents the norm of the dual space of $H^{1/2}(D)$, for any $D \subset \tilde{\Omega}$. According to \cite{GiraultRaviart} there is a constant $C_{\Gamma_0,h}$ depending on $\Omega_h$ such that
\begin{equation}
\label{normaltrace}
 \| {\bf q} \cdot {\bf n}_h \|_{-1/2,\Gamma_{0,h}} \leq C_{\Gamma_0,h} [| {\bf q} |]_h.
\end{equation}
On the other hand we have
\begin{equation}
\label{traceu}
\| \tilde{u} \|_{1/2,\Gamma_{0,h}} = \displaystyle \left[ \sum_{T \in {\mathcal S}_{0,h}} \| \tilde{u} \|_{1/2,e_T}^2 \right]^{1/2}.
\end{equation}
It turns out that $\tilde{u}$ vanishes at the end points of $e_T$. Therefore, by standard approximation results (cf. \cite{Arcangeli}) there is a mesh-independent constant $C_u$ such that 
\begin{equation}
\label{estimeT}
\| \tilde{u} \|_{1/2,e_T} \leq C_u l_T^2 | \tilde{u} |_{5/2,e_T}.
\end{equation}
Plugging \eqref{estimeT} into \eqref{traceu} and combining the resulting expression with \eqref{normaltrace}, after strightforward
 calculations we conclude from \eqref{residual2} that
\begin{equation}
\label{residual3}
|{\mathcal R}_h(\tilde{u},{\bf q})| \leq C_u C_{\Gamma_0,h} h^2 | \tilde{u} |_{5/2,\Gamma_{0,h}} [| {\bf q} |]_h.
\end{equation}
Now, by the Trace Theorem, there exists another constant $C_{\Omega_h}$ depending only on $\Omega_h$ such that 
\begin{equation}
\label{traceubis}
| \tilde{u} |_{5/2,\Gamma_{0,h}} \leq C_{\Omega_h} \| \tilde{u} \|_{3,\Omega_h}.
\end{equation}
Plugging \eqref{traceubis} into \eqref{residual3} we conclude that the residual ${\mathcal R}_h(\tilde{u},{\bf q})$ is likely to be an ${\mathcal O}(h^2)$-term. More precisely, there is a constant $C_{\tilde{\Omega}_h}$ such that 
\begin{equation}
\label{residual4}
|{\mathcal R}_h(\tilde{u},{\bf q})| \leq C_{\tilde{\Omega}_h} h^2 \| \tilde{u} \|_{3,\tilde{\Omega}} [| {\bf q} |]_h.
\end{equation}
Due to the closeness of $\Gamma_0$ and $\Gamma_{0,h}$, it can be strongly conjectured that $C_{\tilde{\Omega}_h}$ is mesh-independent. However we refrain from elaborating any further on this issue for the sake of brevity. Actually, in \textbf{Appendix II} we will show through a numerical example that, in the case of mixed boundary conditions too, our method is indeed of the second order for $k=1$, as long as $u$ is sufficiently smooth.\\

\noindent \underline{\textbf{IV) $\Omega$ is not convex and $length(\Gamma_0)=0$}}\\
\indent First of all we note that here, while $\tilde{\bf p}$ in \eqref{unifymajor} is the extension of ${\bf p}$ to $\tilde{\Omega}$, akin to the convex case, $\tilde{u}$ is not the aforementioned extension of $u$ to the same set, but rather a function $\bar{u}$ which differs from it by a suitable additive constant. Indeed, recalling the operator $\tilde{\pi}_h$ defined in the previous subsection, first we introduce $\bar{u} = \tilde{u} + C_h$ in $\tilde{\Omega}$, where $C_h:=- \int_{\Omega_h} \tilde{\pi}_h \tilde{u}/area(\Omega_h)$. Then we note that the estimate $\| \tilde{u} - \tilde{\pi}_h \tilde{u} \|_{0,h} \leq C_{\pi} h^{k+1} | \tilde{u} |_{k+1,\tilde{\Omega}}$ holds if $u \in H^{k+1}(\Omega)$, but not necessarily for $u \in V^k(\Omega)$. Therefore, defining in 
$\Omega_h$ $\bar{\pi}_h \tilde{u}:= \tilde{\pi}_h \tilde{u} + C_h$ we have $\bar{\pi}_h \tilde{u} \in V_h$ and thus 
\begin{equation} 
\label{errorbaru}
\left\{
\begin{array}{l}
\| \bar{u} - \bar{\pi}_h \tilde{u} \|_{0,h} \leq C_{\pi} h^{k+1} | \tilde{u} |_{k+1,\tilde{\Omega}} \\
\mbox{where}\\ 
\bar{u} = \tilde{u} + C_h \mbox{ and } \bar{\pi}_h \tilde{u}:= \tilde{\pi}_h \tilde{u} + C_h.\\
\mbox{with } C_h:=- \int_{\Omega_h} \tilde{\pi}_h \tilde{u}/area(\Omega_h).
\end{array}
\right.
\end{equation} 
\eqref{errorbaru} can be used to estimate the \textit{inf}-term in \eqref{unifymajor}, that is, 
\begin{equation} 
\label{tildepih}
\displaystyle \inf_{w \in V_h} \| \bar{u} - w \|_{0,h} \leq \| \bar{u} - \bar{\pi}_h \tilde{u} \|_{0,h} =  
\| \tilde{u} - \tilde{\pi}_h \tilde{u} \|_{0,h} \leq C_{\pi} h^{k+1} | \tilde{u} |_{k+1,\tilde{\Omega}}
\end{equation}

\begin{remark}  $\bar{u}$ can be replaced by the regular extension $\tilde{u}$ of $u$ in the \textit{sup}-term of \eqref{unifymajor} since $(C_h,\nabla \cdot {\bf q})_h = 0$ $\forall {\bf q} \in {\bf Q}_h$ if $length(\Gamma_0) =0$. \QED 
\end{remark}

Recalling the treatment of the case $length(\Gamma_0)=0$ for convex domains, an error estimate like \eqref{Estimixnonconv} holds for non convex domains without the term in $h^{3/2}$, though for a pair $({\bf p};\bar{u})$ instead of $({\bf p};u)$, where $\bar{u}$ is a function that differs from $u$ by a suitable additive constant. Indeed now sole $|(\tilde{f},v)_h - F_h(v)|$ appears in the $sup$-term of \eqref{unifymajor}. This leads to  
\begin{theorem}
\label{EstiNeumannonc}
Assume that $\Omega$ is not convex and $length(\Gamma_0)=0$. 
If ${\bf p} \in [H^{k+1+\varepsilon(1-\min[k,1])}(\Omega)]^2$, $u \in H^{k+1}(\Omega)$ and $\nabla \cdot {\bf p} \in H^{k+1}(\Omega)$, then setting $\bar{u} = u + C_h$ we have
\begin{equation}
\label{EstiNeumannonconv}
\begin{array}{l}
||| ({\bf p};\bar{u}) - ({\bf p}_h;u_h) |||_h^{'} 
\leq \max[\tilde{C}_I,C_{\chi,h}] h^{k+1} \\
\times \displaystyle \left[ \| \tilde{f} \|_{k+1+\varepsilon(1-\min[k,1]),\tilde{\Omega}} + 
 \| \tilde{\bf p} \|_{k+1+\varepsilon(1-\min[k,1]),\tilde{\Omega}} + | \nabla \cdot \tilde{\bf p} |_{k+1,\tilde{\Omega}} + | \tilde{u} |_{k+1,\tilde{\Omega}} \right]. \mbox{\QED} 
\end{array}
\end{equation}
\end{theorem}
Theorem \ref{EstiNeumannonc} indicates that, in case pure Neumann conditions are prescribed on the boundary of a non convex domain, our method enjoys an optimal order of convergence for every $k \geq 0$, provided the exact solution of \eqref{Dmix}-\eqref{definec} is sufficiently smooth.
\begin{remark}
If a convergence check is being carried out for a test-problem having a known solution $u$ sufficiently smooth, we can certainly exhibit $\tilde{u}$ in $\tilde{\Omega}_h$ for every triangulation of interest. But in this case it is simpler to define $\bar{u}$ by $u-u(Q)$ where $Q$ is a well chosen point in $\Omega_h$ at which we prescribe $u_h(Q)=0$ for all the triangulations of interest. In this manner the error $\| \bar{u} - u_h \|_{0,h}^{'}$ - or yet, more conveniently, $\| [\tilde{u}-u(Q)] - u_h \|_{0,h}$ -, is the one to be computed for every $h$ under consideration, for it will decrease at the best rate. \QED
\end{remark}

\subsection{Recovering optimal orders for mixed boundary conditions with arbitrary $k$}

\hspace{4mm} As pointed out in Remark 1 of Subsection 7.1, the error estimate \eqref{Esticonvexmix} applying to the case of mixed boundary conditions and convex domains is somewhat optimistic even for $k=0$. However in the non convex case this is no longer true if the intersection between $\bar{\Gamma}_0$ and $\bar{\Gamma}_1$ is empty. This is because in this case the exact solution of \eqref{Mixed} can be arbitrarily smooth, as long as $\Gamma$ and $f$ are sufficiently regular. As we saw in Subsection 7.1, the problem stems from the term in the variational residual involving an integral along $\Gamma_{0,h}$, which is an ${\mathcal O}(h^{3/2})$ whatever $k$. Nevertheless this issue can be settled for $k>0$ if every (straight) triangle in $T \in {\mathcal S}_{0,h}$ is replaced by a triangle $\breve{T}$ having the same straight edges as $T$ except $e_T$, which is replaced by a curvilinear edge $\breve{e}_T$ having the same end-points as $e_T$. The idea here is not to use a parametric representation, since the master element $\widehat{\breve{T}}$ associated with $\breve{T}$ is also a triangle with a curvilinear edge. Actually $\widehat{\breve{T}}$ results from the same affine mapping as the one from $T$ onto the standard master element $\hat{T}$. Therefore the integrations can be easily performed, either in $\breve{T}$ itself or in $\widehat{\breve{T}}$. The point here is that the curvilinear edge $\breve{e}_T$ is generated as the polynomial arc of degree $k+1$ which interpolates $\Gamma_0$ at $k+2$ points including the end points of $e_T$. In this manner, the residual term involving the integral of $u {\bf q} \cdot {\bf n}_h$ on $\Gamma_{0,h}$, with ${\bf q} \in {\bf Q}_h$, becomes 
\begin{equation}
\label{Gamma00}
\breve{\mathcal R}_h(\tilde{u},{\bf q}): = \displaystyle \sum_{T \in {\mathcal S}_{0,h}} (\tilde{u},{\bf q} \cdot \breve{\bf n}_T)_{\breve{e}_T},
\end{equation}
${\bf q}$ being extended to a non empty $\breve{T} \setminus T$, where $\breve{\bf n}_T$ is the outer normal vector to $\breve{e}_T$. \\
Likewise the convex case (for $k \geq 1$), we assume that $u \in H^3(\Omega)$, or yet that $\tilde{u} \in H^3(\tilde{\Omega})$.\\ 
Owing to the construction of $\breve{e}_T$, $\breve{\mathcal R}_h(\tilde{u},{\bf q})$ can be estimated as an ${\mathcal O}(h^{k+3/2})$ provided $\tilde{u} \in H^{3}(\tilde{\Omega})$.
Indeed, akin to \eqref{Gamma01}-\eqref{Gamma02}, we first establish the existence of a constant $\breve{C}_{0I}$ independent of $T$ such that
\begin{equation}
\label{Gamma08}
(\tilde{u},{\bf q} \cdot \breve{\bf n}_h)_{\breve{e}_T} \leq \breve{C}_{0I} h_T^{-1/2} \| \tilde{u} \|_{0,\breve{e}_T} \| {\bf q} \|_{0,T} \; \forall T \in {\mathcal S}_{0,h}.
\end{equation}
Therefore, recalling \eqref{Gamma00} we have
\begin{equation}
\label{residual0}
\breve{\mathcal R}_h(\tilde{u},{\bf q}) \leq  \breve{C}_{0I} \displaystyle \left[\sum_{T \in {\mathcal S}_{0,h}} h_T^{-1} \| \tilde{u} \|^2_{0,\breve{e}_T} \right]^{1/2} [| {\bf q} |]_h.
\end{equation}
Now let $P$ be the point of $\breve{e}_T$ located on the perpendicular to $e_T$ passing through a given point $M \in e_T$. Referring to Figure 1 and recalling the local coordinate system $(A_T;x,y)$ attached to $e_T$, since $\tilde{u}(N)=0$, whatever $M \in e_T$ we have  
\begin{equation}
\label{Gamma09} 
\tilde{u}(P) = \displaystyle \int_N^P \frac{\partial \tilde{u}}{\partial y} \; dy \; \forall P \in \breve{e}_T.  
\end{equation}
Let $\phi_h(x)$ be the function that uniquely represents $\breve{e}_T$ in such a frame. The curvilinear abscissa $s_h$ along $\breve{e}_T$ satisfies $ds_h = \sqrt{1+[\phi_h^{'}(x)]^2} dx$. On the other hand, by basic interpolation theory in one-dimensional space (see e.g. \cite{Blum}, \cite{Quarteroni}), there exists a constant $C_L$ independent of $l_T$ such that $\|\phi_h^{'} \|_{0,\infty,e_T} \leq C_L \|\phi^{'} \|_{0,\infty,e_T} \leq C_{\Gamma} C_L h_T$. Hence 
from \eqref{Gamma09} we easily establish the existence of a constant $\breve{C}_L$ independent of $T$ fulfilling 
\begin{equation}
\label{Gamma10}
\int_{\breve{e}_T} \tilde{u}^2(P) \; ds_h =  \displaystyle \int_{\breve{e}_T} \left[\int_N^P \frac{\partial \tilde{u}}{\partial y} \; dy \right]^2 ds_h \leq \breve{C}_L \displaystyle \int_0^{l_T} \left[\int_N^P \frac{\partial \tilde{u}}{\partial y} \; dy \right]^2 dx.
\end{equation}
Taking an obvious upper bound of the right hand side of \eqref{Gamma10}, we readily obtain $\forall T \in {\mathcal S}_{0,h}$
\begin{equation}
\label{Gamma20}
h_T^{-1} \| \tilde{u} \|_{0,\breve{e}_T}^2 \leq \breve{C}_L \max_{x \in [0,l_T]}|\phi_h(x)-\phi(x)|^2 \| \tilde{u} \|_{1,\infty,\tilde{\Omega}}^2.
\end{equation}
Resorting again to the interpolation theory in one-dimensional space, since a sufficient smoothness of $\Omega$ is assumed, we can assert that there are two constants $C_{L,k}$ and $\breve{C}_{k}$, both independent of $e_T$, such that
\begin{equation} 
\label{Boundphi}
\max_{x \in [0,l_T]}|\phi_h(x)-\phi(x)| \leq C_{L,k} l_T^{k+2} | \phi |_{k+2,\infty,e_T} \leq \breve{C}_{k} l_T^{k+2}.
\end{equation}
Plugging \eqref{Boundphi} into \eqref{Gamma20} and taking into account the continuous embedding of $H^3(\tilde{\Omega})$ into 
$W^{1,\infty}(\tilde{\Omega})$, summing up over ${\mathcal S}_{0,h}$ we readily come up with a mesh-independent constant $C_{\tilde{\Omega}}$ satisfying
\begin{equation}
\label{Gamma30}
\begin{array}{l}
\displaystyle \sum_{T \in {\mathcal S}_{0,h}} h_T^{-1} \| \tilde{u} \|_{0,\breve{e}_T}^2 \leq \breve{C}_L \breve{C}_k^2 h^{2k+3} \| \tilde{u} \|_{1,\infty,\tilde{\Omega}}^2 \displaystyle \sum_{T \in {\mathcal S}_{0,h}} l_T  
\leq C_{\tilde{\Omega}}^2 length(\Gamma_0) h^{2k+3} \| \tilde{u} \|_{3,\tilde{\Omega}}^2. 
\end{array}
\end{equation}
Finally, combining \eqref{Gamma30} and \eqref{residual0}, we infer the existence of a constant $\breve{C}_{\Gamma_0}$ depending only on $\Gamma_0$ and $\tilde{\Omega}$ such that the residual term $\breve{\mathcal R}_h(\tilde{u},{\bf q})$ fulfills
\begin{equation}
\label{brevestimgamma0}
 | \breve{\mathcal R}_h(\tilde{u},{\bf q}) | \leq \breve{C}_{\Gamma_0} h^{k+3/2} \| \tilde{u} \|_{3,\tilde{\Omega}} [| {\bf q} |]_h \; \forall {\bf q} \in {\bf Q}_h.
\end{equation} 
In short, by means of the above described boundary interpolation, in the case of mixed boundary conditions we can expect the method's order to become optimal for $k>0$ too, as long as $\Omega$ is sufficiently smooth, ${\bf p} \in [H^{k+1}(\Omega)]^2$, $u \in H^{\max[3,k+1]}(\Omega)$ and $\nabla \cdot {\bf p} \in H^{k+1}(\Omega)$. Notice however that a formal assertion in this connection requires the demonstration of the uniform stability of the thus modified formulation. We strongly conjecture that this property holds true, since we can take better advantage of material exploited in Section 5, owing to the enhanced regularity of the underlying auxiliary quantities.      
                                       
\section{Final remarks and conclusions}


\subsection{The case of pure Dirichlet conditions}

\hspace{4mm} Although it is not the purpose of our method to handle Dirichlet conditions prescribed on the whole boundary of a curved domain, a comment on this case is in order. More precisely, we mean the problem
\begin{equation}
\label{Dmix0} 
\left\{
\begin{array}{l}
\mbox{Given } f \in L^2(\Omega) \mbox{ find } ({\bf p}; u) \in {\bf H}(div,\Omega) \times L^2(\Omega) \mbox{ such that}\\ 
 -(\nabla \cdot {\bf p}, v)  = (f,v)\; \forall v \in L^2(\Omega); \\
 ({\bf p},{\bf q}) + ( u, \nabla \cdot {\bf q}) = 0 \; \forall {\bf q} \in {\bf H}(div,\Omega).
\end{array}
\right.
\end{equation} 
Mimicking \eqref{Dmixh} and observing that we do not need ${\bf P}_h$ in this case, \eqref{Dmix0} can be approximated by
\begin{equation}
\label{Dmix0h} 
\left\{
\begin{array}{l}
\mbox{Find } ({\bf p}_h; u_h) \in {\bf Q}_h \times V_h \mbox{ such that}\\ 
 -(\nabla \cdot {\bf p}_h, v)  = (f,v)\; \forall v \in V_h; \\
 ({\bf p}_h,{\bf q}) + ( u_h, \nabla \cdot {\bf q}) = 0 \; \forall {\bf q} \in {\bf Q}_h.
\end{array}
\right.
\end{equation}  
According to the predictions in Section 7, everything works fine for $k=0$ and $k=1$, as long as $u \in H^{k+2}(\Omega)$. Otherwise stated, \eqref{Dmix0h} yields optimally converging approximations for $k \leq 1$ and furthermore we only have to assume that $\Omega$ is of the $C^{k}$-class. Indeed, only an \textit{inf-sup} condition involving the product space ${\bf Q}_h \times V_h$ is necessary in this case, i.e., the classical one given in \cite{BrezziFortin}. However, if $k>1$ the variational residual  dominates the interpolation error, and hence the optimal order in the sense of ${\bf H}(div,\Omega_h^{'}) \times L^2(\Omega^{'}_h)$ can only be attained if we interpolate $\Gamma_T$ by a polynomial arc of degree $k+1$ $\forall T \in {\mathcal S}_h$ and integrate in the underlying curved element instead of $T$, as explained in Subsection 7.3.

\subsection{Miscellaneous remarks}

\subsubsection{On the assumed geometric regularity of the problem definition domain}

\hspace{4mm} The ideal geometric regularity for the error estimates given in Section 7 to hold is the assumption that $\Omega$ be of the piecewise $C^{k+2}$-class (cf. \cite{FBRT}). However, in this work we required much more, though in Lemma 3.3 only. We recall that this result was used as a basis for studying the bilinear form $d_h$ and errors for interpolates in ${\bf P}_h$. Actually we believe that this lemma, whose proof heavily relies on properties of Gaussian quadrature in an interval of $\Re$, is of academic interest only. Indeed, we obtained results of comparable accuracy with $RT_0$ by prescribing zero normal components on $\Gamma_T$ for all $T \in {\mathcal S}_{1,h}$ at 
points $P$ shifted away from the points located on the perpendicular to $e_T$ through its (Gaussian quadrature) mid-point. As matter of fact, 
the past experience on the same type of Petrov-Galerkin formulation used in this work, including results reported in \cite{ZAMM}, \cite{IMAJNA} and \cite{JCAM}, indicate that the position of the points on the true boundary at which DOFs are prescribed, does not affect accuracy very much, let alone convergence rates.

\subsubsection{Petrov-Galerkin approach vs. parametric elements} 

Some users might wonder if our Pertov-Galerkin approach with different spaces of trial- and test-fields induces higher implementation cost as compared to a pure Galerkin one with parametric elements. Such a conjecture would be due to the need for calculating a different basis for each triangle having two vertexes on $\Gamma_1$. In this respect we should report that two codes to solve Poisson  Dirichlet problems in curved domains with both approaches, using Lagrange finite elements based either on the Petrov-Galerkin formulation or on the classical isoparametric formulation were compared in \cite{Arxiv2017}. It turns out that both approaches were equivalent in terms of CPU time, with a slight advantage of the former over the later. On the other hand, the Petrov-Galerkin approach showed up significantly more accurate, which might be due to the fact that exact numerical integration can be used to compute the element matrices, in contrast to the parametric case. Finally the implementation of both approaches requires practically the same amount of data defining the curvilinear boundary. Incidentally, this is certainly true of any other method to handle accurately values prescribed on this kind of manifolds. 
   
\subsubsection{Possible extension to the three-dimensional case}

The natural extension of our method to the three-dimensional case in a context of a similar nature consists of considering the Raviart-Thomas method for tetrahedrons. In principle given normal fluxes have to be prescribed at $dim P_k(F_T)=(k+2)(k+1)/2$ points on $\Gamma_1$ for every face $F_T$ of a tetrahedron in a mesh having three vertexes on this boundary portion. However, in this case some preliminary considerations are in order: A priori the definition of the modified Raviart-Thomas space ${\bf P}_h$ using projections of Gaussian quadrature points onto the curvilinear boundary is not suitable for the three-dimensional case. Indeed, in contrast to the analogous two-dimensional situation exploited here, there is no Gauss quadrature formula with $dim P_k(F_T)$ points yielding exact integrals for polynomials of degree $\leq 2k+1$ on a triangular face $F_T$, except for $k=0$. In view of this, a counterpart of Lemma 3.3 must be discarded for arbitrary $k$. Nevertheless, as pointed out in Subsection 8.2.1, this lemma seems to be only a formal requirement, and DOFs attached to more freely chosen boundary points could be used instead, without downgrading neither accuracy nor convergence rates. But of course such a conjecture must be carefully checked, at first numerically, and eventually confirmed in a formal sense at a subsequent stage. That is why the authors intend to address the three-dimensional case in the near future.

\subsubsection{Applicability to BDM methods}

\hspace{4mm} The same formulation considered here can be applied to BDM methods in two-dimensional space introduced in \cite{BDJrM}. In this way the definition of the flux space involving triangles having a curved edge contained on the true boundary could be avoided, although a deeper look into the issue is a must. Notice that the $RT_k$ space strictly contains the space $[P_k]^2$, whose fields can be defined taking only the off-boundary DOFs of $RT_k$. This implies that the way normal fluxes are prescribed on the curved boundary play no role, at least as far as interpolation error estimates in $L^2$ are concerned. In contrast, for BDM triangles the boundary nodes too must be taken into account in order to attain the same optimality. It follows that suboptimal error estimates are expected in this case, unless known normal fluxes are prescribed on the true boundary. In the three-dimensional case (cf. \cite{BDJrDF}) the problem is even more delicate, for geometric issues related to tetrahedrons having three or even two vertexes on the curved boundary come into play. Such complications explained in \cite{IMAJNA} rule out working with tetrahedrons having three plane faces and a curved face contained in the boundary. However, we believe that our method can also be advantageously applied to this case, which we intend to examine in due course.

\subsection{Summary and conclusions}

\hspace{4mm} This work is an outstanding application of the theoretical framework built up in \cite{AsyMVF} to support the reliability analysis of numerical methods based on non symmetric mixed variational formulations with different product spaces for shape- and test-fields. This generalizes the well known theory for symmetric mixed formulations with the same product spaces.\\ 
\indent This paper is rather long, since we authors endeavored to study mixed finite element methods for PDEs posed in curved domains, by providing a rigorous treatment of several delicate issues inherent to them, while keeping in mind accuracy improvement and practical feasibility at a time. This is because, as far as we can see, the simultaneous fulfillment of these three requirements has been overlooked in most works on the subject so far. \\ 
\indent As shown in previous articles listed in the bibliography, the type of formulation adopted in this paper lines up with well established alternatives as a valid one to obtain enhanced accuracy in the solution of elliptic equations with prescribed solution values on the boundary of a curved domain. In the case of the mixed formulation of second order equations our method was compared with the formulation employing the same product spaces of shape- and test-fields, defined upon an approximating polygonal domain. For example, in the case of the $RT_1$ method, the mean-square errors of the flux variable are roughly divided by two with the use of our method, while those of the multiplier remain practically the same for both approaches, as much as their observed second order convergence rates in $H(div)$, as reported in \cite{EXCO}.

\newpage
 

\subsubsection*{APPENDIX I - Independence of $h$ of the continuity constant for Poisson problems in $\Omega_h$}

In this Appendix we establish that the continuity constant $C_{s,h}$ fulfilling \eqref{Csh} admits a strictly positive lower bound $C_s$ independent of $h$. \\
First of all we consider the case where $\Omega$ is convex.
Since in this case $\Omega_h$ is also convex for every $h$ and sufficiently close to $\Omega$, $\Gamma$ and $\Gamma_h$ can be respectively expressed in polar coordinates $(r,\theta)$ and $(r_h,\theta)$ with the same origin, say $O \in \Omega$. We know that $r_h = r R_h(\theta)/R(\theta)$ where $R(\theta)$ and $R_h(\theta)$ are the radial coordinates of $\Gamma$ and $\Gamma_h$ at the azimuthal coordinate $\theta$. We refer to Figure 2 for an illustration of both polar coordinate systems, together with some pertaining quantities exploited in the following lemmata:
\begin{figure}[h]
\label{fig2}
\centerline{\includegraphics[width=3.8in]{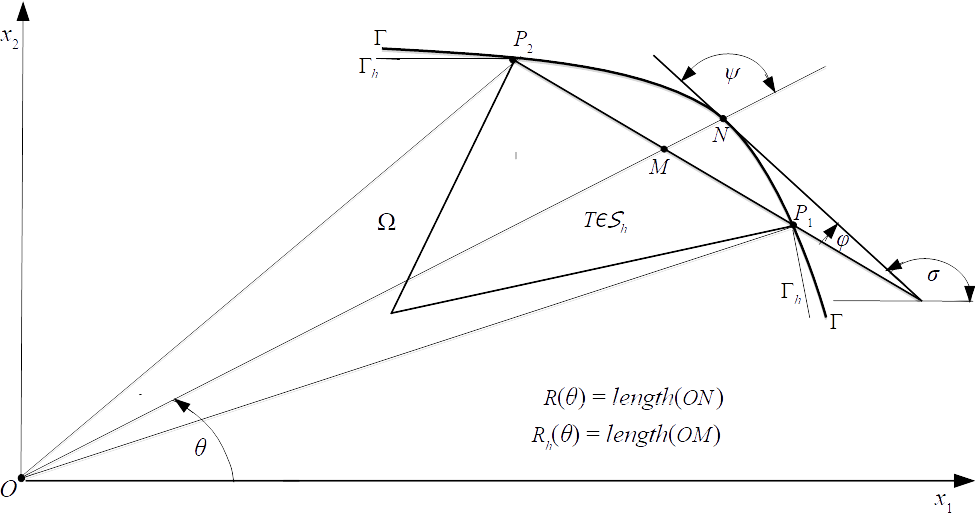}}
\vspace*{8pt}
\caption{A portion of $\Omega$ and $\Omega_h$ and the attached polar coordinate systems with origin $O \in \Omega$.}
\end{figure}  
\begin{lemma}
\label{rhohrho}
Let $\rho_h(\theta):=R_h(\theta)/R(\theta)$ for $\theta \in [0,2\pi)$. Denoting by $\tau^{,}$ the derivative with respect to $\theta$ of any function $\tau$ of $\theta$, provided $h$ is sufficiently small it holds:
\begin{enumerate}
\item There exists a constant $C^{-}$ independent of $h$ such that 
\begin{equation}
\label{rho1}
 C^{-} \leq \rho_h(\theta) \leq 1 \; \forall \theta \in [0,2\pi).
\end{equation} 
\item There exists another constant $C^{,}$ independent of $h$ such that
\begin{equation}
\label{rho2} 
\displaystyle \rho^{,}_h(\theta) \leq C^{,} h \; \forall \theta \in [0,2\pi).
\end{equation}
\end{enumerate}
\end{lemma}
\prov
As for \eqref{rho1}, we first note that $R_h(\theta) \leq R(\theta)$ for all $\theta$. On the other hand we clearly have $R_h(\theta)+ h \geq R(\theta) \; \forall \theta$. Therefore $\rho_h(\theta) \geq 1-h/R_{min}$, where $R_{min}:= \displaystyle \min_{\theta \in [0,2\pi)} R(\theta)$. Therefore under the very reasonable assumption that $h \leq R_{min}/2$ we can take $C^{-} = 1/2$.\\
Next we turn our attention to \eqref{rho2}. Referring to Figure 2, first of all we consider a generic triangle $T$ in ${\mathcal S}_h$. In Figure 2 $e_T \subset \Gamma_h$ is the edge of $T$ whose ends are vertices $P_1$ and $P_2$ of $T$. \\
According to well known properties, we may write:
\begin{equation}
\label{derirho1}
 \tan{\psi}= R(\theta)/R^{,}(\theta), 
\end{equation}
where $\psi$ is the angle between the polar radius $ON$ and the tangent to $\Gamma$ at $N$. \\
Similarly we have
\begin{equation}
\label{derirho1h}
 \tan{\psi_h} = R_h(\theta)/R_h^{,}(\theta), 
\end{equation}
where $\psi_h$ is the angle between the polar radius $OM$ and the segment $e_T$.\\
It happens that $\psi_h = \sigma - \theta$, where $\sigma$ is the angle between the cartesian axis $Ox$ and $e_T$. Moreover $\psi = \psi_h - \varphi$, where $\varphi$ is the angle between the tangent to $\Gamma$ at $N$ and the same edge of triangle $T$; skipping details for the sake of brevity, the angles $\psi$, $\sigma$ and $\varphi$ are oriented in a coherent counterclockwise sense, according to the value of $\theta$.\\

We have
\begin{equation}
\label{derirho5}
\rho_h^{,}(\theta) = \displaystyle \frac{R_h^{,}(\theta)} {R(\theta)} - \displaystyle \frac{R_h(\theta)}{R(\theta)}  
\frac{R^{,}(\theta)}{R(\theta)} = \rho_h(\theta) \displaystyle \left[ \frac{R_h^{,}(\theta)}{R_h(\theta)}-\displaystyle \frac{R^{,}(\theta)}{R(\theta)} \right].  
\end{equation}
Taking into account \eqref{derirho1} and \eqref{derirho1h}, \eqref{derirho5} yields
\begin{equation}
\label{derirho6}
\rho_h^{,}(\theta) = \rho_h(\theta) [(\tan{\psi_h})^{-1} - \tan{\psi})^{-1} = \rho_h(\theta) \{[\tan{(\sigma - \theta)}]^{-1} - [\tan{(\sigma - \varphi - \theta)}]^{-1}\}.
\end{equation}
On the other hand, it is easy to see that
\begin{equation}
\label{derirho7}
[\tan{(\sigma - \theta)}]^{-1} - [\tan{(\sigma - \varphi - \theta)}]^{-1}= \displaystyle \frac{-\sin{\varphi}}{\sin{\psi_h} \sin{\psi}}.
\end{equation}
Plugging \eqref{derirho7} into \eqref{derirho6} we obtain
\begin{equation}
\label{derirho8}
|\rho_h^{,}(\theta)| = | \sin{\varphi} | \rho_h(\theta) \sqrt{1 + \tan^{-2}{(\psi + \varphi)}} \sqrt{1 + \tan^{-2}{\psi}}.  
\end{equation}
Next we develop the first square root in \eqref{derirho8} to rewrite it as
\begin{equation}
\label{derirho9}
|\rho_h^{,}(\theta)| = |\sin{\varphi}| \rho_h(\theta) \displaystyle \sqrt{1 + \frac{(1- \tan{\psi} \tan{\varphi})^2}{(\tan{\psi} + \tan{\varphi})^2}}
\sqrt{1 + \tan^{-2}{\psi}}.
\end{equation}
Recalling \eqref{derirho1} and owing to the smallness of $\varphi$, \eqref{derirho9} yields
\begin{equation}
\label{derirho0}
\begin{array}{ll}
|\rho_h^{,}(\theta)| & \leq  |\tan{\varphi}| \displaystyle \frac{\sqrt{\tan^2{\psi} + \tan^2{\varphi}+1 + \tan^2{\psi} \tan^2{\varphi}} \sqrt{1 + \tan^{-2}{\psi}}}{|\tan{\psi} + \tan{\varphi}|} \\
& = |\tan{\varphi}| \displaystyle \frac{\sqrt{1 + \tan^2{\varphi}} (1 + \tan^{-2}{\psi})}{|1 + \tan^{-1}{\psi} \tan{\varphi}|}
\end{array}
\end{equation}
Since $\tan{\psi} = R(\theta)/R^{,}(\theta)$, setting $Q_R:= \displaystyle \max_{\theta \in [0,2\pi)} \frac{|R^{,}(\theta)|}{R(\theta)}$ 
and recalling that $|\tan{\varphi}| \leq C_{\Gamma} h$ (cf. \eqref{CGamma}), after straightforward manipulations we obtain from \eqref{derirho0} 
\begin{equation}
\label{derirho}
|\rho_h^{,}(\theta)| \leq C_{\Gamma} h \frac{(1+Q_R^2)\sqrt{1+C_{\Gamma}^2 h^2}}{1- C_{\Gamma} Q_R h}
\end{equation}
as long as $h < \displaystyle \frac{1}{C_{\Gamma} Q_R}$.
Finally, for $h \leq \displaystyle \frac{1}{2 C_{\Gamma} Q_R}$, \eqref{rho2} holds with $C^{,} =  \displaystyle \frac{1+3Q_R^2+2Q_R^4}{2 Q_R^2}$.
\QED 
\begin{lemma}
\label{lemmaA}
Let ${\mathcal A}$ be an $n$-component tensor of any order defined in $\Omega_h$ and $\tilde{\mathcal A}$ be its counterpart defined in $\Omega$ by the rule $\tilde{\mathcal A}(r,\theta)={\mathcal A}(r_h,\theta)$ $\forall (r,\theta) \in \Omega$. Assuming that both 
$\tilde{\mathcal A}$ and ${\mathcal A}$ are referred to the same cartesian frame with origin $O \in \Omega$, 
provided $h$ is not too large, for all ${\mathcal T}_h \in {\mathcal F}$ and for every $p \in [1,\infty)$ it holds
\begin{equation}
\label{L2polar}
(C^{-})^{2/p} \| \tilde{\mathcal A} \|_{0,p} \leq \| {\mathcal A} \|_{0,p,h} \leq \| \tilde{\mathcal A} \|_{0,p} \; \forall {\mathcal A} \in [L^p(\Omega_h)]^n.  
\end{equation}
\end{lemma}
\prov We have
\begin{equation}
\label{L2polarbis}
\| {\mathcal A} \|_{0,p,h}^p := \displaystyle \int_{0}^{2 \pi} \left[ \int_0^{R_h(\theta)} | {\mathcal A}(r_h,\theta) |^p r_h dr_h \right] d \theta =  \displaystyle \int_{0}^{2 \pi} [\rho_h(\theta)]^2 \left[ \int_0^{R(\theta)} | \tilde{\mathcal A}(r,\theta) |^p  r dr 
\right] d \theta.
\end{equation}
Then using \eqref{rho1} the result follows.\QED \\

Next we introduce the space $W:= \{w |\; w \in H^1(\Omega_h) \mbox{ with } w=0 \mbox{ on } \Gamma_{0,h} \mbox{ if } length(\Gamma_0)>0  \mbox{ or } \int_{\Omega_h} w=0 \mbox{ otherwise}\}$ together with 
$\tilde{W} := \{\tilde{w} |\; \tilde{w} \in H^1(\Omega) \mbox{ with } \tilde{w}=0 \mbox{ on } \Gamma_{0} \mbox{ if } length(\Gamma_0)>0  \mbox{ or } \int_{\Omega} \tilde{w}=0 \mbox{ otherwise} \}.$
Now, recalling Lemma \ref{Omegah}, given $v \in \bar{V}$, we express the Poisson problem \eqref{Poissaux} in variational form using polar coordinates, that is:
\begin{equation}
\label{Variaux}
\displaystyle \int_{0}^{2\pi} \int_{0}^{R_h(\theta)} \left[ \frac{\partial z}{\partial r_h} \frac{\partial w}{\partial r_h} + \frac{1}{r_h^2} \frac{\partial z}{\partial \theta} \frac{\partial w}{\partial \theta} \right] r_h dr_h d \theta = \displaystyle \int_{0}^{2\pi} \int_{0}^{R_h(\theta)} v w \; r_h dr_h d\theta \; \forall w \in W.
\end{equation}
Let us operate a change of variables in the partial derivatives appearing in \eqref{Variaux}.\\
We have $(\partial {\mathcal A} /\partial r_h; \partial {\mathcal A}/\partial \theta)^{T} = J_h (\partial \tilde{\mathcal A}/\partial r; \partial \tilde{\mathcal A}/\partial \theta)^T$, where $J_h$ is the Jacobian matrix of the transformation from $(r_h, \theta)$ into $(r,\omega)$ with $\omega=\theta$, expressed by
\begin{equation}
\label{Jacob} 
\begin{array}{ll}
J_h= & 
\left[
\begin{array}{ll}
\rho_h(\theta)  & \displaystyle \rho_h^{,}(\theta) r \\
0               & 1
\end{array}
\right].
\end{array}
\end{equation}
Therefore, after straightforward calculations, we come up with 
\begin{equation}
\label{Variatrans} 
\begin{array}{l}
 \displaystyle
 \int_{0}^{2\pi} \int_{0}^{R_h(\theta)} \left[ \frac{\partial z}{\partial r_h} \frac{\partial w}{\partial r_h} + \frac{1}{r_h^2} \frac{\partial z}{\partial \theta} \frac{\partial w}{\partial \theta} \right] r_h dr_h d\theta= \displaystyle \int_{0}^{2\pi} \rho_h^2 \int_{0}^{R(\theta)} \tilde{v} \tilde{w}\; r dr d\theta = \\
\displaystyle \int_{0}^{2\pi} \int_{0}^{R(\theta)} \left\{ \left[1+ (\rho_h^{,})^2 \right] \frac{\partial \tilde{z}}{\partial r} \frac{\partial \tilde{w}}{\partial r} - \frac{1}{r} \rho_h \rho_h^{,} \left[ \frac{\partial \tilde{z}}{\partial r} \frac{\partial \tilde{w}}{\partial \theta} + \frac{\partial \tilde{z}}{\partial \theta} \frac{\partial \tilde{w}}{\partial r}\right] + \frac{1}{r^2} \rho_h^2 \frac{\partial \tilde{z}}{\partial \theta} \frac{\partial \tilde{w}}{\partial \theta} \right\} r dr d\theta \; \forall \tilde{w} \in \tilde{W}.
\end{array}
\end{equation}
\begin{e-proposition}
The function $\tilde{z}$ fulfilling \eqref{Variatrans} is the unique solution of the following mixed second order elliptic equation:
\begin{equation}
\label{mixelliptic}
\left\{
\begin{array}{l}
- \nabla \cdot (\tilde{A} \nabla \tilde{z} ) = \rho_h^2 \tilde{v} \mbox{ in } \Omega \\
\tilde{z} = 0 \mbox{ on } \Gamma_0 \mbox{ if } length(\Gamma_0)>0 \mbox{ and } \int_{\Omega} \rho_h^2 \tilde{z} = 0 \mbox{ otherwise} \\
\partial (\tilde{A} \nabla \tilde{z})/\partial n = 0 \mbox{ on } \Gamma_1,
\end{array}
\right.
\end{equation} 
where, $\tilde{A}$ is a symmetric second order tensor expressed as the matrix $\tilde{A}_{|polar}$ in the frame attached to the pair of polar unit vectors $({\bf e}_r;{\bf e}_{\theta})$:
\[ 
\begin{array}{ll}
\tilde{A}_{|polar} = &  
\left[
\begin{array}{ll}
1+ (\rho_h^{,})^2  & \rho_h \rho_h^{,} \\
\rho_h \rho_h^{,}  & \rho_h^2
\end{array}
\right]
\end{array}
\]
Furthermore for a certain $s \in (2,4]$ there exits a mesh-independent constant $\tilde{C}_s$ such that 
\begin{equation}
\label{GLs}
\| \nabla \tilde{z} \|_{0,s} \leq \tilde{C}_s \| \tilde{v} \|_0. 
\end{equation}
\end{e-proposition}
\prov
First we note that, in case  $length(\Gamma_0)=0$, $\tilde{z}$ differs by an additive constant from any solution 
$\bar{z}$ of the following problem:
\begin{equation}
\label{Neumann}
\left\{
\begin{array}{l}
- \nabla \cdot [\tilde{A} \nabla \bar{z} ] = \bar{v} := \rho_h^2 \tilde{v} \mbox{ in } \Omega \\
\partial (\tilde{A} \nabla \bar{z})/\partial n = 0 \mbox{ on } \Gamma_1,
\end{array}
\right.
\end{equation} 
where both $\bar{v}$ and $\bar{z}$ belong to $L^2_0(\Omega)$. Therefore the properties of $\nabla \tilde{z}$ to be exploited below 
are the same as those of $\nabla \bar{z}$ in case $length(\Gamma_0)=0$.\\  
By a straightforward calculation we infer that the eigenvalues of $\tilde{A}$ are both strictly positive and the smallest one is bounded below by   
$\lambda_h:=\rho_h^2/[1+\rho_h^2+(\rho^{,}_h)^2]$. Thus, recalling \eqref{rho1} and \eqref{rho2}, the uniform coercivity (i.e. independently of $h$) of the bilinear form associated with \eqref{mixelliptic} is guaranteed. Indeed, 
\[ (\tilde{A} \nabla \bar{w},\nabla \tilde{w}) \geq \tilde{\lambda}  \| \nabla \tilde{w} \|_0^2 \; \forall \tilde{w} \in \tilde{W}, \] 
with $\tilde{\lambda} = (C^{-})^2/[2+(C^{,})^2] \leq \lambda_h$, since $h<1$ by assumption. Furthermore, owing to well known inequalities of the Friedrichs-Poincar\'e type (see e.g. \cite{DuvautLions}), there exists a strictly positive constant $\eta$ depending only on $\Omega$ such that $\| \nabla \tilde{w} \|_0 \geq \eta \| \tilde{w} \|_1 \; \forall \tilde{w} \in \tilde{W}$. Hence the uniform coercivity constant $\tilde{\alpha}$ of the operator $- \nabla \cdot [\tilde{A} \nabla (\cdot)]$ holds with $\tilde{\alpha}=\tilde{\lambda} \eta$, that is,
$$(\tilde{A} \nabla \tilde{w}, \nabla \tilde{w}) \geq \tilde{\alpha} \| \tilde{w} \|_1^2, \; \forall \tilde{w} \in \tilde{W}.$$
Similarly, we can easily prove that the above bilinear form is uniformly continuous, since the largest eigenvalue of $\tilde{A}$ is bounded above by $[1+\rho_h^2+(\rho^{,}_h)^2]$, or yet by $2+(C^{,})^2$, according to \eqref{rho1} and \eqref{rho2}. \\
Now resorting to \cite{TaylorKimBrown}, we infer from the above considerations on the uniform boundedness from above and below of the elliptic operator associated with either second order problem \eqref{mixelliptic} or \eqref{Neumann}, that the solution $\tilde{z}$ enjoys the property \eqref{GLs}. Indeed, the constant $\tilde{C}_s$ depends only on $\tilde{\alpha}$ (cf. \cite{TaylorKimBrown}).   \QED \\
\begin{e-proposition}
\label{Lbetah}
There exists a constant $C_s$ independent of $h$ such that the solution $z$ of \eqref{Poissaux} satisfies,
\begin{equation}
\label{GLsh}
\| \nabla z \|_{s,h} \leq C_s \| v \|_{0,h}. 
\end{equation}
\end{e-proposition}
\prov
Expressing both $z$ and $\tilde{z}$ in polar coordinates and recalling the matrix $J_h$ given by \eqref{Jacob}, first we note that $\nabla z_{|({\bf e}_r;{\bf e}_{\theta})}(r_h,\theta)$ $= J_h 
\nabla \tilde{z}_{|({\bf e}_r;{\bf e}_{\theta})}(r,\theta)$.\\
Then from \eqref{L2polar}, it follows that
\begin{equation}
\label{GLsh1}  
\| \nabla z \|_{0,s,h}^s \leq \sqrt{2} \|  |\widetilde{\nabla z}| \|_{0,s}^s = \sqrt{2} \| |J_h \nabla \tilde{z}| \|_{0,s}^s.
\end{equation}
Denoting by $\| J_h \|$ the spectral norm of $J_h$, after straightforward calculations \eqref{GLsh1} further yields,
\begin{equation}
\label{GLsh2}  
\| \nabla z \|_{0,s,h} \leq 2^{1/(2s)} \| \|J_h \| |\nabla \tilde{z}| \|_{0,s} \leq \hat{C}_s(\Omega) \| \nabla \tilde{z} \|_{0,s}, 
\end{equation}
where $\hat{C}_s(\Omega) = (2s)^{-1} \sqrt{2+ (C^{,})^2 diam(\Omega)^2}$.\\
Finally, since $\rho_h \leq 1$, combining \eqref{L2polar} with \eqref{GLs} and taking into account \eqref{GLsh2}, we obtain 
\begin{equation}
\label{GLsh3}  
\| \nabla z \|_{0,s,h} \leq \hat{C}_s(\Omega) \tilde{C}_s \| \tilde{v} \|_{0} \leq  \hat{C}(\Omega) \tilde{C}_s C^{-} \| v \|_{0,h},
\end{equation}
and the result follows with $C_s: = \hat{C}_s(\Omega) \tilde{C}_s C^{-}$. \QED \\

We next examine in main lines how Proposition \ref{Lbetah} extends to the case where $\Omega$ is not convex.\\
The key to the problem is an invertible mapping $\Upsilon_h \in [W^{1,\infty}(\Omega)]^2$ from $\Omega$ onto $\Omega_h$, together with related equations \eqref{rho3} given hereafter. \\
Actually we can consider that in the convex case such a mapping is the one based on the representation of $\Omega$ (resp. $\Omega_h$) in polar coordinates $(r;\theta)$ (resp. $(r_h;\theta)$). Here the same result can be achieved, by subdividing $\Omega$ (resp. $\Omega_h$) into a finite number $\upsilon$ of non overlapping star-shaped sub-domains $\Omega_{\iota}$ (resp. $\Omega_{h,\iota}$), in such a way that each one of them can be represented by a local system of polar coordinates $(r_{\iota};\theta_{\iota})$) (resp. $(r_{h,\iota};\theta_{\iota})$), $\iota=1,\ldots,\upsilon$. After lengthy but in all natural calculations, we can establish two analogs of \eqref{rho1} and \eqref{rho2} in the form of the following pair of estimates: 
\begin{equation}
\label{rho3}
\left\{
\begin{array}{l}
 \| \Upsilon_h - {\mathcal I} \|_{0,\infty} \leq C_{0,\Upsilon} h \\
 | \Upsilon_h |_{1,\infty} \leq C_{1,\Upsilon} h,
\end{array}
\right.
\end{equation}
where ${\mathcal I}$ is the identity operator on $\Omega$ and $C_{0,\Upsilon}$ and $C_{1,\Upsilon}$ are two constants independent of $h$.
More precisely, letting \eqref{rho3} play the role of \eqref{rho1}-\eqref{rho2} the proof of an analog of Proposition \ref{Lbetah} for the non convex case follows in a rather straightforward manner. \\
Incidentally, it is noteworthy that, if $R_{\iota}(\theta_{\iota})$ (resp. $R_{h,\iota}(\theta_{\iota})$) is the local radial coordinate of the boundary of $\Omega_{\iota}$ (resp. $\Omega_{h,\iota}$), at the intersection of the boundary of a star-shaped sub-domain, say $\Omega_{\epsilon}$, among those $\Omega$ is subdivided into, with the boundary of any other such a sub-domain, we necessarily have $R_{\epsilon}(\theta_{\epsilon})=R_{h,\epsilon}(\theta_{\epsilon})$ for all pertinent values of $\theta_{\epsilon}$.

\newpage


\subsubsection*{APPENDIX II - A numerical verification}
The authors refer to \cite{EXCO} for a thorough numerical experimentation of the method studied in this work in the case where $k \leq 1$ and $k=1$, including comparisons with the classical (i.e. polygonal) approach. Nevertheless, just to validate our analysis we report below numerical results obtained for $k=1$ taking a non convex domain. In doing so our aim here is to highlight that optimal second order convergence rates do apply to the case of a mixed Poisson problem with Dirichlet conditions prescribed on a concave boundary portion $\Gamma_0$, which confirms the predictions given at the end of Paragraph III) in Subsection 6.2. \\

More specifically, we checked the performance of the formulation \eqref{Dmixh} for $k=1$, referred to here as the P-G (for Petrov-Galerkin) $RT_1$ method, by solving a test-problem in an annulus with inner radius $r_i=1/2$ and outer radius $r_e = 1$,
 for an exact solution given by $u(x,y)= (r^2-2 r r_e +2 r_i r_e-r_i^2)/2$, where $r = \sqrt{x^2+y^2}$. This function satisfies $u=0$ on $\Gamma_0$, namely, the circle given by $r=r_i$ and $\partial u/\partial r=0$ on $\Gamma_1$, that is, along the circle given by $r=r_e$. We computed with a quasi-uniform family of meshes for a quarter annulus, constructed for a quarter unit disk with $2L^2$ triangles for $L=2^m$, by removing the $L^2/2$ triangles fully contained in the disk with radius $1/2$. For simplicity we set $h=1/L$. \\
In the upper part of Table 1 we supply the approximation errors of $u$, ${\bf p}$ and $\nabla \cdot {\bf p}$ measured in the norm of $L^2(\Omega_h)$, taking $m=2,3,4,5,6$. Although we did not study this kind of errors, in order to further illustrate the properties of the Petrov-Galerkin approximation \eqref{Dmixh}, in the lower part of Table 1 we display the maximum errors of the computed DOFs for both $u_h$ and ${\bf p}_h$, represented in the form of a mesh-dependent semi-norm $| \cdot |_{0,\infty,h}$.

\begin{table*}[h!]
{\small 
\centering
\begin{tabular}{ccccccc} &\\ [-.3cm]  
$h$ & $\longrightarrow$ & $1/4$ & $1/8$ & $1/16$ & $1/32$ & $1/64$ 
\tabularnewline &\\ [-.3cm] \hline &\\ [-.3cm]
$\parallel u_h - u \parallel_{0,h}$ & $\longrightarrow$ &              0.28440E-2 & 0.70676E-3 & 0.17641E-3 & 0.44086E-4 & 0.11021E-4 
\tabularnewline &\\ [-.3cm] \hline &\\ [-.3cm] 
$\parallel {\bf p}_h - {\bf p} \parallel_{0,h}$ & $\longrightarrow$  & 0.38543E-2 & 0.97914E-3 & 0.24598E-3 & 0.61577E-4 & 0.15400E-4
\tabularnewline &\\ [-.3cm] \hline &\\ [-.3cm] 
$\parallel \nabla \cdot({\bf p}_h - {\bf p}) \parallel_{0,h}$ & $\longrightarrow$ & 0.82685E-2 & 0.21250E-2 & 0.53575E-3 & 0.13424E-3 & 0.33579E-4
\tabularnewline &\\ [-.3cm] \hline &\\ [-.3cm] 
\tabularnewline &\\ [-.3cm] \hline &\\ [-.3cm] 
$| u_h - u |_{0,\infty,h}$ & $\longrightarrow$ &             0.12974E-1 & 0.31851E-2 & 0.79104E-3 & 0.19721E-3 & 0.49252E-4 
\tabularnewline &\\ [-.3cm] \hline &\\ [-.3cm] 
$| {\bf p}_h - {\bf p} |_{0,\infty,h}$ & $\longrightarrow$ & 0.66908E-2 & 0.17059E-2 & 0.44211E-3 & 0.11378E-3 & 0.29004E-4 
\tabularnewline &\\ [-.3cm] \hline &\\ [-.3cm] 
\end{tabular}
\caption{The P-G $RT_1$ method: errors for a test-problem in a non convex curved domain} 
}
\label{table1}
\end{table*}

The predicted second order convergence in norm of our method in all the three senses can be observed in Table 1. The same rate also applies to the computed DOFs of both $u_h$ and ${\bf p}_h$.

\end{document}